\documentclass[11pt]{amsart}
\pdfoutput=1
\usepackage[T1]{fontenc}
\usepackage{lmodern}
\usepackage[text={420pt,660pt},centering]{geometry}
\usepackage{esint,amssymb}
\usepackage{graphicx}
\usepackage{MnSymbol}
\usepackage{mathtools}
\usepackage[colorlinks=true,pdfstartview=FitV,linkcolor=blue,
 citecolor=blue,urlcolor=blue,pagebackref=false]{hyperref}
\usepackage{orcidlink}
\usepackage{microtype}

\usepackage{bm}
\usepackage{dsfont}
\usepackage{mathrsfs}
\usepackage{xcolor}
\newtheorem{proposition}{Proposition}
\newtheorem{theorem}[proposition]{Theorem}
\newtheorem{lemma}[proposition]{Lemma}
\newtheorem{corollary}[proposition]{Corollary}
\theoremstyle{remark}
\newtheorem{remark}[proposition]{Remark}
\theoremstyle{definition}
\newtheorem{definition}[proposition]{Definition}
\numberwithin{equation}{section}
\numberwithin{proposition}{section}

\newcommand{\N}{\mathbb N}

\newcommand{\R}{\mathbb R}
\newcommand{\E}{\mathbb E}
\renewcommand{\P}{\mathbb P}
\newcommand{\ep}{\varepsilon}

\renewcommand{\le}{\leqslant}
\renewcommand{\ge}{\geqslant}

\renewcommand{\subset}{\subseteq}
\renewcommand{\bar}{\overline}
\renewcommand{\tilde}{\widetilde}
\renewcommand{\hat}{\widehat}

\newcommand{\Ll}{\left}
\newcommand{\Rr}{\right}
\renewcommand{\d}{\mathrm d}
\newcommand{\dr}{\partial}
\newcommand{\1}{\mathds 1}
\newcommand{\mcl}{\mathcal}
\newcommand{\msf}{\mathsf}
\newcommand{\mfk}{\mathfrak}
\newcommand{\msc}{\mathscr}
\newcommand{\al}{\alpha}
\newcommand{\be}{\beta}
\newcommand{\ga}{\gamma}
\newcommand{\de}{\delta}
\newcommand{\si}{\sigma}
\newcommand{\sS}{\mathscr S}
\newcommand{\la}{\left\langle}
\newcommand{\ra}{\right\rangle}
\newcommand{\cH}{{L^2}}
\newcommand{\cQ}{\mathcal Q}
\newcommand{\cM}{\mathcal M}
\newcommand{\cX}{\mathfrak X}
\newcommand{\cV}{\mathcal V}
\newcommand{\cG}{\mathcal G}
\newcommand{\cF}{\mathcal F}
\newcommand{\cost}{\mathcal C}
\newcommand{\Prob}{\mathcal P}
\newcommand{\bmm}{\boldsymbol m}
\newcommand{\Id}{\mathrm{Id}}
\DeclareMathOperator{\tr}{tr}
\DeclareMathOperator{\supp}{supp}
\DeclareMathOperator{\Law}{Law}
\DeclareMathOperator{\Var}{Var}
\DeclareMathOperator{\Cov}{Cov}
\DeclareMathOperator{\diag}{diag}
\DeclareMathOperator{\Ent}{Ent}
\begin{document}

\raggedbottom

\author{Jeanne Boursier}
\address[Jeanne Boursier]{Department of Mathematics, MIT, Cambridge MA, USA}

\author{Jean-Christophe Mourrat\,\orcidlink{0000-0002-2980-725X}}
\address[Jean-Christophe Mourrat]{Department of Mathematics, ENS Lyon and CNRS, Lyon, France}

\title[Multi-species Ising spin glasses with external fields]
{Free energy of non-convex multi-species\\
Ising spin glasses with external fields}
\date{}

\begin{abstract}
We identify the limit free energy of mean-field multi-species spin glasses with $\pm1$ spins, allowing for the presence of external fields, and without assuming any convexity property of the covariance function of the Hamiltonian. The limit is given by an ``uninverted'' variational formula over Brownian martingales. The proof consists in first identifying the limit free energy of configurations with a fixed magnetization. Once the magnetization is fixed, we can show the convergence of the enriched free energy to the solution of a Hamilton--Jacobi equation. The magnetization constraint is also used to show the convexity of the initial condition in this equation, and this yields a variational representation of the solution, as a Hopf formula. Convex duality then turns this formula into an optimization over martingales whose initial values are the prescribed magnetizations. 
\end{abstract}
\maketitle
%
%
%
%
\section{Introduction}\label{s.introduction}

In this paper, we identify the limit free energy of mean-field multi-species spin glasses with $\pm1$ spins. For models with a convex covariance function, the limit free energy is given by the Parisi formula \cite{barra2015multi, chen2025free, gue03, pan, pan.multi, Tpaper}. For models with arbitrary covariance function but no external field, this limit was also identified, and takes the form of a Hopf formula~\cite{chen2026ising}. However, the variational formulas obtained in these works are known not to generalize to non-convex models with external fields \cite[Section~6]{mourrat2021nonconvex}. In our convention with a minus sign in the definition of the free energy, the Parisi formula (for convex models) and the formula in \cite{chen2026ising} (for models without external field) represent the limit free energy as a supremum. Here, we obtain instead an ``uninverted'' variational representation of the limit, which takes the form of an infimum over a space of Brownian martingales. This representation is that introduced in \cite{mourrat2025uninverting} for spin glasses with a single type of spins and extended in \cite{chen2025convex} to multi-species convex models. 

Let us define the model more precisely. Let $\sS$ be a finite set of species labels. For every $N\in\N$, let $(I_{N,s})_{s\in\sS}$ be a partition of $\{1,\ldots,N\}$, and let
\begin{equation*}
 \lambda_{N,s}:=|I_{N,s}|/N,\qquad \lambda_N := (\lambda_{N,s})_{s \in \sS}.
\end{equation*}
We assume that $\lambda_N$ converges to a vector $\lambda$ with strictly positive entries. We also define the multi-species overlap
\begin{equation*}  
 R_{N,s}(\si,\si'):=\frac1N\sum_{i\in I_{N,s}}\si_i\si'_i, \qquad R_N(\si,\si') := (R_{N,s}(\si,\si'))_{s \in \sS}.
\end{equation*}
Let $\xi:\R^{\sS}\to\R$ admit the expansion
\begin{equation}\label{e.xi.power.series}
 \xi(x)=\sum_{k=1}^{\infty}
 \sum_{s_1,\ldots,s_k\in\sS}
 \beta_{s_1,\ldots,s_k}^2
 \prod_{j=1}^k x_{s_j},
\end{equation}
which we assume to converge absolutely for every $x\in\R^{\sS}$. Let $H_N$ be a centered Gaussian field on $\{-1,1\}^N$ with covariance
\begin{equation}\label{e.def.covariance}
 \E[H_N(\si)H_N(\si')]=N\xi(R_N(\si,\si')).
\end{equation}
As discussed in \cite[Section~6]{mourrat2023free}, the stated assumption on $\xi$ guarantees the existence of such a Gaussian field, and, in the sense made precise there, functions $\xi$ of this form are essentially the only ones for which a Gaussian field exists satisfying \eqref{e.def.covariance}. We impose no convexity assumption on $\xi$.
We fix a deterministic external field $h=(h_s)_{s\in\sS}\in\R^{\sS}$, so that every spin of species $s$ feels the field $h_s$. The aim of this paper is to identify, for each $t \ge 0$, the large-$N$ limit of the free energy
\begin{equation}\label{e.def.free.energy}
 \bar F_N(t,0):=-\frac1N\E\log\sum_{\si\in\{-1,1\}^N}2^{-N}
 \exp\Bigl(\sqrt{2t}H_N(\si)-Nt\xi(\lambda_N)
 +\sum_{s\in\sS}h_s\sum_{i\in I_{N,s}}\si_i\Bigr).
\end{equation}
We will later extend $\bar F_N$ so that its second argument can be a function taking values in~$\R_+^\sS$; in \eqref{e.def.free.energy} we just take this function to be constant equal to zero. 

We next introduce the objects necessary to the description of the limit free energy. Let
\begin{equation*}
 I(z):=\frac{1+z}{2}\log(1+z)+\frac{1-z}{2}\log(1-z),
 \qquad -1\le z\le1,
\end{equation*}
with $0\log0=0$. This is the convex conjugate of $\log\cosh$. Let $T=(T_s)_{s\in\sS} \in (0,\infty)^\sS$, let $(B_s)_{s\in\sS}$ be independent standard Brownian motions, and let $\cF^s=(\cF^s_v)_{v\ge0}$ be the completed natural filtration of $B_s$. We denote by $\cM_T$ the set of processes $\al=(\al_s)_{s\in\sS}$ such that $\al_s$ is an $\cF^s$-martingale on $[0,T_s]$ with values in $[-1,1]$. In particular, $\al_s(0)=\E\al_s(T_s)$ is deterministic. For every $\al \in \cM_T$, we write
\begin{equation*}
 \cost_T(\al):=\sum_s\lambda_s\left\{
 \E\bigl[I(\al_s(T_s))-\sqrt2\al_s(T_s)B_s(T_s)\bigr]
 +\int_0^{T_s}\bigl(1-\E[\al_s(v)^2]\bigr)\,\d v\right\},
\end{equation*}
and for every $x\in\R_+^{\sS}$, we set
\begin{equation}\label{e.def.K}
 \mcl K^{T}_{t,\al}(x):=\sup_{z\in[0,1]^{\sS}}
 \left\{x\cdot(\lambda z)+t\xi(\lambda z)
 -\sum_s\lambda_s\int_0^{T_s}\bigl(z_s-\E[\al_s(v)^2]\bigr)_+\,\d v\right\},
\end{equation}
where $\lambda z:=(\lambda_sz_s)_s$ and $u_+:=\max(u,0)$. 

\begin{theorem}[Identification of the limit free energy]\label{t.main}
For every $t\ge0$, $h\in\R^{\sS}$, and $T$ satisfying $T_s\ge t\dr_s\xi(\lambda)$ for every $s\in\sS$, we have
\begin{equation}\label{e.main.martingale}
 \lim_{N\to\infty}\bar F_N(t,0)
 =\inf_{\al\in\cM_T}
 \left\{\cost_T(\al)-\sum_s\lambda_sh_s\al_s(0)
                         +\mcl K^{T}_{t,\al}(0)\right\},
\end{equation}
and the infimum is attained.
\end{theorem}
With the help of \eqref{e.explicit.conjugate}, one can check that this formula coincides with that shown in \cite[Theorem~1.5]{chen2025convex} for convex models. In Theorem~\ref{t.potential}, we also give a more abstract representation of the limit free energy that was conjectured in \cite{issa2024hopflike} and shown in \cite[Theorem~1.4]{chen2025convex} for convex models.

Our approach to the proof of Theorem~\ref{t.main} consists in first identifying the free energy of configurations with a constrained magnetization. For every integer vector $k=(k_s)_{s\in\sS}$ with $0\le k_s\le |I_{N,s}|$, we define
\begin{equation*}
 \Omega_{N,k}:=\left\{\si\in\{-1,1\}^N:
       |\{i\in I_{N,s}:\si_i=1\}|=k_s\text{ for every }s\right\},
 \qquad m_{N,k,s}:=\frac{2k_s}{|I_{N,s}|}-1,
\end{equation*}
with $m_{N,k} := (m_{N,k,s})_{s \in \sS}$. We call the set $\Omega_{N,k}$ a slice. The free energy restricted to this slice is defined by
\begin{equation}\label{e.def.slice.zero}
 \bar F_{N,k}(t,0):=-\frac1N\E\log
 \sum_{\si\in\Omega_{N,k}}2^{-N}
       \exp\bigl(\sqrt{2t}H_N(\si)-Nt\xi(\lambda_N)\bigr).
\end{equation}
The contribution of the deterministic field is constant on $\Omega_{N,k}$, and we removed it in the definition of $\bar F_{N,k}$. Since there are at most $(N+1)^{|\sS|}$ slices, the limit free energies over restricted slices determine that of the unconstrained model.

To state our results concerning these constrained models, we first introduce enriched free energies where the second argument, set at $0$ in \eqref{e.def.free.energy} and \eqref{e.def.slice.zero}, is more general. Let $\cQ$ denote the set of right-continuous increasing maps $q : [0,1) \to \R_+$, and let $\cQ_p:=\cQ\cap L^p[0,1)$ for $1\le p\le\infty$. Throughout, the word ``increasing'' is used in the sense of non-strict inequalities, that is, here we are asking that for every $u \le v \in [0,1)$, we have $q(u) \le q(v)$ (and likewise for ``decreasing'').

For every $t\ge0$, we extend $\bar F_{N,k}(t,\cdot)$ to a function on $\cQ_1^{\sS}$, the enriched free energy of the slice, and we do the same for $\bar F_N(t,\cdot)$. The path $q$ encodes the law of an additional random external field acting on the spins, with an ultrametric structure prescribed by $q$; we will comment further on this shortly after stating the next theorem.

The limit of the enriched free energy of a slice is described in terms of the following scalar quantities. For every $T \in (0,+\infty)$ and probability measure $\nu$ on $[0,T]$, let $\Phi_\nu = \Phi_\nu(v,x) : [0,T] \times \R \to \R$ solve
\begin{equation}\label{e.def.scalar.pde}
 -\dr_v\Phi_\nu=\dr_x^2\Phi_\nu+
              \nu([0,v])(\dr_x\Phi_\nu)^2,
 \qquad \Phi_\nu(T,x)=\log\cosh x,
\end{equation}
and let
\begin{equation}\label{e.def.scalar.initial}
 \psi_\circ(\nu;y):=T-\Phi_\nu(0,y),\qquad
 g_m(\nu):=\sup_{y\in\R}\{\psi_\circ(\nu;y)+ym\},
 \qquad m\in[-1,1].
\end{equation}
These quantities do not depend on $T$, as long as $[0,T]$ contains the support of $\nu$. We also evaluate $\psi_\circ$ and $g_m$ at a scalar path $q$, by evaluating them at the law of $q(U)$, where $U$ is uniform on $[0,1]$. The initial condition for the slice $m\in[-1,1]^{\sS}$ is
\begin{equation}\label{e.def.slice.initial}
 \psi_m(q):=\sum_s\lambda_s g_{m_s}(q_s),
 \qquad q\in\cQ_1^{\sS}.
\end{equation}
The extension to $\cQ_1^{\sS}$ follows from the Lipschitz estimate proved in Lemma~\ref{l.legendre.m}. For $m=0$, one can show that $g_0=\psi_\circ(\cdot;0)$, and this is the centered initial condition that showed up in \cite{chen2026ising, ho2026concavity}. The function $g_m$ is the limit of the free energy, at $t=0$, of a system with a single species constrained to have magnetization $m$; see Proposition~\ref{p.scalar.initial}.
We define $\cM_T^m:=\{\al\in\cM_T:\al_s(0)=m_s\text{ for every }s \in \sS\}$.

\begin{theorem}[Limit free energy on a slice]\label{t.slices}
For every $t\ge0$ and $q\in\cQ_2^{\sS}$, there is a continuous function $m\mapsto f_m(t,q)$ on $[-1,1]^{\sS}$ such that
\begin{equation}\label{e.slice.convergence}
\lim_{N \to \infty} \sup_k\bigl|\bar F_{N,k}(t,q)-f_{m_{N,k}}(t,q)\bigr|
        = 0.
\end{equation}
For every $m\in[-1,1]^{\sS}$, the function $f_m$ is the Lipschitz viscosity solution of
\begin{equation}\label{e.slice.hj}
 \dr_t f_m-\int_0^1\xi(\dr_q f_m(r))\,\d r=0,
 \qquad f_m(0,\cdot)=\psi_m,
\end{equation}
on the closed cone $\cQ_2^{\sS}$, in the sense specified in Section~\ref{s.hj}. It has the Hopf representation
\begin{equation}\label{e.slice.hopf}
 f_m(t,q)=\sup_{p\in\cQ_\infty^{\sS}}\inf_{q'\in\cQ_\infty^{\sS}}
 \left\{\psi_m(q')+\int_0^1\bigl(q(r)-q'(r)\bigr)\cdot p(r)\,\d r
       +t\int_0^1\xi(p(r))\,\d r\right\}.
\end{equation}
For every $q \in \mcl Q_\infty^\sS$ and every $T \in (0,\infty)^\sS$ with $T_s\ge\|q_s\|_\infty+t\dr_s\xi(\lambda)$ for every $s$, we have
\begin{equation}\label{e.slice.martingale}
 f_m(t,q)=\inf_{\al\in\cM_T^m}
 \left\{\cost_T(\al)+\int_0^1\mcl K^{T}_{t,\al}(q(r))\,\d r\right\},
\end{equation}
and the infimum is attained.
\end{theorem}

To recover the limit free energy of the full system, we will show that, for every $t \ge 0$ and $q \in \mcl Q_2^\sS$, 
\begin{equation}\label{e.full.slices}
 \lim_{N\to\infty}\bar F_N(t,q)
 =\inf_{m\in[-1,1]^{\sS}}
       \left\{f_m(t,q)-\sum_s\lambda_sh_sm_s\right\}.
\end{equation}
Combining this identity with \eqref{e.slice.martingale} releases the constraint on $\al(0)$ and proves Theorem~\ref{t.main}. It also gives its enriched version, with $\mcl K^{T}_{t,\al}(0)$ replaced by $\int_0^1\mcl K^{T}_{t,\al}(q(r))\,\d r$.

\medskip

\noindent \textbf{Proof ideas.} The core of the proof of Theorem~\ref{t.slices} consists of two facts: the enriched free energy of a slice satisfies the Hamilton--Jacobi equation \eqref{e.slice.hj} up to an error term that is small, and the initial condition $\psi_m$ is convex as a function of the path variable, so that the solution has the Hopf representation \eqref{e.slice.hopf}. 

In the usual formulation of the Hamilton--Jacobi approach \cite{chen2024ms, chen2026ising, chen2025free, chen2026spherical, HJbook, mourrat2021nonconvex, mourrat2022parisi, mourrat2023free, mourrat2020extending}, the enriched free energy is defined by coupling the spins to a Gaussian field indexed by a Poisson--Dirichlet cascade, the path $q$ prescribing the covariance between the fields at two leaves of the cascade as a function of their overlap. The Gibbs measure then features a cascade variable in addition to the spin configuration, and the error term in the PDE vanishes if and only if the overlap of two configurations sampled from the Gibbs measure is determined by the overlap of the cascade variables; see for instance \cite[(6.79)]{HJbook}. Showing that the error term is small thus requires the overlap of the configurations to synchronize with the overlap of the cascade variables. This synchronization is derived from the ultrametricity of the Gibbs measure, by the mechanism of \cite{pan.multi, pan.potts, pan.vec}, and ultrametricity itself is only known after a perturbation of the Hamiltonian, using the approach of \cite{pan.aom} based on the Ghirlanda--Guerra identities. What we obtain then is comparable to knowing that a function solves a Hamilton--Jacobi equation only almost everywhere, which does not suffice to identify it in general, and overcoming this difficulty is a major challenge.

In this paper, we use instead an equivalent definition of the enriched free energy that does not resort to any cascade and that features a single random field, which we call the branch. The branch is a stochastic process $X$ on $[0,1]$, taking values in $\R^N$, and the enriched Hamiltonian features a term $X_1 \cdot \si$ in addition to the terms already present in~\eqref{e.def.free.energy}. The law of the branch corresponds, in the setting of the cascade formulation, to the law of the Gaussian field attached to a leaf that is selected according to the Gibbs weights. This law can be described explicitly as the solution of a stochastic differential equation in $\R^N$; we refer to Section~\ref{s.branch} for the precise definition. Appendix~\ref{s.cascade} is not necessary for the purposes of showing any result in this paper, but we provide it to explain in more detail this correspondence between the classical definition of the enriched free energy and the one we adopt here. 

With this definition, the derivatives of the free energy necessarily take a different form from the classical one. They involve a family of Gibbs measures indexed by the time variable of the branch; this family corresponds to a disintegration of the usual Gibbs measure from the cascade representation according to the value of the overlap between the cascade variables. Controlling the error term in the PDE then amounts to showing that the overlap concentrates under these Gibbs measures. Adapting the approach of~\cite{aronow2026quantitative} for obtaining a quantitative proof of the Parisi formula for the SK model without external field, we show that the fluctuations of the overlap under these Gibbs measures are indeed small, with an explicit bound in terms of $N$ and of a lower bound on the derivative of the path (Theorem~\ref{t.remainder}). One key observation is that conditionally on a sampled configuration, the Gaussian data have a log-concave law, and one can therefore use the Brascamp--Lieb inequality to control the fluctuations of functionals of the disorder. No perturbation of the Hamiltonian is involved. This argument bounds the fluctuations conditionally on the sampled configuration. By working on a slice, we ensure that intra-species permutations of the coordinates act transitively on the model, and thus for observables that are invariant under this transformation, the conditioning can be removed. In the model without external field studied in \cite{aronow2026quantitative}, the spin-flip symmetry is used instead to cover the full configuration space, but the presence of an external field breaks the invariance under such transformations. 

The second role played by the restriction of the free energy to slices concerns the initial condition. The Hopf representation of the solution to a Hamilton--Jacobi equation requires a convex initial condition, and the initial condition $\psi_\circ(\cdot;h)$ of the model with a fixed field $h$ is not convex in the path when $h$ is large, as was pointed out in \cite[Section~6]{mourrat2021nonconvex}. Fixing the magnetization changes this: we prove that~$g_m$ is convex as a function of the path variable. That $\psi_\circ(\cdot;0)=g_0$ satisfies this property has already been shown in \cite{chen2026ising, ho2026concavity}, and here we extend this property to every $m$. This yields a variational representation for the limit free energy over each slice, which we turn into the uninverted formula \eqref{e.slice.martingale} using convex-duality arguments, and the full free energy is finally recovered by minimizing over the slices, as in \eqref{e.full.slices}. 

To close this introduction, we also mention that the existence of the limit free energy of a broader class of spin-glass models was very recently shown in \cite{ho2026convergence}, also using ideas from~\cite{aronow2026quantitative}. Now that the existence of the limit free energy is established for this broader class of models, we can appeal to \cite{chen2024ms, chen2025free} to assert that this limit can be represented in terms of a (possibly non-unique) critical point of an explicit functional.

\medskip

\noindent \textbf{Organization of the paper.} Sections~\ref{s.branch} to~\ref{s.derivatives} define the enriched free energy and derive the identities for its derivatives, as outlined above, both for the full system and for the one restricted to a slice. This yields that the error term in the Hamilton--Jacobi equation is controlled by an average of the variances of the overlap under the family of Gibbs measures discussed above. Section~\ref{s.concentration} shows that these variances are small on a slice, uniformly over the slices, once the slopes of the path are bounded below, following the method of~\cite{aronow2026quantitative}. Section~\ref{s.scalar} presents several properties of the scalar function~$g_m$, including that it can be represented as an infimum over martingales started at $m$, and that it is a convex function of the path variable. Section~\ref{s.initial} identifies the limit of the initial condition $\bar F_{N,k}(0,\cdot)$ as $\psi_{m_{N,k}}$, uniformly over the slices. Section~\ref{s.hj} uses this to show \eqref{e.slice.convergence}--\eqref{e.slice.hopf} of Theorem~\ref{t.slices}. Section~\ref{s.duality} turns the Hopf formula into the martingale formula via a convex-duality argument, yielding \eqref{e.slice.martingale}. Section~\ref{s.conclusion} shows \eqref{e.full.slices} and completes the proof of Theorem~\ref{t.main}. Finally, Appendix~\ref{s.cascade} relates the enriched free energy used here to its classical definition through Poisson--Dirichlet cascades; the main text does not depend on the results of this appendix.

\medskip

\noindent \textbf{Statement of AI use.} The proofs of our main results were first obtained by prompting ChatGPT 6 Astra. The initial arguments thus obtained were however difficult for us to parse, and our main contribution has been to turn these arguments into a form that we find conceptually appealing. We were assisted by ChatGPT 6 Astra and Claude Fable~5.1 in this process. We have manually checked all statements and proofs and polished the paper, and we take full responsibility for its contents.

\medskip

\noindent \textbf{Acknowledgements.} JCM acknowledges the support of the MSCA grant SLOHD (101203974).

%
%
%
%
\section{The free energy on one branch}
\label{s.branch}

In this section, we define the enriched free energy of a slice, and more generally of any weighted set of configurations. Throughout this section and the next three, $N$ is fixed, and every statement holds for every finite set $\sS$, every $N\ge1$, and every partition $(I_{N,s})_{s\in\sS}$ of $\{1,\ldots,N\}$ into nonempty sets.

\subsection{Setting}
\label{s.branch.setting}

Let $\Omega$ be a nonempty subset of $\Sigma_N:=\{-1,1\}^N$, and let $w:\Omega\to(0,\infty)$ be a weight function. Two cases are of primary interest. In the first, $\Omega$ is a slice $\Omega_{N,k}$ and $w\equiv2^{-N}$, so that the slice free energy retains the entropy cost of the constraint, as in \eqref{e.def.slice.zero}. In the second, $\Omega$ is the full cube $\Sigma_N$ and $w(\si)=2^{-N}e^{h^{(N)}\cdot\si}$ is the weight in \eqref{e.def.free.energy}, where $h^{(N)}\in\R^N$ denotes the vector with $h^{(N)}_i:=h_s$ for $i\in I_{N,s}$, so that $h^{(N)}\cdot\si=\sum_sh_s\sum_{i\in I_{N,s}}\si_i$. A third case, a coupling of two slices, appears in Section~\ref{s.initial}. Since the configuration set is finite, the positive semidefinite kernel $(\si,\si')\mapsto N\xi(R_N(\si,\si'))$ on $\Sigma_N$ admits a Gram representation: there exist a finite-dimensional Euclidean space and vectors $(v_\si)_{\si\in\Sigma_N}$ in it such that
\begin{equation}\label{e.gram}
 H_N(\si)=g\cdot v_\si,\qquad
 v_\si\cdot v_{\si'}=N\xi(R_N(\si,\si')),\qquad
 |v_\si|^2=N\xi(\lambda_N),
\end{equation}
where $g$ is a standard Gaussian vector in this space. We fix one such representation on the full cube and use it for every $\Omega$, so that the disorder is one finite Gaussian vector~$g$, shared by all the slices. We write $P_s$ for the orthogonal projection of $\R^N$ onto the coordinates in $I_{N,s}$, and, for a vector $a=(a_s)_{s\in\sS}\in\R^{\sS}$, we write $\diag(a)$ for the diagonal $N\times N$ matrix whose $i$-th entry is $a_s$ for $i\in I_{N,s}$. For every realization of $g$, every $t\ge0$ and every $x\in\R^N$, we set
\begin{equation}\label{e.def.terminal}
 U(1,x;g):=\log\sum_{\si\in\Omega}w(\si)
                  \exp\bigl(\sqrt{2t}\,g\cdot v_\si+x\cdot\si\bigr).
\end{equation}
For every fixed $g$, the function $U(1,\cdot;g)$ is smooth and convex, its gradient is the mean of $\si$ under the Gibbs measure on $\Omega$ with weights $w$ and Hamiltonian $\sqrt{2t}\,g\cdot v_\si+x\cdot\si$, its Hessian is the covariance of $\si$ under this measure, and all its derivatives are bounded, being cumulants of the bounded vector $\si$. In particular $|\dr_{x_i}U(1,\cdot;g)|\le1$ for every $i$.

\subsection{The branch equation}
\label{s.branch.equation}

We say that a path $q:[0,1]\to\R_+^{\sS}$ is \emph{regular} if it is continuous and piecewise affine with finitely many pieces, and if
\begin{equation}\label{e.def.regular}
 q_s(0)>0
 \qquad\text{and}\qquad
 \dot q_s>0\ \text{on every piece},
\end{equation}
for every species $s$; restricted to $[0,1)$, such a path belongs to $\cQ_\infty^{\sS}$. For a regular path, let
\begin{equation}\label{e.def.D}
 D(r):=\diag(\dot q(r)),
\end{equation}
so that $D$ is piecewise constant and positive definite, and $\diag(q(0))$ is positive definite. The free energy is first defined for regular paths. Its extension to all of $\cQ_1^{\sS}$ is obtained in Proposition~\ref{p.extension} below, after the Lipschitz estimate that justifies it. For a fixed realization of $g$, we consider the backward equation
\begin{equation}\label{e.branch.pde}
 \dr_rU+\tr\bigl(D(r)\nabla^2U\bigr)+r\,\nabla U\cdot D(r)\nabla U=0
 \quad\text{on }[0,1)\times\R^N,
\end{equation}
with the terminal condition \eqref{e.def.terminal}. All spatial gradients below are with respect to $x\in\R^N$, and we write
\begin{equation}\label{e.def.branch.operator}
 \mcl L[W]:=\dr_rW+\tr\bigl(D(r)\nabla^2W\bigr)+r\,\nabla W\cdot D(r)\nabla W
\end{equation}
for the nonlinear operator in \eqref{e.branch.pde}. This equation is the finite-$N$ analogue of the Parisi equation \eqref{e.def.scalar.pde}, up to a time change discussed further in Section~\ref{s.scalar.dictionary} for a single spin.

\begin{lemma}[Well-posedness of the branch equation]
\label{l.wellposed}
Let $q$ be regular, and fix a realization of $g$. There is a unique function $U=U(\cdot,\cdot\,;g)$ on $[0,1]\times\R^N$, continuous and with $|\dr_{x_i}U|\le1$ for every $i$, which is a classical solution of \eqref{e.branch.pde} on each piece of $q$ and satisfies the terminal condition \eqref{e.def.terminal}. Its spatial derivatives of all orders are continuous and bounded on $[0,1]\times\R^N$, and $\dr_rU$ is bounded and continuous away from the breakpoints of $q$.
\end{lemma}

\begin{proof}
We argue by backward induction on the pieces of $q$. On a piece $[s_j,s_{j+1}]$, the matrix $D$ is constant and positive definite, so that \eqref{e.branch.pde} is a uniformly parabolic semilinear equation with coefficients that are smooth in $r$ and a nonlinearity that is quadratic in the gradient. The terminal datum at $s_{j+1}$ is \eqref{e.def.terminal} for the last piece, and otherwise the value already constructed at $s_{j+1}$; in both cases it has bounded derivatives of all orders. The standard theory of such equations, see for instance \cite[Chapter~V]{LSU}, provides a unique classical solution with bounded derivatives of all orders on the piece, the gradient bound below supplying the a priori estimate needed for global existence; the solution converges to its terminal datum in every $C^k$ norm as $r\uparrow s_{j+1}$. Hence the spatial derivatives of the glued function are continuous across $s_{j+1}$. Each partial derivative $w:=\dr_{x_i}U$ solves the linear equation $\dr_rw+\tr(D\nabla^2w)+2r\,\nabla U\cdot D\nabla w=0$, with terminal datum in $[-1,1]$, so the maximum principle gives $|w|\le1$. Uniqueness among continuous solutions with bounded gradient is the case $\phi_1=\phi_2=0$ of Lemma~\ref{l.comparison} below.
\end{proof}

\begin{remark}[Other terminal data]\label{r.terminal.data}
The proof uses only two properties of the terminal datum \eqref{e.def.terminal}: its derivatives of all orders are bounded, and its first derivatives take values in $[-1,1]$. Lemma~\ref{l.wellposed} therefore holds for every terminal datum with these two properties. 
\end{remark}

The following comparison principle is used repeatedly: it gives the uniqueness in Lemma~\ref{l.wellposed}, the stability of the solution with respect to the path (Lemmas~\ref{l.zero} and~\ref{l.block}), and the comparison between the free energies of the full model and of the slices (Lemma~\ref{l.gluing}).

\begin{lemma}[Comparison principle]
\label{l.comparison}
Let $q$ be regular. Let $W_1$ and $W_2$ be continuous functions on $[0,1]\times\R^N$ with bounded spatial gradients, which on each piece of $q$ are classical solutions of $\mcl L[W_i]=\phi_i$, where $\phi_i$ is bounded on $[0,1)\times\R^N$ and continuous on each piece. For every $r\in[0,1]$, we have
\begin{equation}\label{e.comparison}
 \sup_{x\in\R^N}|W_1(r,x)-W_2(r,x)|
 \le\sup_{x\in\R^N}|W_1(1,x)-W_2(1,x)|
  +\int_r^1\sup_{x\in\R^N}|\phi_1(u,x)-\phi_2(u,x)|\,\d u .
\end{equation}
If moreover $\phi_1\le\phi_2$ on $[0,1)\times\R^N$ and $W_1(1,\cdot)\ge W_2(1,\cdot)$, then $W_1\ge W_2$ on $[0,1]\times\R^N$.
\end{lemma}

\begin{proof}
The difference $w:=W_1-W_2$ satisfies, on each piece of $q$,
\begin{equation*}
 \dr_rw+\tr(D\nabla^2w)+r\,(\nabla W_1+\nabla W_2)\cdot D\nabla w=\phi_1-\phi_2 ,
\end{equation*}
which is a linear equation with the bounded drift coefficient $rD(\nabla W_1+\nabla W_2)$. Fix $(r,x)$. By Girsanov's theorem, there is a weak solution $Y$ on $[r,1]$ of $\d Y_u=uD(u)(\nabla W_1+\nabla W_2)(u,Y_u)\,\d u+\sqrt{2D(u)}\,\d B_u$ with $Y_r=x$. Since $w$ grows at most linearly and $\nabla w$ is bounded, It\^o's formula on each piece and the continuity of $w$ across the breakpoints give
\begin{equation*}
 w(r,x)=\E\bigl[w(1,Y_1)\bigr]-\E\int_r^1(\phi_1-\phi_2)(u,Y_u)\,\d u ,
\end{equation*}
the stochastic integral having zero mean. Both assertions follow.
\end{proof}

The arguments used in the proof of Lemma~\ref{l.wellposed} also give smooth dependence on parameters: if the terminal condition or the coefficients $D$ and $\diag(q(0))$ depend smoothly on a real parameter, with bounded derivatives, then so does $U$, and each parameter derivative solves a linear parabolic equation with bounded coefficients, obtained by differentiating \eqref{e.branch.pde}. 
This is used for the response identity of Section~\ref{s.gibbs} and for the derivatives of Section~\ref{s.derivatives}.

\subsection{The branch and the definition}
\label{s.branch.definition}

We fix $t\ge0$ and a regular path $q$. We sample~$g$ and $X_0\sim\mcl N(0,2\diag(q(0)))$ independently, and we let $B$ be an independent standard $N$-dimensional Brownian motion. Since $\nabla^2U$ is bounded, the stochastic differential equation
\begin{equation}\label{e.branch.sde}
 \d X_r=2rD(r)\nabla U(r,X_r;g)\,\d r+\sqrt{2D(r)}\,\d B_r,
 \qquad 0\le r\le1,
\end{equation}
has a unique strong solution. We call $X$ the \emph{branch}, and set
\begin{equation}\label{e.def.m.A}
 \bmm_r:=\nabla U(r,X_r;g)\in[-1,1]^N,
 \qquad
 A_q:=\int_0^1r\,\bmm_r\cdot D(r)\bmm_r\,\d r.
\end{equation}
Throughout, $\E$ denotes the expectation with respect to all sources of randomness: so far this means the disorder $g$, the initial field $X_0$, and the Brownian motion $B$.

\begin{definition}[Enriched free energy on one branch]
\label{d.free.energy}
For $t\ge0$ and a regular path $q$, let
\begin{equation*}
 H^{\mathrm{br}}_N(\si):=\sqrt{2t}H_N(\si)+X_1\cdot\si
       -Nt\xi(\lambda_N)-N\lambda_N\cdot q(1)-A_q,
\end{equation*}
and define
\begin{equation}\label{e.def.branch.free.energy}
 \bar F_{N,\Omega}(t,q):=-\frac1N\E\log
          \sum_{\si\in\Omega}w(\si)e^{H^{\mathrm{br}}_N(\si)},
\end{equation}
the dependence on the choice of the weight $w$ being kept implicit. 
We write $\bar F_{N,k}:=\bar F_{N,\Omega_{N,k}}$ for a slice with the weights $2^{-N}$, and $\bar F_N:=\bar F_{N,\Sigma_N}$ for the full cube with the weights $2^{-N}e^{h^{(N)}\cdot\si}$.
\end{definition}

The three subtracted terms in $H^{\mathrm{br}}_N$ do not depend on the configuration $\si$ that is summed over, and will therefore not have any influence on the associated Gibbs measure. The first two are the usual self-overlap corrections: $t\xi(\lambda_N)$ is half the variance of $\sqrt {2t}H_N(\si)/\sqrt N$, and $\lambda_N\cdot q(1)$ is half the variance per site of the Gaussian part of the field. The term $A_q$ records a relative-entropy cost, see Remark~\ref{r.entropy}. The definition  \eqref{e.def.branch.free.energy} can be rewritten as
\begin{equation}\label{e.free.energy.unpacked}
 \bar F_{N,\Omega}(t,q)=t\xi(\lambda_N)+\lambda_N\cdot q(1)
                -\frac1N\E\bigl[U(1,X_1;g)-A_q\bigr].
\end{equation}

\begin{lemma}[It\^o's formula along the branch]\label{l.ito}
For every $g$, $X_0$, and almost every~$B$, we have
\begin{equation}\label{e.ito.identity}
 U(1,X_1;g)-A_q=U(0,X_0;g)+\int_0^1\bmm_r\cdot\sqrt{2D(r)}\,\d B_r.
\end{equation}
Consequently
\begin{equation}\label{e.free.energy.V}
 \bar F_{N,\Omega}(t,q)=t\xi(\lambda_N)+\lambda_N\cdot q(1)-\frac1N\E\,U(0,X_0;g).
\end{equation}
\end{lemma}

\begin{proof}
By Lemma~\ref{l.wellposed}, It\^o's formula applies on each piece of $q$, and gives
\begin{align*}
 \d U(r,X_r;g)
 &=\bigl(\dr_rU+\tr(D\nabla^2U)+2r\,\nabla U\cdot D\nabla U\bigr)(r,X_r;g)\,\d r
   +\bmm_r\cdot\sqrt{2D(r)}\,\d B_r\\
 &=r\,\bmm_r\cdot D(r)\bmm_r\,\d r+\bmm_r\cdot\sqrt{2D(r)}\,\d B_r,
\end{align*}
where we used \eqref{e.branch.pde}. Integration over $[0,1]$ gives \eqref{e.ito.identity}. The stochastic integral has zero mean since $|\bmm_r|\le\sqrt N$, and \eqref{e.free.energy.unpacked} gives \eqref{e.free.energy.V}.
\end{proof}

\begin{remark}[Control characterization]\label{r.control}
The same calculation, applied to the controlled process $Y^v_r:=X_0+\int_0^r2uD(u)v_u\,\d u+\int_0^r\sqrt{2D(u)}\,\d B_u$ for a progressively measurable $v$ with values in $[-1,1]^N$, gives
\begin{multline*}
 U(0,X_0;g)-\E_B\Ll[U(1,Y_1^v;g)-\int_0^1r\,v_r\cdot D(r)v_r\,\d r\Rr]\\
 =\E_B\int_0^1r\,\bigl|v_r-\nabla U(r,Y_r^v;g)\bigr|_{D(r)}^2\,\d r\ge0,
\end{multline*}
where $|w|_D^2:=w\cdot Dw$, with equality exactly for the branch. Thus $U(0,X_0;g)$ is the value of a stochastic control problem, the branch is its optimal state, and $\bar F_{N,\Omega}(t,q)$ could equivalently be defined by maximizing the expected terminal log partition function minus the quadratic running cost, separately for each realization of $(g,X_0)$. This is the finite-$N$ counterpart of the control representation of the Parisi functional \cite{auffinger2015parisi, jagtob}.
\end{remark}
We refer to Appendix~\ref{s.cascade} for the relationship between the enriched free energy as defined here and that based on Poisson--Dirichlet cascades. 

\subsection{Symmetries of the branch equation}
\label{s.branch.symmetries}

Three elementary symmetries of \eqref{e.branch.pde} will be used: translations of the field and constant shifts, the separation of the equation over blocks of coordinates, and the permutations of the sites within a species. For the last one, let $S$ be a permutation of $\{1,\ldots,N\}$ that preserves each $I_{N,s}$, acting on $\R^N$ by permutation of the coordinates. It preserves every species overlap, so the Gram matrix in \eqref{e.gram} is invariant under $S$, and there is an orthogonal map $O_S$ of the disorder space with $O_Sv_\si=v_{S\si}$ for every $\si\in\Sigma_N$. Such a permutation preserves the slices and the full cube.

\begin{lemma}[Symmetries]\label{l.symmetries}
Let $q$ be regular, and fix a realization of $g$.
\begin{enumerate}
\item For every $y\in\R^N$ and $c\in\R$, the solution of \eqref{e.branch.pde} with the terminal datum $x\mapsto U(1,x+y;g)+c$ is $(r,x)\mapsto U(r,x+y;g)+c$.
\item Let $\{1,\ldots,N\}=B_1\cup\cdots\cup B_p$ be a partition into blocks, each contained in a single species, and suppose that the terminal datum is of the form $U(1,x;g)=\sum_\ell U_\ell(1,x_{B_\ell})$, where $x_{B_\ell}$ denotes the restriction of $x$ to the coordinates in $B_\ell$ and each $U_\ell(1,\cdot)$ has bounded derivatives of all orders. Then $U(r,x;g)=\sum_\ell U_\ell(r,x_{B_\ell})$, where $U_\ell$ solves \eqref{e.branch.pde} in dimension $|B_\ell|$, with $D$ replaced by its restriction to the coordinates in $B_\ell$, and with the terminal datum $U_\ell(1,\cdot)$.
\item Suppose that $S\Omega=\Omega$ and $w\circ S=w$. Then, for every such permutation $S$,
\begin{equation}\label{e.permutation.symmetry}
\begin{gathered}
 U(r,Sx;O_Sg)=U(r,x;g),\qquad
 \nabla U(r,Sx;O_Sg)=S\nabla U(r,x;g),\\
 \nabla^2U(r,Sx;O_Sg)=S\nabla^2U(r,x;g)S^{\mathsf T}.
\end{gathered}
\end{equation}
\end{enumerate}
\end{lemma}

\begin{proof}
In each case, the function on the right-hand side is continuous with bounded gradient, so by the uniqueness in Lemma~\ref{l.wellposed} it suffices to check that it solves \eqref{e.branch.pde} with the stated terminal datum. For (1), the equation only involves derivatives of $U$ and is invariant under translations in $x$. For (2), since $D$ is diagonal and constant on each species, the function $\sum_\ell U_\ell(r,x_{B_\ell})$ satisfies $\tr(D\nabla^2\sum_\ell U_\ell)=\sum_\ell\tr(D_{B_\ell}\nabla^2U_\ell)$, and the quadratic term separates as well, because the gradients of the different blocks are supported on disjoint coordinates. For (3), let $W(r,x):=U(r,Sx;O_Sg)$. Then $\nabla W(r,x)=S^{\mathsf T}\nabla U(r,Sx;O_Sg)$ and $\nabla^2W(r,x)=S^{\mathsf T}\nabla^2U(r,Sx;O_Sg)S$. Since $D$ is constant on each species, $SDS^{\mathsf T}=D$, so that $\tr(D\nabla^2W)=\tr(D\nabla^2U)$ and $\nabla W\cdot D\nabla W=\nabla U\cdot D\nabla U$, evaluated at $(r,Sx;O_Sg)$: the function $W$ solves \eqref{e.branch.pde}. As for the terminal datum, since $O_S$ is orthogonal, $O_Sg\cdot v_\si=g\cdot O_S^{-1}v_\si=g\cdot v_{S^{-1}\si}$, and $Sx\cdot\si=x\cdot S^{-1}\si$. Substituting $\si'=S^{-1}\si$ in \eqref{e.def.terminal}, which is a bijection of $\Omega$ preserving the weights, gives $U(1,Sx;O_Sg)=U(1,x;g)$. The identities for the derivatives follow by differentiating the first identity in \eqref{e.permutation.symmetry}.
\end{proof}
%
%
%
%
\section{The Gibbs measures along the branch}\label{s.gibbs}

The quantities that drive the argument are Gibbs averages of the sampled configuration, conditioned on the branch up to an intermediate level. We introduce them now, together with the response identity that expresses the derivatives of $U$ along the branch in terms of these averages. Throughout this section, $(\Omega,w)$ is as in Section~\ref{s.branch}, $t\ge0$ and $q$ is a regular path.

\subsection{The sampled configuration and the Gibbs measure at level \texorpdfstring{$r$}{r}}
\label{s.gibbs.terminal}

After sampling~$g$ and the branch~$X$, we sample one configuration $\si$ from the normalized Gibbs measure of $H^{\mathrm{br}}_N$, that is,
\begin{equation}\label{e.sampled.configuration}
 \P\bigl(\si=\tau\bigm|g,(X_r)_{0\le r\le1}\bigr)
 =w(\tau)\exp\bigl(\sqrt{2t}\,g\cdot v_\tau+X_1\cdot\tau-U(1,X_1;g)\bigr),
 \qquad\tau\in\Omega.
\end{equation}
We call $\si$ the \emph{sampled configuration}. Let $\cG_r$ be the $\sigma$-algebra generated by $g$, $X_0$ and $(B_u)_{u\le r}$; it does not contain the sampled configuration, even for $r=1$. For $r\in[0,1]$ and a function $b$ of $n$ configurations $\si^1,\ldots,\si^n$, we write $\langle b\rangle_r$ for the integral of $b$ against the $n$-fold product of the conditional law of $\si$ given $\cG_r$. Thus the replicas at level $r$ are conditionally independent given $\cG_r$, each with the conditional law of $\si$ given $\cG_r$, and $\langle b\rangle_r=\E[b(\si)\mid\cG_r]$ when $n=1$. For a function $b$ of one configuration, we also write
\begin{equation*}
 \Var_r(b):=\langle b^2\rangle_r-\langle b\rangle_r^2
\end{equation*}
for its conditional variance given $\cG_r$. The tower property gives, for $r\le u$,
\begin{equation}\label{e.tower}
 \langle b\rangle_r=\E\bigl[\langle b\rangle_u\bigm|\cG_r\bigr],
\end{equation}
so that $r\mapsto\langle b\rangle_r$ is a martingale. At $r=1$, $\langle\cdot\rangle_1$ is the Gibbs measure in \eqref{e.sampled.configuration}.

We note that the conditional laws of $\si$ given $\cG_r$ thus form a measure-valued martingale that evokes the stochastic localization scheme of \cite{chen2022localization, eldan2013thin}. However, the presence of the coefficient $r$ in \eqref{e.branch.pde} and \eqref{e.branch.sde} in place of $1$ makes this interpretation difficult to flesh out precisely. 

We also point out that, for $r<1$, the average $\langle\cdot\rangle_r$ integrates over the unobserved continuation $(X_u)_{u>r}$ of the branch, and it is in general not the Gibbs measure obtained by replacing $X_1$ with $X_r$ in the terminal Hamiltonian. To sample two replicas at level~$r$, one uses two independent continuations of the branch after level $r$, each followed by an independent draw from the Gibbs measure \eqref{e.sampled.configuration} in its terminal field. In the cascade picture of Appendix~\ref{s.cascade}, two replicas at level $r$ are two replicas of the cascade Gibbs measure, conditioned on the event that the cascade variables separate at level $r$; see Lemma~\ref{l.selected.pair} and Proposition~\ref{p.cascade.pair}.

Lemma~\ref{l.response} below shows that the Gibbs mean $\langle\si\rangle_r$ is the drift vector $\bmm_r=\nabla U(r,X_r;g)$ of \eqref{e.def.m.A}. Anticipating this, we call $\bmm_r$ the posterior mean, and we define
\begin{equation}\label{e.def.posterior.quantities}
 C_r:=\Cov(\si\mid\cG_r)=\bigl\langle(\si-\bmm_r)(\si-\bmm_r)^{\mathsf T}\bigr\rangle_r,
 \qquad
 K_r:=\nabla^2U(r,X_r;g)
\end{equation}
for the conditional covariance of the sampled configuration and for the Hessian of $U$ along the branch, which we call the response matrix. We write $C_{r,ss}:=P_sC_rP_s$ and $K_{r,ss}:=P_sK_rP_s$ for their diagonal blocks. All matrix comparisons below are in the positive semidefinite order $\preceq$. Since $\si_i^2=1$, we have
\begin{equation}\label{e.trace.covariance}
 \tr C_{r,ss}=|I_{N,s}|-|P_s\bmm_r|^2,
 \qquad
 \tr C_r=N-|\bmm_r|^2\le N .
\end{equation}
Finally, for two replicas $\si^1,\si^2$, we write $R_N:=R_N(\si^1,\si^2)$ and $R_{N,s}:=R_{N,s}(\si^1,\si^2)$. We set
\begin{equation}\label{e.def.S.p}
 S_{r,s}:=\langle R_{N,s}\rangle_r=\frac{|P_s\bmm_r|^2}{N},
 \qquad
 p_{N,s}(r):=\E\langle R_{N,s}\rangle_r=\E S_{r,s},
\end{equation}
and $S_r$, $p_N(r)$ denote the corresponding vectors in $\R^{\sS}$. The quantity $S_{r,s}$ is the random Gibbs mean of the overlap at level $r$, and $p_{N,s}(r)$ is its expectation; both lie in $[0,\lambda_{N,s}]$. 

\subsection{The response identity}
\label{s.gibbs.response}

For a function $b:\Omega\to\R$, possibly depending on $g$, and $y\in\R$, let $U^y$ denote the solution of \eqref{e.branch.pde} with the terminal condition obtained from \eqref{e.def.terminal} by adding $yb(\si)$ to the Hamiltonian,
\begin{equation*}
 U^y(1,x;g):=\log\sum_{\si\in\Omega}w(\si)
            \exp\bigl(\sqrt{2t}\,g\cdot v_\si+x\cdot\si+yb(\si)\bigr),
\end{equation*}
so that $U^0=U$. As for \eqref{e.def.terminal}, the $x$-derivatives of this terminal datum are cumulants of $\si$ under a Gibbs measure on $\Omega$, its gradient being the mean, so it satisfies the assumptions of Remark~\ref{r.terminal.data}, and thus Lemma~\ref{l.wellposed} applies to $U^y$; by the parameter dependence noted after Lemma~\ref{l.comparison}, it depends smoothly on $y$, and its $y$-derivatives solve the linearized equations written below. We define the response to the observable $b$ at level $r$ by
\begin{equation}\label{e.def.response}
 \Xi_r(b):=\dr_y^2U^y(r,X_r;g)\big|_{y=0}.
\end{equation}

\begin{lemma}[Gibbs averages and response along the branch]
\label{l.response}
Let $b:\Omega\to\R$. For every $r\in[0,1]$, we have 
\begin{align}
 \dr_yU^y(r,X_r;g)\big|_{y=0}&=\langle b\rangle_r,
      \label{e.response.first}\\
 \Xi_r(b)&=r\Var_r(b)
       +\int_r^1\E\bigl[\Var_u(b)\bigm|\cG_r\bigr]\,\d u.
      \label{e.response.second}
\end{align}
In particular, $\langle\si\rangle_r=\nabla U(r,X_r;g)=\bmm_r$, and the response matrix satisfies
\begin{equation}\label{e.response.matrix}
 K_r=rC_r+\int_r^1\E[C_u\mid\cG_r]\,\d u,
 \qquad
 rC_r\preceq K_r\preceq C_r\preceq N\Id.
\end{equation}
Moreover,
\begin{equation}\label{e.mean.martingale}
 \d\bmm_r=K_r\sqrt{2D(r)}\,\d B_r .
\end{equation}
\end{lemma}

\begin{proof}
Let $\mcl A:=\dr_r+\tr(D\nabla^2)+2r(D\nabla U)\cdot\nabla$ be the backward generator of the branch \eqref{e.branch.sde}. Differentiating \eqref{e.branch.pde} once and twice in $y$ at $y=0$, we find that $W:=\dr_yU^y|_{y=0}$ and $W_2:=\dr_y^2U^y|_{y=0}$ solve
\begin{equation*}
 \mcl AW=0,\qquad
 \mcl AW_2+2r\,\nabla W\cdot D\nabla W=0,
\end{equation*}
with terminal conditions $W(1,x)=\dr_yU^y(1,x;g)|_{y=0}$ and $W_2(1,x)=\dr_y^2U^y(1,x;g)|_{y=0}$. These are the mean and the variance of $b$ under the Gibbs measure with field $x$, so that $W(1,X_1)=\langle b\rangle_1$ and $W_2(1,X_1)=\Var_1(b)$ by \eqref{e.sampled.configuration}.

By It\^o's formula, $M_r:=W(r,X_r)$ is a bounded continuous martingale with terminal value $\langle b\rangle_1=\E[b(\si)\mid\cG_1]$. Hence $M_r=\E[b(\si)\mid\cG_r]=\langle b\rangle_r$, which is \eqref{e.response.first}, and $\d[M]_r=2\nabla W\cdot D\nabla W(r,X_r)\,\d r$. Set $V_r:=\Var_r(b)$. We claim that
\begin{equation}\label{e.response.martingale}
 \Xi_r(b)-rV_r+\int_0^rV_u\,\d u\qquad\text{is a martingale.}
\end{equation}
Since $\Xi_1(b)=W_2(1,X_1)=V_1$, the terminal value of this martingale is $\int_0^1V_u\,\d u$, and \eqref{e.response.second} follows by conditioning on $\cG_r$. To prove the claim, we note that three processes are martingales: first, $\Xi_r(b)+\int_0^ru\,\d[M]_u$, by It\^o's formula applied to $W_2(r,X_r)$ and the equation for $W_2$; second, $V_r+[M]_r=\E[b(\si)^2\mid\cG_r]-(M_r^2-[M]_r)$, by the tower property and the martingale property of $M^2-[M]$; third, $rY_r-\int_0^rY_u\,\d u$ for every bounded martingale $Y$, as a direct computation of conditional expectations shows. Applying the third fact to $Y:=V+[M]$ and subtracting the result from the first martingale, we find that
\begin{equation*}
 \Xi_r(b)-rV_r+\int_0^rV_u\,\d u
 +\Bigl(\int_0^ru\,\d[M]_u-r[M]_r+\int_0^r[M]_u\,\d u\Bigr)
\end{equation*}
is a martingale, and the expression in parentheses vanishes identically, by integration by parts for the continuous increasing process $[M]$. This proves \eqref{e.response.martingale}.

For $b(\si)=c\cdot\si$ with $c\in\R^N$, we have $U^y(1,x;g)=U(1,x+yc;g)$, so $U^y(r,x;g)=U(r,x+yc;g)$ by Lemma~\ref{l.symmetries}(1), and thus $W=c\cdot\nabla U$ and $W_2=c\cdot\nabla^2Uc$. Now \eqref{e.response.first} gives $c\cdot\nabla U(r,X_r;g)=\langle c\cdot\si\rangle_r$, and \eqref{e.response.second} gives $c\cdot K_rc=r\,c\cdot C_rc+\int_r^1\E[c\cdot C_uc\mid\cG_r]\,\d u$; polarization yields \eqref{e.response.matrix}. The lower bound on $K_r$ follows from $C_u\succeq0$. For the upper bound, the conditional variance decomposition of the vector $\si$ reads
\begin{equation*}
 C_r=\E[C_u\mid\cG_r]+\E\bigl[(\bmm_u-\bmm_r)(\bmm_u-\bmm_r)^{\mathsf T}\bigm|\cG_r\bigr]
 \succeq\E[C_u\mid\cG_r]
 \qquad(r\le u),
\end{equation*}
and $c\cdot C_rc\le\langle(c\cdot\si)^2\rangle_r\le N|c|^2$. Finally, differentiating \eqref{e.branch.pde} in $x$ shows that each component of $\nabla U$ solves $\mcl A(\dr_{x_i}U)=0$, and It\^o's formula gives \eqref{e.mean.martingale}.
\end{proof}

We note that at $r=0$, the identity \eqref{e.response.second} reads
\begin{equation}\label{e.response.zero}
 \Xi_0(b)=\int_0^1\E\bigl[\Var_u(b)\bigm|\cG_0\bigr]\,\d u,
\end{equation}
which will yield the formula for the time derivative of the free energy.

Lemma~\ref{l.response} and its proof apply verbatim to the branch started from any point $(r_0,x_0)$ in place of $(0,X_0)$. In particular, for every $(r,x)\in[0,1]\times\R^N$,
\begin{equation}\label{e.hessian.bounds}
 0\preceq\nabla^2U(r,x;g)\preceq N\Id,
 \qquad
 \tr\nabla^2U(r,x;g)\le N,
 \qquad
 |\nabla U(r,x;g)|^2\le N,
\end{equation}
all three bounds coming from \eqref{e.response.matrix} and \eqref{e.trace.covariance}, applied at the initial level of the branch started from $(r,x)$: there $\nabla^2U(r,x;g)\preceq C\preceq N\Id$ and $\tr\nabla^2U(r,x;g)\le\tr C=N-|\nabla U(r,x;g)|^2\le N$, where $C$ is the covariance of $\si$ at that level. In particular, $U(r,\cdot\,;g)$ is convex for every $r$.

\subsection{The gain in predicting the sampled configuration}
\label{s.gibbs.gain}
We define
\begin{equation}\label{e.def.Gamma}
 \Gamma(r):=\frac1N\E|\bmm_r|^2=\sum_{s\in\sS}p_{N,s}(r)
 =1-\frac1N\E\tr C_r=1-\frac1N\E|\si-\bmm_r|^2 ,
\end{equation}
where we used \eqref{e.trace.covariance}. Since $\bmm_r=\E[\si\mid\cG_r]$ is the best $\cG_r$-measurable prediction of $\si$ in mean square, and since the trivial prediction $0$ has normalized error $1$, the number $\Gamma(r)$ is the expected reduction of the normalized squared prediction error afforded by the information available at level $r$; equivalently, $1-\Gamma(r)$ measures how much of the sampled configuration remains unresolved at level $r$. Since $\bmm$ is a martingale, its increments over disjoint intervals are orthogonal, and for $u\le r$,
\begin{equation}\label{e.gain.increment}
 p_{N,s}(r)-p_{N,s}(u)=\frac1N\E\bigl|P_s(\bmm_r-\bmm_u)\bigr|^2,
 \qquad
 \Gamma(r)-\Gamma(u)=\frac1N\E|\bmm_r-\bmm_u|^2 ,
\end{equation}
so that the gain afforded by the levels between $u$ and $r$ is the mean square displacement of the prediction over these levels; in particular, $p_{N,s}$ and $\Gamma$ are increasing.

\begin{lemma}[Bounded gain in prediction]\label{l.energy}
The function $\Gamma$ is absolutely continuous on $[0,1]$ and
\begin{equation}\label{e.energy.identity}
 \Gamma'(r)=\frac2N\E\tr\bigl(K_rD(r)K_r\bigr),
 \qquad 0\le\Gamma(0)\le\Gamma(1)\le1 .
\end{equation}
Each component $p_{N,s}$ is increasing with values in $[0,\lambda_{N,s}]$, and $\sum_sp_{N,s}\le1$. If $\dot q_s\ge\kappa>0$ on $[0,1)$ for every species $s$, then
\begin{equation}\label{e.energy.bound}
 \int_0^1\E\tr K_r^2\,\d r\le\frac{N}{2\kappa}.
\end{equation}
\end{lemma}

\begin{proof}
By \eqref{e.mean.martingale}, It\^o's formula gives
\begin{equation*}
 \d|P_s\bmm_r|^2=2P_s\bmm_r\cdot\d\bmm_r
 +2\tr\bigl(P_sK_rD(r)K_rP_s\bigr)\d r,
\end{equation*}
whose drift is nonnegative. Taking expectations and summing over $s$ gives \eqref{e.energy.identity} and the monotonicity of $p_{N,s}$; the bounds follow from $|P_s\bmm_r|^2\le|I_{N,s}|$. Under the lower bound $D(r)\succeq\kappa\Id$ on the slopes, we have $\tr(K_rDK_r)\ge\kappa\tr K_r^2$, and integrating \eqref{e.energy.identity} gives \eqref{e.energy.bound}.
\end{proof}

The bound \eqref{e.energy.bound} turns the lower bound in the slope of the path into a control of the response matrices, and it will play a key role in the derivation of Theorem~\ref{t.remainder}. 
%
%
%
%
\section{Free-energy derivatives}\label{s.derivatives}

The main result of this section is an explicit representation of the derivatives of the free energy in terms of the Gibbs averages introduced above. From this follow an explicit representation of the remainder in the Hamilton--Jacobi equation for the free energy, and the Lipschitz continuity of the free energy as a function of the path.

\subsection{The derivative identities}
\label{s.derivatives.identities}
We call a continuous, piecewise affine $h:[0,1]\to\R^{\sS}$ an \emph{admissible variation}; for a regular path $q$, the path $q+\ep h$ is again regular for all $\ep$ small enough, and the solution of \eqref{e.branch.pde} depends smoothly on $\ep$, by the parameter dependence noted after Lemma~\ref{l.comparison}. We recall that we denote the overlap between two replicas by $R_N = R_N(\sigma^1, \sigma^2)$. 
\begin{theorem}[Derivative identities]\label{t.derivatives}
Let $t\ge0$ and let $q$ be a regular path. For every admissible variation $h$, we have
\begin{equation}\label{e.derivatives}
 \begin{aligned}
 \frac{\d}{\d\ep}\bar F_{N,\Omega}(t,q+\ep h)\Big|_{\ep=0}
   &=\int_0^1h(r)\cdot\E\langle R_N\rangle_r\,\d r,\\[2pt]
 \dr_t\bar F_{N,\Omega}(t,q)
   &=\int_0^1\E\langle\xi(R_N)\rangle_r\,\d r .
 \end{aligned}
\end{equation}
\end{theorem}

We write $\dr_q\bar F_{N,\Omega}(t,q)(r)=p_N(r)=\E\langle R_N\rangle_r$, as in \eqref{e.def.S.p}, for the density of the first variation. Both identities average the Gibbs measures over a level $r$ that is uniformly distributed on $[0,1]$. In the cascade representation of the free energy, this uniform distribution on $[0,1]$ arises as the law of the overlap between two cascade variables.

\begin{proof}[Proof of the time derivative]
We write $H_N(\si)=g\cdot v_\si=\sum_ag_ah_a(\si)$, where $h_a(\si)$ denotes the $a$-th coordinate of $v_\si$ in an orthonormal basis of the disorder space, so that $\sum_ah_a(\si)h_a(\tau)=N\xi(R_N(\si,\tau))$ by \eqref{e.gram}. For a deterministic vector $G$ in the disorder space, let $\Psi(G)$ denote $\E\,U(0,X_0;g)$ computed with the terminal condition $\log\sum_{\si\in\Omega}w(\si)e^{G\cdot v_\si+x\cdot\si}$; then $\E\,U(0,X_0;g)=\E\Psi(\sqrt{2t}g)$. Since $\sqrt{2t}g$ is a centered Gaussian vector with covariance~$2t\Id$, and since $\Psi$ is smooth with bounded second derivatives (by the polarized form of \eqref{e.response.zero}, $|\dr_{G_a}\dr_{G_b}\Psi|\le\|h_a\|_\infty\|h_b\|_\infty$), the heat equation gives, for every $t\ge0$,
\begin{equation*}
 \dr_t\E\,U(0,X_0;g)=\E\sum_a\dr^2_{G_a}\Psi(\sqrt{2t}g)
 =\sum_a\E\Xi_0(h_a),
\end{equation*}
because a change of $G_a$ perturbs the terminal Hamiltonian by the fixed observable $h_a$, and the second variation of $U(0,X_0;g)$ in that direction is the response \eqref{e.def.response}. By \eqref{e.response.zero} and \eqref{e.gram},
\begin{align*}
 \sum_a\E\Xi_0(h_a)
 &=\int_0^1\E\sum_a\Var_u(h_a)\,\d u\\
 &=\int_0^1\E\Bigl\langle N\xi\bigl(R_N(\si,\si)\bigr)
                  -N\xi\bigl(R_N(\si^1,\si^2)\bigr)\Bigr\rangle_u\,\d u\\
 &=N\xi(\lambda_N)-N\int_0^1\E\langle\xi(R_N)\rangle_u\,\d u ,
\end{align*}
since $R_N(\si,\si)=\lambda_N$ for every $\si\in\Sigma_N$. The self-overlap correction $t\xi(\lambda_N)$ in \eqref{e.free.energy.V} cancels the first term, and the second identity in \eqref{e.derivatives} follows. At $t=0$, the derivative of $t\mapsto\E\Psi(\sqrt{2t}g)$ is a right derivative, and the right-hand side is continuous in $t$.
\end{proof}

\begin{proof}[Proof of the path derivative]
Let $q^\ep:=q+\ep h$, let $U^\ep$ be the corresponding solution of \eqref{e.branch.pde}, and let $X^\ep_0\sim\mcl N(0,2\diag(q(0))+2\ep\diag(h(0)))$. Let $D_h(r):=\diag(\dot h(r))$. By \eqref{e.free.energy.V},
\begin{equation}\label{e.path.derivative.split}
 \frac{\d}{\d\ep}\bar F_{N,\Omega}(t,q^\ep)\Big|_{\ep=0}
 =\lambda_N\cdot h(1)
 -\frac1N\frac{\d}{\d\ep}\E\,U^\ep(0,X_0;g)\Big|_{\ep=0}
 -\frac1N\frac{\d}{\d\ep}\E\,U(0,X^\ep_0;g)\Big|_{\ep=0}.
\end{equation}
The last term is the derivative of a Gaussian expectation with respect to its covariance, and equals
\begin{equation}\label{e.root.term}
 -\frac1N\E\tr\bigl(\nabla^2U(0,X_0;g)\diag(h(0))\bigr)
 =-\frac1N\sum_sh_s(0)\,\E\tr K_{0,ss}.
\end{equation}
For the middle term, $W:=\dr_\ep U^\ep|_{\ep=0}$ solves the linearized equation
\begin{equation*}
 \mcl AW=-\tr(D_h\nabla^2U)-r\,\nabla U\cdot D_h\nabla U,
 \qquad W(1,\cdot)=0,
\end{equation*}
where $\mcl A$ is the backward generator of the branch, as in the proof of Lemma~\ref{l.response}. It\^o's formula along the branch gives
\begin{equation*}
 \E\,W(0,X_0)=\int_0^1\E\bigl[\tr(D_h(r)K_r)+r\,\bmm_r\cdot D_h(r)\bmm_r\bigr]\,\d r.
\end{equation*}
We now express the integrand through the Gibbs measures, for a general diagonal matrix. Let $c\in\R^{\sS}$. Since $C_r+\bmm_r\bmm_r^{\mathsf T}=\langle\si\si^{\mathsf T}\rangle_r$ has unit diagonal, the response identity \eqref{e.response.matrix} gives
\begin{equation*}
\begin{split}
 \tr(\diag(c)K_r)+r\,\bmm_r\cdot\diag(c)\bmm_r
 &=r\tr\bigl(\diag(c)\langle\si\si^{\mathsf T}\rangle_r\bigr)
     +\int_r^1\tr\bigl(\diag(c)\E[C_u\mid\cG_r]\bigr)\,\d u\\
 &=Nc\cdot\Ll(r\lambda_N
     +\int_r^1\E\bigl[\lambda_N-S_u\bigm|\cG_r\bigr]\,\d u\Rr),
\end{split}
\end{equation*}
where we used \eqref{e.trace.covariance} in the last step. Taking expectations and recalling that $\E S_u=p_N(u)$, we obtain, for every $c\in\R^{\sS}$ and $r\in[0,1]$,
\begin{equation}\label{e.a.identity}
 \frac1N\E\bigl[\tr(\diag(c)K_r)+r\,\bmm_r\cdot\diag(c)\bmm_r\bigr]
 =c\cdot\mfk a(r),
 \qquad
 \mfk a(r):=\lambda_N-\int_r^1p_N(u)\,\d u .
\end{equation}
We apply this identity twice. With $c=\dot h(r)$, so that $\diag(c)=D_h(r)$, it gives the middle term in \eqref{e.path.derivative.split}:
\begin{equation*}
 -\frac1N\frac{\d}{\d\ep}\E\,U^\ep(0,X_0;g)\Big|_{\ep=0}
 =-\frac1N\E\,W(0,X_0)
 =-\int_0^1\dot h(r)\cdot\mfk a(r)\,\d r .
\end{equation*}
With $c=h(0)$ and $r=0$, the quadratic term in $\bmm_r$ vanishes, and it gives the root term \eqref{e.root.term}:
\begin{equation*}
 -\frac1N\E\tr\bigl(\nabla^2U(0,X_0;g)\diag(h(0))\bigr)
 =-\frac1N\E\tr\bigl(\diag(h(0))K_0\bigr)
 =-h(0)\cdot\mfk a(0).
\end{equation*}
The last two terms in \eqref{e.path.derivative.split} together form the derivative of $\E\,U^\ep(0,X^\ep_0;g)$, in which $\ep$ enters both through the equation and through the initial law; we have shown that
\begin{equation}\label{e.V.path.derivative.raw}
 \frac1N\frac{\d}{\d\ep}\E\,U^\ep(0,X^\ep_0;g)\Big|_{\ep=0}
 =h(0)\cdot\mfk a(0)+\int_0^1\dot h(r)\cdot\mfk a(r)\,\d r .
\end{equation}
Since $h$ and $\mfk a$ are Lipschitz, an integration by parts gives
\begin{align*}
h(0)\cdot\mfk a(0)+  \int_0^1\dot h(r)\cdot\mfk a(r)\,\d r
 & =h(1)\cdot\mfk a(1)-\int_0^1h(r)\cdot\mfk a'(r)\,\d r,
\end{align*}
and using also that $\mfk a'=p_N$ and $\mfk a(1)=\lambda_N$, we thus obtain
\begin{equation}\label{e.V.path.derivative}
 \frac1N\frac{\d}{\d\ep}\E\,U^\ep(0,X^\ep_0;g)\Big|_{\ep=0}
 =\lambda_N\cdot h(1)-\int_0^1h(r)\cdot p_N(r)\,\d r .
\end{equation}
Inserting this into \eqref{e.path.derivative.split}, the self-overlap terms $\lambda_N\cdot h(1)$ cancel, and
\begin{equation*}
 \frac{\d}{\d\ep}\bar F_{N,\Omega}(t,q^\ep)\Big|_{\ep=0}
 =\int_0^1h(r)\cdot p_N(r)\,\d r
 =\int_0^1h(r)\cdot\E\langle R_N\rangle_r\,\d r ,
\end{equation*}
by the definition \eqref{e.def.S.p} of $p_N$. This is the first identity in \eqref{e.derivatives}.
\end{proof}

\subsection{The Hamilton--Jacobi remainder}
\label{s.derivatives.remainder}
For $p \in \mcl Q_\infty^\sS$, we define $\mcl H(p):=\int_0^1\xi(p(r))\,\d r$, and for a regular path $q$, we set
\begin{equation}\label{e.def.remainder}
 \msc E_N(t,q):=\dr_t\bar F_{N,\Omega}(t,q)
                     -\mcl H\bigl(\dr_q\bar F_{N,\Omega}(t,q)\bigr).
\end{equation}
Theorem~\ref{t.derivatives} gives at once
\begin{equation}\label{e.defect.identity}
 \msc E_N(t,q)=\int_0^1\Bigl\{\E\langle\xi(R_N)\rangle_r
                 -\xi\bigl(\E\langle R_N\rangle_r\bigr)\Bigr\}\,\d r .
\end{equation}
The main result of the next section is that this term is small. We note here that it is controlled by the variance of the overlap. Let $M_\xi:=\sup_{x\in[-1,1]^{\sS}}\|\nabla^2\xi(x)\|$. Since the linear term of the Taylor expansion of $\xi$ around $p_N(r)$ has zero mean under $\E\langle\cdot\rangle_r$, we have
\begin{equation}\label{e.defect.variance}
 |\msc E_N(t,q)|\le\frac{M_\xi}{2}\msc W_N(t,q),\qquad \text{ with } \ 
 \msc W_N(t,q):=\int_0^1\E\bigl\langle|R_N-p_N(r)|^2\bigr\rangle_r\,\d r .
\end{equation}
The conditional variance decomposition separates two sources of fluctuations:
\begin{equation}\label{e.W.split}
 \E\bigl\langle|R_N-p_N(r)|^2\bigr\rangle_r
 =\E\bigl\langle|R_N-S_r|^2\bigr\rangle_r
  +\E\bigl|S_r-p_N(r)\bigr|^2 .
\end{equation}
The first term concerns two replicas within the Gibbs measure at a fixed disorder and branch; the second concerns the fluctuations of the Gibbs mean overlap across the disorder and the branch. Writing $P_s\si^a=P_s\bmm_r+\theta^a_s$ for the two replicas, we have $NR_{N,s}=|P_s\bmm_r|^2+P_s\bmm_r\cdot\theta^2_s+\theta^1_s\cdot P_s\bmm_r+\theta^1_s\cdot\theta^2_s$, and the three fluctuating terms are uncorrelated under $\langle\cdot\rangle_r$, with variances $\bmm_r\cdot C_{r,ss}\bmm_r$ (twice) and $\tr(C_{r,ss}^2)$. With \eqref{e.W.split}, this gives
\begin{equation}\label{e.overlap.decomposition}
 \E\bigl\langle(R_{N,s}-p_{N,s}(r))^2\bigr\rangle_r
 =\frac{\E\tr(C_{r,ss}^2)
       +2\,\E\bigl[\bmm_r\cdot C_{r,ss}\bmm_r\bigr]}{N^2}
       +\Var(S_{r,s}).
\end{equation}
The first term in \eqref{e.overlap.decomposition} will be controlled by the response bound $K_r\succeq rC_r$ and \eqref{e.energy.bound}. Obtaining a control of the other two terms will be the main focus of the next section.

\subsection{Lipschitz continuity and extension to all paths}
\label{s.derivatives.extension}

We set $\mathbf1:=(1,\ldots,1)\in\R^{\sS}$.

\begin{proposition}[Lipschitz bound and extension]\label{p.extension}
For all regular paths $q,q'$ and all $t,t'\ge0$, we have
\begin{equation}\label{e.finite.lipschitz}
 |\bar F_{N,\Omega}(t,q)-\bar F_{N,\Omega}(t',q')|
 \le\sum_{s\in\sS}\lambda_{N,s}\|q_s-q'_s\|_{L^1}+\xi(\mathbf1)\,|t-t'|
 \le\|q-q'\|_{L^1}+\xi(\mathbf1)\,|t-t'|.
\end{equation}
Consequently, $\bar F_{N,\Omega}$ extends uniquely to a function on $[0,\infty)\times\cQ_1^{\sS}$ satisfying \eqref{e.finite.lipschitz}, and we denote this extension by the same symbol. For a slice, this extension is the enriched slice free energy $\bar F_{N,k}(t,q)$ appearing in Theorem~\ref{t.slices}.
\end{proposition}

\begin{proof}
The set of regular paths is convex, so $q+\ep(q'-q)$ is regular for $\ep\in[0,1]$, with the admissible variation $h=q'-q$. By Theorem~\ref{t.derivatives} and $0\le p_{N,s}(r)\le\lambda_{N,s}$, the derivative along this segment is bounded in absolute value by $\sum_s\lambda_{N,s}\|q'_s-q_s\|_{L^1}$, and since $\lambda_{N,s}\ge0$ and $\sum_s\lambda_{N,s}=1$, we have $\sum_s\lambda_{N,s}|x_s|\le\max_s|x_s|\le|x|$ for every $x\in\R^{\sS}$, which gives the second inequality in \eqref{e.finite.lipschitz} after integration in $r$. Since $R_N(\si^1,\si^2)\in[-1,1]^{\sS}$ and the coefficients of $\xi$ are nonnegative, $|\dr_t\bar F_{N,\Omega}|\le\sup_{[-1,1]^{\sS}}|\xi|=\xi(\mathbf1)$. Regular paths are dense in $\cQ_1^{\sS}$: a path in $\cQ_1^{\sS}$ is approximated in $L^1$ by bounded paths, which are approximated by continuous piecewise affine increasing paths, and adding $\eta(1+r)\mathbf1$ with small $\eta>0$ makes them regular.
\end{proof}

The path $q=0$ is not regular, but substituting it formally in Definition~\ref{d.free.energy} gives $X\equiv0$ and $A_q=0$, whatever $U$ is, hence the ordinary free energy. The extension indeed takes this value. To see this, we will use the heat interpolation identity
\begin{equation}\label{e.heat.interpolation}
 \E f(V+\sqrt{2\ep}Z)-\E f(V)=\int_0^\ep\E\,\Delta f(V+\sqrt{2s}Z)\,\d s,
\end{equation}
valid for $\ep\ge0$, a random vector $V$, an independent standard Gaussian vector $Z$, and a $C^2$ function $f$ with bounded derivatives. For the functions $U(r,\cdot\,;g)$, whose Laplacian lies in $[0,N]$ by \eqref{e.hessian.bounds}, the left-hand side lies in $[0,N\ep]$.

\begin{lemma}[The value at $q=0$]\label{l.zero}
For every $t\ge0$,
\begin{equation}\label{e.zero}
 \bar F_{N,\Omega}(t,0)=-\frac1N\E\log\sum_{\si\in\Omega}w(\si)
 \exp\bigl(\sqrt{2t}H_N(\si)-Nt\xi(\lambda_N)\bigr),
\end{equation}
which is the free energy \eqref{e.def.slice.zero} for a slice, and the free energy \eqref{e.def.free.energy} for the full cube with the weights $2^{-N}e^{h^{(N)}\cdot\si}$.
\end{lemma}

\begin{proof}
For $\ep>0$, the path $q^\ep(r):=\ep(1+r)\mathbf1$ is regular, with $D^\ep=\ep\Id$ and initial field $X^\ep_0=\sqrt{2\ep}\,Z$, and $\bar F_{N,\Omega}(t,q^\ep)\to\bar F_{N,\Omega}(t,0)$ as $\ep\to0$ by \eqref{e.finite.lipschitz}. Let $U^\ep$ be the solution of \eqref{e.branch.pde} for $q^\ep$, with operator $\mcl L^\ep$, and let $W(r,x):=U(1,x;g)$ be the terminal function \eqref{e.def.terminal}, held constant in $r$. Since $\nabla^2W$ is a covariance matrix of $\si$, of trace at most $N$, and $|\nabla W|^2\le N$, we have $\mcl L^\ep[W]=\ep(\tr\nabla^2W+r|\nabla W|^2)\in[0,2\ep N]$, while $\mcl L^\ep[U^\ep]=0$ and both functions have the same terminal datum, so Lemma~\ref{l.comparison} gives $|U^\ep(0,\cdot)-U(1,\cdot\,;g)|\le2\ep N$. By \eqref{e.heat.interpolation}, $0\le\E\,U(1,\sqrt{2\ep}\,Z;g)-\E\,U(1,0;g)\le N\ep$. Hence $|\E\,U^\ep(0,X^\ep_0;g)-\E\,U(1,0;g)|\le3\ep N$, and \eqref{e.free.energy.V}, in which $\lambda_N\cdot q^\ep(1)=2\ep$, gives
\begin{equation*}
 \Bigl|\bar F_{N,\Omega}(t,q^\ep)-t\xi(\lambda_N)+\frac1N\E\,U(1,0;g)\Bigr|\le5\ep .
\end{equation*}
Letting $\ep\to0$ yields $\bar F_{N,\Omega}(t,0)=t\xi(\lambda_N)-\frac1N\E\,U(1,0;g)$, which is \eqref{e.zero}, since $U(1,0;g)=\log\sum_{\si\in\Omega}w(\si)e^{\sqrt{2t}H_N(\si)}$.
\end{proof}
%
%
%
%
\section{Concentration of the overlap on a slice}
\label{s.concentration}

In this section, $\Omega=\Omega_{N,k}$ is a slice, with the weights $w\equiv2^{-N}$, and we show that the overlap variance $\msc W_N(t,q)$ of \eqref{e.defect.variance} is small, uniformly over the slices and over $t$, as soon as the slopes of the path are bounded below. By \eqref{e.defect.variance}, this controls the remainder $\msc E_N(t,q)$ of \eqref{e.def.remainder}, where we recall that $ M_\xi=\sup_{x\in[-1,1]^{\sS}}\|\nabla^2\xi(x)\|$.

\begin{theorem}[Quantitative remainder estimate on a slice]\label{t.remainder}
Let $\kappa>0$, and let $q$ be a regular path with $\dot q_s\ge\kappa$ on every piece, for every species $s$. For every slice $\Omega=\Omega_{N,k}$, every $\de\in\{1,\frac12,\frac13,\ldots\}$ and every $t\ge0$, we have
\begin{equation}\label{e.W.bound}
 \msc W_N(t,q)\le13\de+\frac{5}{2\kappa N\de^2}.
\end{equation}
Consequently, for every $N$, every slice, every $t\ge0$ and every such path, we have
\begin{equation}\label{e.small.defect}
 \Bigl|\dr_t\bar F_{N,k}(t,q)
              -\mcl H\bigl(\dr_q\bar F_{N,k}(t,q)\bigr)\Bigr|
                          \le15M_\xi(\kappa N)^{-1/3}.
\end{equation}
\end{theorem}

The argument is adapted from \cite{aronow2026quantitative}, where it is used to prove the Parisi formula for the SK model. As noted after \eqref{e.overlap.decomposition}, the first term there is controlled by the response bound $K_r\succeq rC_r$ of \eqref{e.response.matrix} and by the bounded gain in prediction of Lemma~\ref{l.energy}, which makes the response matrices small in an integrated sense once the slopes are bounded below. The other two terms, which involve the posterior mean $\bmm_r$, are the main concern of this section. We first show that they are bounded by the variance of a single observable which we call the realized gain in prediction~$L_{r,s}$, defined in~\eqref{e.def.L}. To bound this variance, we reverse the direction of the conditioning: given the sampled configuration, the Gaussian data have a log-concave law, and the Brascamp--Lieb inequality bounds the variance of averages of the realized gain over levels in terms of the response matrices; the symmetry of the slice under permutations of the sites within each species then removes the conditioning on the configuration. This is Proposition~\ref{p.posterior}. A last step, based on the martingale property of the posterior mean, passes from averages over levels to the integrated variance of the realized gain; the parameter $\de$ in \eqref{e.W.bound} is the length of the intervals over which the realized gain is averaged.

\subsection{The realized gain in prediction}
\label{s.concentration.observable}

For a level $r$ and a species $s$, we define
\begin{equation}\label{e.def.L}
 L_{r,s}
 :=\lambda_{N,s}-\frac1N\bigl|P_s(\si-\bmm_r)\bigr|^2
 =\frac1N\Ll(2\si\cdot P_s\bmm_r-|P_s\bmm_r|^2\Rr),
\end{equation}
where we used $|P_s\si|^2=|I_{N,s}|$. This is the gain in predicting $P_s\si$ by $P_s\bmm_r$ rather than by $0$, realized along the branch rather than averaged, and we call it the realized gain in prediction: summing over the species gives $\sum_sL_{r,s}=1-\frac1N|\si-\bmm_r|^2$, whose expectation is $\Gamma(r)$ by \eqref{e.def.Gamma}, and likewise $\E L_{r,s}=p_{N,s}(r)$ by \eqref{e.trace.covariance}. The realized gain is thus the quenched, per-species counterpart of the expected gain of Section~\ref{s.gibbs.gain}. It is determined by the Gibbs measure at level $r$ and the sampled configuration: with $\si$ held fixed and two further replicas $\tau^1,\tau^2$ at level $r$, averaging only over $\tau^1,\tau^2$, we have
\begin{equation*}
 L_{r,s}=\Bigl\langle 2R_{N,s}(\si,\tau^1)
                    -R_{N,s}(\tau^1,\tau^2)\Bigr\rangle_r .
\end{equation*}
The conditional mean of $L_{r, s}$ given $\cG_r$ is $S_{r,s}$, and its conditional variance given $\cG_r$ equals $4N^{-2}\bmm_r\cdot C_{r,ss}\bmm_r$. Hence
\begin{equation}\label{e.L.decomposition}
 \Var(L_{r,s})
 =\frac{4\,\E\bigl[\bmm_r\cdot C_{r,ss}\bmm_r\bigr]}{N^2}
                         +\Var(S_{r,s}).
\end{equation}
Comparing \eqref{e.L.decomposition} with \eqref{e.overlap.decomposition}, we obtain 
\begin{equation}\label{e.overlap.reduction}
\begin{split}
 \E\bigl\langle(R_{N,s}-p_{N,s}(r))^2\bigr\rangle_r
 &=\frac{\E\tr(C_{r,ss}^2)}{N^2}+\Var(L_{r,s})
 -\frac{2\,\E\bigl[\bmm_r\cdot C_{r,ss}\bmm_r\bigr]}{N^2}\\
 &\le\frac{\E\tr(C_{r,ss}^2)}{N^2}+\Var(L_{r,s}).
\end{split}
\end{equation}
Integrating over $r$ and summing over $s$ leads to
\begin{equation}\label{e.W.reduction}
 \msc W_N(t,q)\le\frac1{N^2}\int_0^1\E\sum_s\tr(C_{r,ss}^2)\,\d r
                         +\int_0^1\sum_s\Var(L_{r,s})\,\d r .
\end{equation}

\subsection{The joint law of the disorder, the branch and the sampled configuration}
\label{s.concentration.law}

Let $Z_r:=X_0+\int_0^r\sqrt{2D(u)}\,\d B_u$ be the Gaussian part of the branch, and let $\ga$ denote the law of the pair $(g,Z)$, that is, the Gaussian disorder together with the Gaussian field path with initial covariance $2\diag(q(0))$ and increment covariance $2D(r)\,\d r$.

\begin{lemma}[Joint density]\label{l.joint.density}
The law of $(g,X,\si)$ has density
\begin{equation}\label{e.joint.density}
 2^{-N}\exp\Ll(\sqrt{2t}\,g\cdot v_\si+x_1\cdot\si
                         -\int_0^1U(r,x_r;g)\,\d r\Rr)
\end{equation}
at $(g,x,\si)$ with respect to the product of $\ga$ and the counting measure on $\Omega$.
\end{lemma}

\begin{proof}
Given $g$, the branch has drift $2rD\nabla U$ relative to the Gaussian path. Since $|\nabla U|\le\sqrt N$, Girsanov's formula gives the density
\begin{equation*}
 \exp\Ll(\int_0^1r\,\nabla U\cdot\sqrt{2D}\,\d\tilde B_r
          -\int_0^1r^2\,\nabla U\cdot D\nabla U\,\d r\Rr)
\end{equation*}
of the branch law with respect to the Gaussian path law. Here $x=(x_r)_{0\le r\le1}$ is the path variable, $\nabla U$ is evaluated at $(r,x_r;g)$, and $\tilde B_r:=\int_0^r(2D(u))^{-1/2}\,\d x_u$. Since $D$ is piecewise constant and invertible by \eqref{e.def.regular}, $\tilde B$ is a linear function of the path. Under the Gaussian path law, $\tilde B$ is a standard $N$-dimensional Brownian motion and $\d x_r=\sqrt{2D(r)}\,\d\tilde B_r$, whereas along the branch $\tilde B_r=B_r+\int_0^ru\sqrt{2D(u)}\,\nabla U(u,X_u;g)\,\d u$. Under the Gaussian law, It\^o's formula and \eqref{e.branch.pde} give, as in the proof of Lemma~\ref{l.ito}, $\d U(r,x_r;g)=-r\,\nabla U\cdot D\nabla U\,\d r+\nabla U\cdot\sqrt{2D}\,\d\tilde B_r$, so that the exponent equals
\begin{equation*}
 \int_0^1r\,\d U(r,x_r;g)=U(1,x_1;g)-\int_0^1U(r,x_r;g)\,\d r .
\end{equation*}
Multiplying by the Gibbs probability \eqref{e.sampled.configuration} cancels $U(1,x_1;g)$.
\end{proof}

For a fixed configuration $\si$, the exponent in \eqref{e.joint.density} is the sum of a linear function of the Gaussian data and of $-\int_0^1U(r,x_r;g)\,\d r$. This second part is concave:

\begin{lemma}[Convexity in the Gaussian data]\label{l.convexity}
For every $r\in[0,1]$, the map $(x,g)\mapsto U(r,x;g)$ is smooth, with derivatives of all orders bounded and continuous in $(r,x,g)$, and it is convex, with
\begin{equation}\label{e.joint.hessian.bound}
 0\preceq\nabla^2_{(x,g)}U(r,x;g)\preceq N\bigl(1+2t\xi(\lambda_N)\bigr)\Id .
\end{equation}
\end{lemma}

\begin{proof}
The terminal datum \eqref{e.def.terminal} is smooth in $(x,g)$ with bounded derivatives of all orders, so the regularity follows from Lemma~\ref{l.wellposed} and the smooth dependence on parameters noted after Lemma~\ref{l.comparison}. Fix a direction $(\hat x,\hat g)$ with $|\hat x|^2+|\hat g|^2=1$. By Lemma~\ref{l.symmetries}(1) and the uniqueness in Lemma~\ref{l.wellposed}, we have $U(r,x+y\hat x;g+y\hat g)=U^y(r,x;g)$, where $U^y$ solves the equation with the terminal Hamiltonian perturbed by $y\,b(\si)$ with $b(\si):=\sqrt{2t}\,\hat g\cdot v_\si+\hat x\cdot\si$. The proof of Lemma~\ref{l.response}, applied to the branch started from $(r,x)$, shows that $\dr_y^2U^y(r,x;g)|_{y=0}$ is a nonnegative combination of conditional variances of $b$, with total weight one, as in \eqref{e.response.second}. This gives the convexity, and the upper bound in \eqref{e.joint.hessian.bound} follows from $b(\si)^2\le(|\hat x|^2+|\hat g|^2)(|\si|^2+2t|v_\si|^2)=N(1+2t\xi(\lambda_N))$, by the Cauchy--Schwarz inequality and \eqref{e.gram}.
\end{proof}

\subsection{The averaged concentration estimate}
\label{s.concentration.conditioning}

\begin{proposition}[Variance of level averages]\label{p.posterior}
Let $q$ be regular. For every slice $\Omega_{N,k}$, every deterministic $a\in L^2([0,1])$ and every species $s$, we have
\begin{equation}\label{e.posterior.variance}
 \Var\Ll(\int_0^1a(r)L_{r,s}\,\d r\Rr)
 \le\frac4{N^2}\int_0^1a(r)^2\,
                   \E\tr\bigl(K_{r,ss}C_{r,ss}\bigr)\,\d r .
\end{equation}
\end{proposition}

The main step of the proof is Lemma~\ref{l.conditional.variance} below, which bounds the variance of the level average conditionally on the sampled configuration. Given $\si$, the law of the Gaussian data $(g,X)$ has, by Lemma~\ref{l.joint.density}, a density with respect to $\ga$ which is the exponential of a linear function of $(g,x)$ minus $\int_0^1U(r,x_r;g)\,\d r$, and $U$ is convex in $(x,g)$ by Lemma~\ref{l.convexity}. The conditional law is therefore log-concave, and the Hessian of minus its log-density is the inverse covariance of $\ga$ plus the integral over $r$ of the Hessian of $U(r,\cdot\,;\cdot)$. By the Brascamp--Lieb inequality, the conditional variance of a function of the Gaussian data is at most the expectation of its squared gradient, measured in the inverse of this Hessian. The level average is a quadratic function of the residuals $P_s(\si-\bmm_r)$, and $\bmm_r=\nabla U(r,X_r;g)$; its gradient in the Gaussian data is therefore an integral over $r$ of the Hessian of $U(r,\cdot\,;\cdot)$ applied to the residual. The squared gradient thus contains two such Hessians, and the next lemma shows that the inverse Hessian cancels one of them, leaving the quadratic forms $P_s(\si-\bmm_r)\cdot K_rP_s(\si-\bmm_r)$. The symmetry of the slice then makes the level average independent of $\si$, so that its conditional variance is its variance, and the expectation of the quadratic forms is computed by conditioning on $\cG_r$, where the residual is centered with covariance $C_{r,ss}$.

\begin{lemma}[Cancellation of one Hessian]\label{l.algebra}
Let $(A_r)_{0\le r\le1}$ be a bounded measurable family of positive semidefinite matrices of the same size, let $a\in L^2([0,1])$, let $(v_r)_{0\le r\le1}$ be a bounded measurable family of vectors, and let $Q$ be a positive definite matrix with $Q\succeq\int_0^1A_r\,\d r$. The vector $z:=\int_0^1a(r)A_rv_r\,\d r$ satisfies
\begin{equation}\label{e.curvature.algebra}
 z\cdot Q^{-1}z\le\int_0^1a(r)^2\;v_r\cdot A_rv_r\,\d r .
\end{equation}
\end{lemma}

\begin{proof}
Using the assumption on $Q$, we can write, for every vector $u$,
\begin{equation*}
 2u\cdot z-u\cdot Qu
 \le\int_0^1\bigl(2a(r)\,u\cdot A_rv_r-u\cdot A_ru\bigr)\,\d r
 \le\int_0^1a(r)^2\;v_r\cdot A_rv_r\,\d r ,
\end{equation*}
where we used that $(u-a(r)v_r)\cdot A_r(u-a(r)v_r)\ge0$ for the second inequality. 
Taking $u=Q^{-1}z$, the left-hand side equals $z\cdot Q^{-1}z$, and we thus obtain the announced result.
\end{proof}

\begin{lemma}[Concentration given the sampled configuration]\label{l.conditional.variance}
Let $q$ be regular, let $\Omega=\Omega_{N,k}$ be a slice, let $a\in L^2([0,1])$ be deterministic, and let $s$ be a species. For every $\tau\in\Omega$,
\begin{equation}\label{e.conditional.variance}
 \Var\Ll(\int_0^1a(r)L_{r,s}\,\d r\Bigm|\si=\tau\Rr)
 \le\frac4{N^2}\,\E\Ll[\int_0^1a(r)^2\,
   P_s(\si-\bmm_r)\cdot K_rP_s(\si-\bmm_r)\,\d r\Bigm|\si=\tau\Rr].
\end{equation}
\end{lemma}

\begin{proof}
\emph{Step 1: the conditional law.}
Fix $\tau\in\Omega$, and let $f_\tau(g,x)$ denote the density \eqref{e.joint.density} at $(g,x,\tau)$. Then $c_\tau:=\int f_\tau\,\d\ga=\P(\si=\tau)>0$, and $\mu_\tau:=c_\tau^{-1}f_\tau\,\ga$ is the conditional law of $(g,X)$ given $\si=\tau$. For a realization $g$ of the disorder and a continuous path $x$, set
\begin{equation}\label{e.def.conditional.observables}
\begin{gathered}
 \theta_r(g,x):=P_s\bigl(\tau-\nabla U(r,x_r;g)\bigr),
 \qquad
 F(g,x):=\frac1N\int_0^1a(r)\Ll(|I_{N,s}|-\bigl|\theta_r(g,x)\bigr|^2\Rr)\d r,\\
 \mcl E(g,x):=\int_0^1a(r)^2\,\theta_r(g,x)\cdot\nabla^2U(r,x_r;g)\,\theta_r(g,x)\,\d r .
\end{gathered}
\end{equation}
On the event $\{\si=\tau\}$, by \eqref{e.def.L}, $F(g,X)$ is the level average $\int_0^1a(r)L_{r,s}\,\d r$, and $\mcl E(g,X)$ is the integral on the right-hand side of \eqref{e.conditional.variance}. Thus \eqref{e.conditional.variance} reads
\begin{equation}\label{e.conditional.variance.mu}
 \Var_{\mu_\tau}(F)\le\frac4{N^2}\,\E_{\mu_\tau}[\mcl E] .
\end{equation}
Since $\nabla U\in[-1,1]^N$, $|\theta_r|\le2\sqrt{|I_{N,s}|}$ and $0\preceq\nabla^2U\preceq N\Id$ by \eqref{e.hessian.bounds}, we have
\begin{equation}\label{e.observable.bounds}
 |F|\le3\lambda_{N,s}\|a\|_{L^1},
 \qquad
 0\le\mcl E\le4N|I_{N,s}|\,\|a\|_{L^2}^2 .
\end{equation}
Moreover, $\nabla U$ and $\nabla^2U$ are Lipschitz continuous in $x$, uniformly in $r$, by Lemma~\ref{l.wellposed}, so that $F(g,\cdot)$ and $\mcl E(g,\cdot)$ are continuous for the uniform norm on paths.

\emph{Step 2: finite-dimensional approximations of the conditional law.}
For $n\ge1$ and $r_j:=j/n$, let $x^{[n]}$ denote the polygonal interpolation of the path $x$ on the mesh $(r_j)_{0\le j\le n}$, that is, the continuous path which agrees with $x$ at the mesh points and is affine on each $[r_j,r_{j+1}]$. It is a linear function of $(x_{r_0},\ldots,x_{r_n})$, it satisfies $\|x^{[n]}\|_\infty\le\|x\|_\infty$, and $x^{[n]}\to x$ uniformly as $n\to\infty$ for every continuous path $x$. Let
\begin{equation}\label{e.approximate.density}
 f_{n,\tau}(g,x):=2^{-N}\exp\Ll(\sqrt{2t}\,g\cdot v_\tau+x_1\cdot\tau
                         -\int_0^1U\bigl(r,x^{[n]}_r;g\bigr)\,\d r\Rr),
\end{equation}
which is $f_\tau$ with the path $x$ replaced by $x^{[n]}$ inside the integral, the function $U$ being unchanged. We first show that these functions are integrable. By Lemma~\ref{l.convexity}, the first derivatives of $U$ in $(x,g)$ are bounded, and $r\mapsto U(r,0;0)$ is continuous on $[0,1]$, so that $|U(r,x;g)|\le C(1+|x|+|g|)$ for some constant $C$. Since $\|x^{[n]}\|_\infty\le\|x\|_\infty$, we obtain
\begin{equation}\label{e.density.envelope}
 f_{n,\tau}(g,x)\le\exp\bigl(C'(1+|g|+\|x\|_\infty)\bigr)
\end{equation}
for a constant $C'$ independent of $n$, and the same bound holds for $f_\tau$. Under $\ga$, the vector~$g$ is Gaussian, and each coordinate $i\in I_{N,s}$ of the path is a Gaussian initial value plus a Brownian motion run with the clock $2(q_s(r)-q_s(0))$, so $|g|$ and $\|x\|_\infty$ have Gaussian tails, and the right-hand side of \eqref{e.density.envelope} is $\ga$-integrable. Hence $c_{n,\tau}:=\int f_{n,\tau}\,\d\ga$ is finite and positive, and $\mu_{n,\tau}:=c_{n,\tau}^{-1}f_{n,\tau}\,\ga$ is a probability measure. Since $|U(r,x;g)-U(r,x';g)|\le\sqrt N|x-x'|$, we also have $f_{n,\tau}\to f_\tau$ pointwise, and dominated convergence gives $c_{n,\tau}\to c_\tau$. Consequently, if $(G_n)$ is a uniformly bounded sequence of measurable functions of $(g,x)$ converging pointwise to $G$, then
\begin{equation}\label{e.moments.convergence}
 \E_{\mu_{n,\tau}}[G_n]=\frac{\int G_nf_{n,\tau}\,\d\ga}{c_{n,\tau}}
 \longrightarrow\frac{\int Gf_\tau\,\d\ga}{c_\tau}=\E_{\mu_\tau}[G],
\end{equation}
again by dominated convergence. We write $F_n(g,x):=F(g,x^{[n]})$ and $\mcl E_n(g,x):=\mcl E(g,x^{[n]})$ for the functions in \eqref{e.def.conditional.observables} evaluated along the interpolated path.

\emph{Step 3: Gaussian coordinates and the log-density.}
Let $Y:=(g,x_{r_0},\ldots,x_{r_n})\in\R^d$, where $d$ is the dimension of the disorder space plus $(n+1)N$. Under $\ga$, this is a centered Gaussian vector whose covariance matrix $\Sigma$ is positive definite: $g$ is a standard Gaussian vector, and, independently of $g$, the initial value $x_{r_0}$ and the increments $x_{r_j}-x_{r_{j-1}}$ are independent and centered, with covariances $2\diag(q(0))$ and $2\diag(q(r_j)-q(r_{j-1}))$, which are positive definite since $q$ is regular. The interpolated path $x^{[n]}$ is a linear function of $Y$, so the density $f_{n,\tau}$ and the functions $F_n$ and $\mcl E_n$ depend on $(g,x)$ only through $Y$. The law of $Y$ under $\mu_{n,\tau}$ therefore has, relative to Lebesgue measure, a density proportional to $e^{-\Phi_{n,\tau}}$, where
\begin{equation*}
 \Phi_{n,\tau}(Y):=\frac12Y\cdot\Sigma^{-1}Y-\sqrt{2t}\,g\cdot v_\tau-x_1\cdot\tau
                   +\int_0^1U\bigl(r,x^{[n]}_r;g\bigr)\,\d r .
\end{equation*}
By Lemma~\ref{l.convexity}, $U(r,\cdot\,;\cdot)$ is smooth and convex, with bounded derivatives, and $Y\mapsto(g,x^{[n]}_r)$ is linear. Hence $A_r:=\nabla_Y^2[U(r,x^{[n]}_r;g)]$ is positive semidefinite, bounded and continuous in $(r,Y)$, the function $\Phi_{n,\tau}$ is smooth, and, differentiating under the integral sign,
\begin{equation}\label{e.conditional.curvature}
 \nabla^2\Phi_{n,\tau}=Q:=\Sigma^{-1}+\int_0^1A_r\,\d r\succeq\Sigma^{-1} .
\end{equation}
The Hessian in \eqref{e.conditional.curvature} exceeds the inverse covariance $\Sigma^{-1}$ of the Gaussian coordinates by the positive semidefinite term $\int_0^1A_r\,\d r$, which involves differentiation in the disorder $g$ as well as in the field coordinates. 

\emph{Step 4: the gradient of the level average.}
Let $\msf J:\R^N\to\R^d$ be the linear map $\msf Ju:=(0,u,\ldots,u)$, which leaves $g$ unchanged and shifts every $x_{r_j}$, hence the whole path $x^{[n]}$, by $u$. Thus $U(r,x^{[n]}_r;g)$ evaluated at $Y+\msf Ju$ equals $U(r,x^{[n]}_r+u;g)$, and differentiating once and twice in $u$ at $u=0$ gives, with $\bmm^{[n]}_r:=\nabla U(r,x^{[n]}_r;g)$ and $K^{[n]}_r:=\nabla^2U(r,x^{[n]}_r;g)$,
\begin{equation}\label{e.shift.identities}
 \msf J^{\mathsf T}\nabla_Y\bigl[U(r,x^{[n]}_r;g)\bigr]=\bmm^{[n]}_r,
 \qquad
 \msf J^{\mathsf T}A_r\msf J=K^{[n]}_r .
\end{equation}
Let $\theta^{[n]}_r:=P_s(\tau-\bmm^{[n]}_r)=\theta_r(g,x^{[n]})$ be the residual along the interpolated path. By the first identity in \eqref{e.shift.identities}, the Jacobian matrix of $\theta^{[n]}_r$ as a function of $Y$ is $-P_s\msf J^{\mathsf T}A_r$. Since $A_r$ is symmetric and $P_s\theta^{[n]}_r=\theta^{[n]}_r$, differentiating $F_n$ with $\tau$ fixed gives, with $v_r:=\msf J\theta^{[n]}_r$,
\begin{equation}\label{e.gradient.matching}
 N\nabla_YF_n=2\int_0^1a(r)A_rv_r\,\d r,
 \qquad
 v_r\cdot A_rv_r=\theta^{[n]}_r\cdot K^{[n]}_r\theta^{[n]}_r,
\end{equation}
the second identity being the second one in \eqref{e.shift.identities}. The positive initial variance $q(0)>0$ was used to make $\Sigma$ invertible, but $\Sigma$ will disappear from the estimate once Lemma~\ref{l.algebra} is applied. Since $A_r$ and $\theta^{[n]}_r$ are bounded, $F_n$ is a smooth function of $Y$ with a bounded gradient.

\emph{Step 5: the Brascamp--Lieb inequality.}
We use the Brascamp--Lieb inequality \cite{brascamp1976extensions}, in the following form: if a probability measure on $\R^d$ has a density proportional to $e^{-\Phi}$, with $\Phi$ smooth and $\nabla^2\Phi\succeq\lambda\Id$ for some $\lambda>0$, then $\Var(G)\le\E[\nabla G\cdot(\nabla^2\Phi)^{-1}\nabla G]$ for every smooth $G$ with bounded gradient. Applied to $\Phi_{n,\tau}$ and $F_n$, together with \eqref{e.conditional.curvature}, \eqref{e.gradient.matching} and Lemma~\ref{l.algebra} applied to $z:=\frac N2\nabla_YF_n$, it gives
\begin{equation}\label{e.finite.estimate}
\begin{split}
 \Var_{\mu_{n,\tau}}(F_n)
 &\le\E_{\mu_{n,\tau}}\bigl[\nabla_YF_n\cdot Q^{-1}\nabla_YF_n\bigr]\\
 &\le\frac4{N^2}\,\E_{\mu_{n,\tau}}\int_0^1a(r)^2\,
     \theta^{[n]}_r\cdot K^{[n]}_r\theta^{[n]}_r\,\d r
 =\frac4{N^2}\,\E_{\mu_{n,\tau}}[\mcl E_n] .
\end{split}
\end{equation}

\emph{Step 6: passage to the limit.}
As $n\to\infty$, $F_n\to F$ and $\mcl E_n\to\mcl E$ pointwise, by the uniform convergence of $x^{[n]}$ to $x$ and the continuity noted in Step~1, and these functions are bounded uniformly in $n$, by \eqref{e.observable.bounds} applied along the interpolated paths. By \eqref{e.moments.convergence}, applied to $F_n$, $F_n^2$ and $\mcl E_n$, both sides of \eqref{e.finite.estimate} converge to the corresponding sides of \eqref{e.conditional.variance.mu}.
\end{proof}

\begin{proof}[Proof of Proposition~\ref{p.posterior}]
Let $F:=\int_0^1a(r)L_{r,s}\,\d r$. We first show that $F$ is independent of $\si$. Let $S$ be a permutation of $\{1,\ldots,N\}$ that preserves each $I_{N,s}$, and let $O_S$ be as in Lemma~\ref{l.symmetries}. Consider the transformation
\begin{equation}\label{e.symmetry}
 (g,x,\si)\longmapsto(O_Sg,Sx,S\si)
\end{equation}
of the disorder, the path and the configuration. It preserves $\ga$, because $O_S$ is orthogonal and $S$ commutes with $\diag(q(0))$ and with $D(r)$; it preserves the counting measure on the slice; and it preserves the density \eqref{e.joint.density}, because $O_Sg\cdot v_{S\si}=g\cdot v_\si$, $Sx_1\cdot S\si=x_1\cdot\si$, and $U(r,Sx;O_Sg)=U(r,x;g)$ by Lemma~\ref{l.symmetries}(3), the weights $2^{-N}$ being invariant. Hence the law of $(g,X,\si)$ is invariant under \eqref{e.symmetry}. Since $\nabla U(r,Sx;O_Sg)=S\nabla U(r,x;g)$ and $P_s$ commutes with $S$, each $L_{r,s}$, and therefore $F$, is invariant under \eqref{e.symmetry}. It follows that $\E[\varphi(F)\1_{\{\si=\tau\}}]=\E[\varphi(F)\1_{\{\si=S\tau\}}]$ for every bounded measurable function $\varphi$, and that $\P(\si=\tau)=\P(\si=S\tau)$, for every $\tau\in\Omega_{N,k}$. These permutations act transitively on $\Omega_{N,k}$, so the conditional law of $F$ given $\si=\tau$ does not depend on $\tau$, that is, $F$ is independent of $\si$. In particular,
\begin{equation}\label{e.gauge}
 \Var(F)=\Var(F\mid\si=\tau)\qquad(\tau\in\Omega_{N,k}).
\end{equation}
Averaging \eqref{e.conditional.variance} over the law of $\si$ and using \eqref{e.gauge}, we obtain
\begin{equation*}
 \Var(F)\le\frac4{N^2}\int_0^1a(r)^2\,
 \E\bigl[P_s(\si-\bmm_r)\cdot K_rP_s(\si-\bmm_r)\bigr]\,\d r .
\end{equation*}
Conditionally on $\cG_r$, the matrix $K_r$ is fixed, and the residual $P_s(\si-\bmm_r)$ is centered with covariance $P_sC_rP_s$, whose only nonzero block is $C_{r,ss}$, by Lemma~\ref{l.response}. Hence $\E[P_s(\si-\bmm_r)\cdot K_rP_s(\si-\bmm_r)\mid\cG_r]=\tr(K_rP_sC_rP_s)=\tr(K_{r,ss}C_{r,ss})$, which proves \eqref{e.posterior.variance}.
\end{proof}

The symmetry step \eqref{e.gauge} is the only place in this section where the fact that we work with a slice $\Omega = \Omega_{N,k}$ is used. What it requires is a group of transformations of $(g,x,\si)$ that preserves the law of $(g,X,\si)$ and the function $F$, and acts transitively on $\Omega$. The permutations \eqref{e.symmetry} preserve the law whenever $\Omega$ and $w$ are invariant under them, which includes the full cube with the weights $2^{-N}e^{h^{(N)}\cdot\si}$, but their orbits on the full cube are the slices. On the full cube with an external field, the argument would therefore only show that $\E[F\mid\si]$ is a function of the magnetization of $\si$, and the term $\Var(\E[F\mid\si])$ in the decomposition $\Var(F)=\E[\Var(F\mid\si)]+\Var(\E[F\mid\si])$ is not controlled by the Brascamp--Lieb inequality, which bounds only the first term. Without external field, the spin flips $\si\mapsto\tau\si$, the coordinatewise products with $\tau\in\Sigma_N$, supply the missing transformations, as in \cite{aronow2026quantitative}. A flip preserves every species overlap, so there is an orthogonal map $O_\tau$ of the disorder space with $O_\tau v_\si=v_{\tau\si}$, as for $O_S$; the matrix $T_\tau:=\diag(\tau)$ commutes with $\diag(q(0))$ and $D(r)$, so $(g,x,\si)\mapsto(O_\tau g,T_\tau x,\tau\si)$ preserves $\ga$; and the proof of Lemma~\ref{l.symmetries}(3) applies verbatim, because the uniform weights are invariant under $\si\mapsto\tau\si$, and gives $U(r,T_\tau x;O_\tau g)=U(r,x;g)$ and $\nabla U(r,T_\tau x;O_\tau g)=T_\tau\nabla U(r,x;g)$. Hence the density \eqref{e.joint.density}, each $L_{r,s}$ and $F$ are invariant, and since the flips act transitively on $\Sigma_N$, the proof of Proposition~\ref{p.posterior} goes through on the full cube without field, with no other change; Lemma~\ref{l.conditional.variance} and its proof do not use the slice. The weights $2^{-N}e^{h^{(N)}\cdot\si}$ are not invariant under the flips when $h\neq0$: by Lemma~\ref{l.symmetries}(1), a flip maps the model with the field $h$ to the model with the site-dependent field $\tau_ih_s$ at $i\in I_{N,s}$. On a slice, the field is constant and plays no role in this section. 

\subsection{Proof of Theorem~\ref{t.remainder}}
\label{s.small}

We bound the two terms of \eqref{e.W.reduction} in turn. Two consequences of the response bound $K_r\succeq rC_r$ will be used for $r\ge\de$. Principal blocks preserve the order of symmetric matrices, so $K_{r,ss}\succeq\de C_{r,ss}$; since the trace of a product of two positive semidefinite matrices is nonnegative, this gives $\tr(C_{r,ss}^2)\le\de^{-1}\tr(K_{r,ss}C_{r,ss})\le\de^{-2}\tr(K_{r,ss}^2)$, and $\sum_s\tr(K_{r,ss}^2)\le\tr K_r^2$. Hence
\begin{equation}\label{e.block.bounds}
 \sum_s\tr(C_{r,ss}^2)\le\de^{-2}\tr K_r^2,
 \qquad
 \sum_s\tr(K_{r,ss}C_{r,ss})\le\de^{-1}\tr K_r^2
 \qquad(r\ge\de).
\end{equation}

\emph{Gibbs fluctuations of the overlap.}
For $r<\de$, we use $\sum_s\tr(C_{r,ss}^2)\le\sum_s(\tr C_{r,ss})^2\le N^2\sum_s\lambda_{N,s}^2\le N^2$. For $r\ge\de$ we use \eqref{e.block.bounds} and \eqref{e.energy.bound}. Thus
\begin{equation}\label{e.covariance.bound}
 \frac1{N^2}\int_0^1\E\sum_s\tr(C_{r,ss}^2)\,\d r
 \le\de+\frac1{2\kappa N\de^2}.
\end{equation}

\emph{From averages to individual levels.}
For $u\le r$, the variation of the realized gain across levels is controlled by the increments of the expected gain \eqref{e.gain.increment}: writing $L_{r,s}-L_{u,s}=N^{-1}P_s(\bmm_r-\bmm_u)\cdot P_s(2\si-\bmm_r-\bmm_u)$, the second factor has norm at most $4\sqrt{|I_{N,s}|}$, and \eqref{e.gain.increment} with $\lambda_{N,s}\le1$ gives
\begin{equation}\label{e.modulus}
 \sum_s\E|L_{r,s}-L_{u,s}|^2
 \le16\sum_s\lambda_{N,s}\bigl(p_{N,s}(r)-p_{N,s}(u)\bigr)
 \le16\bigl(\Gamma(r)-\Gamma(u)\bigr).
\end{equation}
We partition $[\de,1]$ into intervals $\Delta$ of length $\de$, and write $\bar L_{\Delta,s}:=\de^{-1}\int_\Delta L_{r,s}\,\d r$. The exact variance identity
\begin{equation}\label{e.orthogonal.decomposition}
 \int_\Delta\Var(L_{r,s})\,\d r
 =\de\Var(\bar L_{\Delta,s})
   +\frac1{2\de}\int_\Delta\int_\Delta
                  \Var(L_{r,s}-L_{u,s})\,\d r\,\d u ,
\end{equation}
which splits the variance over $\Delta$ into that of the average and that of the fluctuations around the average, together with \eqref{e.modulus}, bounds the sum over $\Delta$ and $s$ of the second terms by
\begin{equation*}
 8\de\sum_\Delta\bigl(\Gamma(\sup\Delta)-\Gamma(\inf\Delta)\bigr)\le8\de .
\end{equation*}
For the first terms, Proposition~\ref{p.posterior} with $a=\de^{-1}\1_\Delta$, summed over $\Delta$ and $s$, followed by \eqref{e.block.bounds} and \eqref{e.energy.bound}, gives a total of at most
\begin{equation*}
 \frac4{\de N^2}\int_\de^1
       \E\sum_s\tr(K_{r,ss}C_{r,ss})\,\d r
 \le\frac2{\kappa N\de^2}.
\end{equation*}
Finally, $L_{r,s}$ takes values in an interval of length $4\lambda_{N,s}$, so $\sum_s\Var(L_{r,s})\le4$ and the discarded levels contribute at most $4\de$. Altogether,
\begin{equation}\label{e.L.variance.bound}
 \int_0^1\sum_s\Var(L_{r,s})\,\d r
 \le12\de+\frac2{\kappa N\de^2}.
\end{equation}
Inserting \eqref{e.covariance.bound} and \eqref{e.L.variance.bound} into \eqref{e.W.reduction} gives \eqref{e.W.bound}.

For \eqref{e.small.defect}, recall from \eqref{e.def.remainder} and \eqref{e.defect.variance} that the left-hand side is $|\msc E_N(t,q)|\le\frac{M_\xi}2\msc W_N(t,q)$. If $\kappa N\ge1$, choose $\de=1/\lfloor(\kappa N)^{1/3}\rfloor$, so that $\de\le2(\kappa N)^{-1/3}$ and $\de^{-2}\le(\kappa N)^{2/3}$; the right-hand side of \eqref{e.W.bound} is then at most $(26+\frac52)(\kappa N)^{-1/3}$. If $\kappa N<1$, we use the trivial bound $\msc W_N\le\sum_s(2\lambda_{N,s})^2\le4$, which gives $|\msc E_N|\le2M_\xi\le2M_\xi(\kappa N)^{-1/3}$. \qed
%
%
%
%
\section{The scalar initial condition at fixed magnetization}
\label{s.scalar}

In this section, we study the scalar function $g_m$ appearing in the definition \eqref{e.def.slice.initial} of the slice initial condition $\psi_m$. We will see in Section~\ref{s.initial} that $g_m$ is the limit, as $n\to\infty$, of the free energy at $t=0$ of a one-species system with $n$ spins, when the configurations are constrained to have magnetization $m$. Three properties of $g_m$ are needed later, each for a different interpolation of its argument.

The first property is that the map $m \mapsto g_m(\nu) - I(m)$ is convex and Lipschitz. It is used in Section~\ref{s.initial} to identify the limit of the slice initial conditions. Since $I''\ge1$, the convexity makes $g_m$ strongly convex in $m$, which forces the minimizing slice to have a magnetization close to the prescribed one. The Lipschitz bound also gives the continuity of $f_m$ in $m$, which is used in Section~\ref{s.conclusion} to pass from the minimum over slices to the infimum over $m$. 

The second property is that, for every $T>0$ and $\nu_0,\nu_1\in\Prob([0,T])$, the map $[0,1]\ni\theta\mapsto g_m((1-\theta)\nu_0+\theta\nu_1)$ is concave. It comes from a representation of $g_m(\nu)$ as an infimum, over martingales started at $m$, of functionals that are affine in $\nu$; we deduce this representation from the variational formula of \cite{mourrat2025uninverting}, which itself is based on \cite{auffinger2015parisi}. In Section~\ref{s.duality}, we use this representation to pass from the Hopf formula to the martingale formula \eqref{e.slice.martingale}. 

Recall that, for a path $q$, we write $g_m(q)$ to denote the function $g_m$ evaluated at the law of $q(U)$, where $U$ is a uniform random variable on $[0,1]$. The last property we show is that, for every $q_0,q_1\in\cQ_\infty$, the map $[0,1]\ni\theta\mapsto g_m((1-\theta)q_0+\theta q_1)$ is convex. Since the path is the quantile function of the law it represents, the pointwise interpolation of paths is the displacement interpolation of laws, and not their mixture. The fact that the function $\psi_\circ(\cdot;0)=g_0$ satisfies this convexity property in the path variable has already been shown in \cite{chen2026ising, ho2026concavity}, and here we extend this result to every $m$. It makes the slice initial condition $\psi_m$ convex, which allows us to represent the solution of the Hamilton--Jacobi equation \eqref{e.slice.hj} by the Hopf formula \eqref{e.slice.hopf} in Section~\ref{s.hj}. As recalled in the introduction, the unconstrained initial condition $\psi_\circ(\cdot;h)$ is not convex in the path when the field $h$ is large; it is the constraint on the magnetization that restores this property.

\subsection{Paths, laws, and the Parisi equation}
\label{s.scalar.setup}

We first fix the correspondence between paths and probability measures. Recall that $\cQ$ is the set of right-continuous increasing maps from $[0,1)$ to $\R_+$. For $q\in\cQ_1$, we denote by $\nu_q$ the law of $q(U)$, where $U$ is uniform on $[0,1]$, so that $g_m(q)=g_m(\nu_q)$ with the convention recalled above. Conversely, for a probability measure $\nu$ on $\R_+$ with finite first moment, its quantile function $r\mapsto\inf\{v\ge0:\nu([0,v])>r\}$ is the unique element $q\in\cQ_1$ with $\nu_q=\nu$. Since the coupling $(q(U),q'(U))$ of $\nu_q$ and $\nu_{q'}$ is monotone, in the sense that both coordinates are increasing functions of the same variable~$U$ (see \cite[Section~2]{mourrat2023free}), the $L^1$ distance between two elements of $\cQ_1$ is the $1$-Wasserstein distance between their laws: $W_1(\nu_q,\nu_{q'})=\|q-q'\|_{L^1}$ for all $q,q'\in\cQ_1$. In particular, a functional of laws that is $1$-Lipschitz for $W_1$ defines a functional of paths that is $1$-Lipschitz for the $L^1$ norm, and conversely. We write $\Prob([0,T])$ for the set of probability measures on $[0,T]$, and $\Prob_1(\R_+)$ for the set of probability measures on $\R_+$ with finite first moment. For $\nu\in\Prob_1(\R_+)$, we write
\begin{equation*}
 \zeta_\nu(v)=\nu([0,v]),\qquad v\ge0,
\end{equation*}
for its cumulative distribution function. When $\nu=\nu_q$ for a path $q$, we have $\zeta_\nu(v)=|\{r\in[0,1):q(r)\le v\}|$, and for $\nu,\nu'\in\Prob([0,T])$ we have $W_1(\nu,\nu')=\int_0^T|\zeta_\nu-\zeta_{\nu'}|$. We say that a law $\nu\in\Prob([0,T])$ has finite support if it is a finite combination of Dirac masses, and that a path is a finite-step path if it takes finitely many values. Finite-step paths correspond to laws with finite support.

We now define $\Phi_\nu$ for every $\nu\in\Prob([0,T])$, using the well-posedness theory of \cite{jagtob}. A weak solution of \eqref{e.def.scalar.pde} is a continuous function $\Phi$ on $[0,T]\times\R$ with bounded weak derivative $\dr_x\Phi$, which satisfies the equation in the sense of distributions on $(0,T)\times\R$ and the terminal condition $\Phi(T,\cdot)=\log\cosh$. Let $\mu$ be the image of $\nu$ under $v\mapsto v/T$, so that $\mu([0,t])=\zeta_\nu(Tt)$, and let $\xi_T(t):=Tt^2$. Then $\Phi$ is a weak solution of \eqref{e.def.scalar.pde} if and only if $u(t,x):=\Phi(Tt,x)$ is a weak solution of
\begin{equation}\label{e.parisi.normalization}
 -\dr_tu=\frac{\xi_T''(t)}2\Bigl(\dr_x^2u+\mu([0,t])(\dr_xu)^2\Bigr)
 \quad\text{on }[0,1]\times\R,\qquad u(1,x)=\log\cosh x,
\end{equation}
which is the Parisi equation of \cite{jagtob,mourrat2025uninverting} for the covariance function $\xi_T$ and the measure $\mu\in\Prob([0,1])$. Since the results of these papers are valid for every measure and every covariance function with nonnegative coefficients, they apply to \eqref{e.def.scalar.pde} after this change of time. In particular, we have the following.

\begin{lemma}[Well-posedness and regularity]\label{l.scalar.wellposed}
Let $\nu\in\Prob([0,T])$. The equation \eqref{e.def.scalar.pde} has a unique weak solution, which we denote by $\Phi_\nu$. For every $j\ge1$, the derivative $\dr_x^j\Phi_\nu$ is continuous and bounded on $[0,T]\times\R$, and Lipschitz in $v$ uniformly in $x$. The equation obtained by differentiating \eqref{e.def.scalar.pde} $j\ge0$ times in $x$ holds at every point of $(0,T)\times\R$, with one-sided derivatives in $v$ and the corresponding one-sided limits of $\zeta_\nu$; in particular, it holds classically at every $v$ at which $\zeta_\nu$ is continuous. Moreover,
\begin{equation}\label{e.scalar.gradient.range}
 |\dr_x\Phi_\nu|<1,\qquad
 0<\dr_x^2\Phi_\nu\le1-(\dr_x\Phi_\nu)^2\le1.
\end{equation}
Finally, for every $j\ge0$, It\^o's formula applies to $\dr_x^j\Phi_\nu(v,Y_v)$ along every It\^o process $Y$ with bounded progressively measurable drift and diffusion coefficients, the time derivative being understood in the one-sided sense above.
\end{lemma}

\begin{proof}
The existence, uniqueness and regularity statements, including the pointwise form of the differentiated equations, are \cite[Theorem~4 and Corollary~11]{jagtob}, and It\^o's formula at this level of regularity is \cite[Proposition~22]{jagtob}, see also the remark following it. The bounds $|\dr_x\Phi_\nu|<1$ and $0<\dr_x^2\Phi_\nu\le1$ are \cite[Lemma~16]{jagtob}; the sharper upper bound in \eqref{e.scalar.gradient.range} is proved in Lemma~\ref{l.scalar.ito} below, by the same argument.
\end{proof}

For a finitely supported law, the solution is explicit. For $r\in[0,1]$ and $a\ge0$, we define the operators
\begin{equation}\label{e.def.gaussian.operator}
 T_{r,a}f(x):=\frac1r\log\E\exp\bigl(rf(x+\sqrt{2a}Z)\bigr)
 \quad(r>0),\qquad
 T_{0,a}f(x):=\E f(x+\sqrt{2a}Z),
\end{equation}
where $Z$ is a standard Gaussian random variable. These operators act on functions with at most linear growth. By H\"older's inequality, each $T_{r,a}$ is monotone, commutes with the addition of constants, and is a contraction for the uniform norm; moreover $T_{r,a}T_{r,a'}=T_{r,a+a'}$. If $f$ is smooth with bounded derivatives of all orders, then so is $T_{r,a}f$, and $(r,a,x)\mapsto T_{r,a}f(x)$ is smooth on $[0,1]\times[0,\infty)\times\R$: for $r>0$, $T_{r,a}f(x)=K(r)/r$, where $K(r):=\log\E\exp(rf(x+\sqrt{2a}Z))$ is smooth in $(r,a,x)$ with $K(0)=0$, so that $K(r)/r=\int_0^1K'(\theta r)\,\d\theta$ is smooth up to $r=0$.

Let $\nu\in\Prob([0,T])$ have finite support, and let $0\le v_0<v_1<\cdots<v_\ell\le T$ be the atoms of $\nu$. The function $\zeta_\nu$ is then constant on each interval $[v_i,v_{i+1})$, with value $\zeta_\nu(v_i)$, and it equals $0$ on $[0,v_0)$ and $1$ on $[v_\ell,T]$. Define $\Phi:[0,T]\times\R\to\R$ by backward recursion: set $\Phi(T,x):=\log\cosh x$, and, on each interval on which $\zeta_\nu$ is constant with value $r$,
\begin{equation}\label{e.def.scalar.recursion}
 \Phi(v,\cdot):=T_{r,v'-v}\Phi(v',\cdot),
 \qquad v\le v'\text{ in this interval}.
\end{equation}
The semigroup property of the operators shows that this is consistent. By construction, $\Phi$ is continuous, with $|\dr_x\Phi|\le1$, it is smooth on each strip $[v_i,v_{i+1}]\times\R$ with bounded derivatives of all orders in $x$, and it solves there the Parisi equation \eqref{e.def.scalar.pde}, which we recall:
\begin{equation}\label{e.scalar.pde}
 -\dr_v\Phi_\nu=\dr_x^2\Phi_\nu+\zeta_\nu(v)(\dr_x\Phi_\nu)^2.
\end{equation}
Indeed, if $r>0$, the function $e^{r\Phi}$ solves the backward heat equation $-\dr_vh=\dr_x^2h$ on the strip, and if $r=0$ this is true of $\Phi$ itself. Integrating by parts on each strip, the boundary terms cancel by the continuity of $\Phi$, so $\Phi$ is a weak solution, and $\Phi=\Phi_\nu$ by Lemma~\ref{l.scalar.wellposed}. On $[v_\ell,T]$, where $\zeta_\nu=1$, we have $T_{1,a}\log\cosh=\log\cosh+a$, so that
\begin{equation}\label{e.scalar.top}
 \Phi_\nu(v,x)=\log\cosh x+T-v,\qquad v\in[v_\ell,T].
\end{equation}
This identity holds for every $\nu\in\Prob([0,T])$ on $[T_0,T]\times\R$, where $T_0:=\sup\supp\nu$: the function equal to $\Phi^{T_0}_\nu(v,x)+T-T_0$ for $v\le T_0$, where $\Phi^{T_0}_\nu$ is the solution with $T_0$ in place of $T$, and to $\log\cosh x+T-v$ for $v\ge T_0$ is a weak solution of \eqref{e.def.scalar.pde}, because $\zeta_\nu=1$ on $[T_0,T]$ and $(\log\cosh)''+((\log\cosh)')^2=1$, so it is $\Phi_\nu$ by uniqueness. In particular, the quantity $T-\Phi_\nu(0,x)$ does not depend on the choice of $T\ge T_0$, and we may always enlarge~$T$.

The regularity of Lemma~\ref{l.scalar.wellposed} allows us to use the Parisi equation through It\^o's formula, along the diffusion whose drift is the magnetization $\dr_x\Phi_\nu$. Let $\nu\in\Prob([0,T])$, and write $\Phi:=\Phi_\nu$ and $\zeta:=\zeta_\nu$. The drift $2\zeta(v)\dr_x\Phi(v,x)$ is bounded, measurable in $v$, and $2$-Lipschitz in $x$ by \eqref{e.scalar.gradient.range}, so that for every $s\in[0,T]$ and $x\in\R$, the stochastic differential equation
\begin{equation}\label{e.scalar.diffusion}
 \d X_v=2\zeta(v)\,\dr_x\Phi(v,X_v)\,\d v+\sqrt2\,\d B_v,
 \qquad v\in[s,T],
\end{equation}
has a unique strong solution started from $x$ at time $s$, which we denote by $X^{s,x}$; here $B$ is a standard Brownian motion, with its completed natural filtration. We write $X^y:=X^{0,y}$. This is the diffusion of \cite{auffinger2015parisi}, in the present normalization.

\begin{lemma}[It\^o identities along the field diffusion]\label{l.scalar.ito}
Let $\nu\in\Prob([0,T])$, $s\in[0,T]$, $x\in\R$, and $X:=X^{s,x}$. For every $v\in[s,T]$,
\begin{equation}\label{e.magnetization.martingale}
 \dr_x\Phi(v,X_v)=\dr_x\Phi(s,x)+\sqrt2\int_s^v\dr_x^2\Phi(w,X_w)\,\d B_w,
\end{equation}
so that $M_v:=\dr_x\Phi(v,X_v)$ is a martingale with $M_T=\tanh X_T$. Moreover,
\begin{gather}
 \dr_x\Phi(s,x)=\E\bigl[\tanh X_T\bigr],\notag
 \\
 \dr_x^2\Phi(s,x)=1-\dr_x\Phi(s,x)^2
 -2\int_s^T\bigl(1-\zeta(w)\bigr)\E\bigl[\dr_x^2\Phi(w,X_w)^2\bigr]\,\d w,
 \label{e.gradient.representation}\\
 \Phi(s,x)=\E\bigl[\log\cosh X_T\bigr]
 -\int_s^T\zeta(w)\,\E\bigl[M_w^2\bigr]\,\d w.
 \label{e.scalar.optimal.value}
\end{gather}
In particular, \eqref{e.scalar.gradient.range} holds, $\Phi(v,\cdot)$ is strictly convex, $\dr_x\Phi$ is $1$-Lipschitz in $x$, and $M$ takes values in $(-1,1)$.
\end{lemma}

\begin{proof}
By Lemma~\ref{l.scalar.wellposed}, It\^o's formula applies to $\dr_x^j\Phi(w,X_w)$ for $j\le2$, and the equations obtained by differentiating \eqref{e.scalar.pde} $j$ times in $x$ hold at every $w$ at which $\zeta$ is continuous, hence for almost every $w$. The drift of $\Phi(w,X_w)$ is $\dr_w\Phi+\dr_x^2\Phi+2\zeta(\dr_x\Phi)^2=\zeta(\dr_x\Phi)^2$, that of $\dr_x\Phi(w,X_w)$ is $\dr_w\dr_x\Phi+\dr_x^3\Phi+2\zeta\,\dr_x\Phi\,\dr_x^2\Phi=0$, and that of $\dr_x^2\Phi(w,X_w)$ is $\dr_w\dr_x^2\Phi+\dr_x^4\Phi+2\zeta\,\dr_x\Phi\,\dr_x^3\Phi=-2\zeta(\dr_x^2\Phi)^2$, all evaluated at $(w,X_w)$. The stochastic integrals have bounded integrands, hence zero mean. The second identity gives \eqref{e.magnetization.martingale}, and the first gives \eqref{e.scalar.optimal.value}, since $\Phi(T,\cdot)=\log\cosh$. Taking expectations in \eqref{e.magnetization.martingale} at $v=T$ gives $\dr_x\Phi(s,x)=\E[\tanh X_T]$, so $|\dr_x\Phi|<1$, and It\^o's isometry gives
\begin{equation*}
 \E\bigl[\tanh^2X_T\bigr]=\dr_x\Phi(s,x)^2+2\int_s^T\E\bigl[\dr_x^2\Phi(w,X_w)^2\bigr]\,\d w.
\end{equation*}
The third identity gives $\dr_x^2\Phi(s,x)=\E[1-\tanh^2X_T]+2\int_s^T\zeta(w)\E[\dr_x^2\Phi(w,X_w)^2]\,\d w$, which is positive, and substituting the previous display yields  \eqref{e.gradient.representation}, hence $0 < \dr_x^2\Phi\le1-(\dr_x\Phi)^2$. The remaining assertions follow.
\end{proof}

\begin{lemma}[Continuity in the law]\label{l.scalar.continuity}
For all $\nu,\nu'\in\Prob([0,T])$ and $v\in[0,T]$,
\begin{equation}\label{e.scalar.continuity}
\begin{gathered}
 \sup_{x\in\R}|\Phi_\nu(v,x)-\Phi_{\nu'}(v,x)|
 \le\int_v^T|\zeta_\nu-\zeta_{\nu'}|,\\
 \sup_{x\in\R}|\dr_x\Phi_\nu(v,x)-\dr_x\Phi_{\nu'}(v,x)|
 \le2\Bigl(\int_0^T|\zeta_\nu-\zeta_{\nu'}|\Bigr)^{1/2}.
\end{gathered}
\end{equation}
In particular, if $\zeta_{\nu_n}\to\zeta_\nu$ in $L^1([0,T])$, then $\Phi_{\nu_n}\to\Phi_\nu$ and $\dr_x\Phi_{\nu_n}\to\dr_x\Phi_\nu$ uniformly on $[0,T]\times\R$.
\end{lemma}

\begin{proof}
The difference $w:=\Phi_\nu-\Phi_{\nu'}$ vanishes at $v=T$ and solves, for almost every $v$,
\begin{equation*}
 -\dr_vw=\dr_x^2w+\zeta_\nu(\dr_x\Phi_\nu+\dr_x\Phi_{\nu'})\dr_xw
 +(\zeta_\nu-\zeta_{\nu'})(\dr_x\Phi_{\nu'})^2.
\end{equation*}
Fix $(v,x)$, and let $Y$ solve $\d Y_u=\zeta_\nu(u)(\dr_x\Phi_\nu+\dr_x\Phi_{\nu'})(u,Y_u)\,\d u+\sqrt2\,\d B_u$ on $[v,T]$ with $Y_v=x$; the drift is bounded, measurable in $u$ and Lipschitz in $x$, by \eqref{e.scalar.gradient.range}, so this equation has a unique strong solution. By It\^o's formula (Lemma~\ref{l.scalar.wellposed}), the drift of $w(u,Y_u)$ is $-(\zeta_\nu-\zeta_{\nu'})(u)(\dr_x\Phi_{\nu'})^2(u,Y_u)$, and since $w(T,\cdot)=0$,
\begin{equation*}
 w(v,x)=\E\int_v^T\bigl(\zeta_\nu(u)-\zeta_{\nu'}(u)\bigr)(\dr_x\Phi_{\nu'})^2(u,Y_u)\,\d u.
\end{equation*}
Since $|\dr_x\Phi_{\nu'}|\le1$, this gives the first bound. For the second, we use the elementary interpolation inequality $\|f'\|_\infty^2\le4\|f\|_\infty\|f''\|_\infty$, valid for every bounded $f\in C^2(\R)$ with bounded second derivative; it follows from the Taylor bound $h|f'(x)|\le2\|f\|_\infty +\frac{h^2}2\|f''\|_\infty$, optimized in $h>0$. We apply it to $f=w(v,\cdot)$, whose second derivative is bounded by $1$ in absolute value, by \eqref{e.scalar.gradient.range}. The last assertion is immediate.
\end{proof}

The derivative $\dr_x\Phi_\nu$ is the magnetization of a single spin in the field $x$ at the covariance~$v$. The following two-sided bound compares it with the magnetization of a free spin, in shifted fields. It will give the behavior of $\dr_x\Phi_\nu$ at infinity, and the Lipschitz dependence of $g_m$ on $m$.

\begin{lemma}[Barriers for the magnetization]\label{l.barriers}
For every $\nu\in\Prob([0,T])$, $v\in[0,T]$ and $x\in\R$,
\begin{equation}\label{e.barriers}
 \tanh\bigl(x-2b_\nu(v)\bigr)\le\dr_x\Phi_\nu(v,x)\le\tanh\bigl(x+2b_\nu(v)\bigr),
 \qquad b_\nu(v):=\int_v^T\bigl(1-\zeta_\nu(u)\bigr)\,\d u.
\end{equation}
In particular, $\dr_x\Phi_\nu(v,x)\to\pm1$ as $x\to\pm\infty$, for every $v\in[0,T]$, so that $\dr_x\Phi_\nu(v,\cdot)$ is an increasing bijection from $\R$ onto $(-1,1)$. Moreover, $b_\nu(0)=\int u\,\nu(\d u)$.
\end{lemma}

\begin{proof}
Let $u:=\dr_x\Phi_\nu$, which is continuous, with $\dr_xu$ and $\dr_x^2u$ continuous and bounded, and which solves, for almost every $v$,
\begin{equation}\label{e.magnetization.pde}
 -\dr_vu=\dr_x^2u+2\zeta_\nu\,u\,\dr_xu,
\end{equation}
with $u(T,x)=\tanh x$, by Lemma~\ref{l.scalar.wellposed}. Let $u_\pm(v,x):=\tanh(x\pm2b_\nu(v))$, which is smooth in $x$ and Lipschitz in $v$, with $b_\nu'=-(1-\zeta_\nu)$ almost everywhere. Writing $\theta:=1-u_\pm^2$, we have $\dr_xu_\pm=\theta$, $\dr_x^2u_\pm=-2u_\pm\theta$ and $\dr_vu_\pm=\mp2(1-\zeta_\nu)\theta$, so that
\begin{equation*}
 -\dr_vu_\pm-\dr_x^2u_\pm-2\zeta_\nu u_\pm\dr_xu_\pm
 =\pm2(1-\zeta_\nu)\theta+2(1-\zeta_\nu)u_\pm\theta
 =\pm2(1-\zeta_\nu)(1\pm u_\pm)\theta,
\end{equation*}
which is nonnegative for $u_+$ and nonpositive for $u_-$. Hence the difference $d:=u_+-u$ satisfies, for almost every $v$,
\begin{equation*}
 -\dr_vd=\dr_x^2d+\zeta_\nu(u_++u)\,\dr_xd+\zeta_\nu\,\dr_x(u_++u)\,d+S,
 \qquad S\ge0,
\end{equation*}
with $d(T,\cdot)=0$ and $|d|\le2$, where we used $2(u_+\dr_xu_+-u\dr_xu)=\dr_x[(u_+-u)(u_++u)]$. Fix $(v,x)$, let $Y$ solve $\d Y_w=\zeta_\nu(w)(u_++u)(w,Y_w)\,\d w+\sqrt2\,\d B_w$ on $[v,T]$ with $Y_v=x$, whose drift is bounded and Lipschitz in $x$, and let $c:=\zeta_\nu\dr_x(u_++u)$, which is bounded by $2$. By It\^o's formula, which applies to $u$ by Lemma~\ref{l.scalar.wellposed} and to $u_+$ by its regularity, the drift of $e^{\int_v^wc(r,Y_r)\,\d r}d(w,Y_w)$ is $-e^{\int_v^wc(r,Y_r)\,\d r}S(w,Y_w)\le0$. Hence 
\begin{equation*}
d(v,x)\ge\E\bigl[e^{\int_v^Tc(r,Y_r)\,\d r}d(T,Y_T)\bigr]=0,
\end{equation*} 
which is the upper bound in \eqref{e.barriers}; the lower bound is proved in the same way with $u_-$. The bijection statement follows from \eqref{e.scalar.gradient.range}, and the last identity is Fubini's theorem, since $1-\zeta_\nu(u)=\nu((u,T])$.
\end{proof}

We now define the constrained initial condition. For $\nu\in\Prob([0,T])$, $y\in\R$ and $m\in[-1,1]$, recall that
\begin{equation}\label{e.def.g.m}
 \psi_\circ(\nu;y)=T-\Phi_\nu(0,y),\qquad
 g_m(\nu)=\sup_{y\in\R}\{\psi_\circ(\nu;y)+ym\}
 =T+\sup_{y\in\R}\{ym-\Phi_\nu(0,y)\}.
\end{equation}
Neither quantity depends on $T$, by \eqref{e.scalar.top}. In words, with the sign convention of \eqref{e.def.free.energy}, $\psi_\circ(\nu;y)$ is the free energy of a single spin in the field $y$ and in the random hierarchical field encoded by $\nu$, and $g_m(\nu)$ is the convex conjugate of $y\mapsto-\psi_\circ(\nu;y)$. The function $y\mapsto ym-\Phi_\nu(0,y)$ is strictly concave, by Lemma~\ref{l.scalar.ito}, and for $|m|<1$ it tends to $-\infty$ at $\pm\infty$, by Lemma~\ref{l.barriers}. For $|m|<1$, we therefore denote by
\begin{equation}\label{e.def.y.m}
 y(\nu,m)\in\R\qquad\text{the unique point with }\dr_x\Phi_\nu(0,y(\nu,m))=m,
\end{equation}
which is the unique maximizer in \eqref{e.def.g.m}. We will see in Proposition~\ref{p.scalar.martingale} that $g_{\pm1}(\nu)=\log2+\int v\,\nu(\d v)$.

\subsection{The one-spin branch in covariance time}
\label{s.scalar.dictionary}

We now relate the Parisi equation \eqref{e.scalar.pde}, which is parametrized by the covariance $v$, to the branch equation \eqref{e.branch.pde} for a single spin at $t=0$, which is parametrized by the level $r$. This dictionary connects the scalar initial condition $\psi_\circ$ to the free energies of Sections~\ref{s.branch} to~\ref{s.concentration}, and it is used in Section~\ref{s.initial.factorization}. The two parametrizations discussed here are two instances of a more general freedom in the choice of the time variable; see for instance the decompositions of paths in \cite{chen2024simultaneous}.

Let $q:[0,1]\to\R_+$ be a regular path with a single species, so that $q$ is continuous, piecewise affine and strictly increasing with $q(0)>0$, and let $\nu:=\nu_q$ be the law of $q(U)$, which is supported in $[q(0),q(1)]$. Its cumulative function is
\begin{equation*}
 \zeta_\nu(v)=\bigl|\{r\in[0,1):q(r)\le v\}\bigr|
 =\begin{cases}
   0&\text{if }0\le v<q(0),\\
   q^{-1}(v)&\text{if }q(0)\le v\le q(1),
  \end{cases}
\end{equation*}
which is continuous and piecewise affine on $[0,q(1)]$. Let $u$ be the solution of the one-spin branch equation
\begin{equation}\label{e.one.spin.branch}
 \dr_ru+\dot q(r)\,\dr_x^2u+r\,\dot q(r)\,(\dr_xu)^2=0
 \quad\text{on }[0,1)\times\R,
 \qquad u(1,x)=\log\cosh x,
\end{equation}
given by Lemma~\ref{l.wellposed} with $N=1$, $t=0$, $\Omega=\{-1,1\}$ and $w\equiv\frac12$, for which the terminal function \eqref{e.def.terminal} is $\log\frac12(e^x+e^{-x})=\log\cosh x$.

\begin{lemma}[Dictionary between the two parametrizations]\label{l.dictionary}
Let $q$, $\nu$ and $u$ be as above, let $T:=q(1)$, and let $Z$ be a standard Gaussian random variable. We have
\begin{equation}\label{e.def.dictionary}
 \Phi_\nu(v,x)=
 \begin{cases}
 \E\,u\bigl(0,x+\sqrt{2(q(0)-v)}\,Z\bigr)&\text{if }0\le v\le q(0),\\
 u\bigl(q^{-1}(v),x\bigr)&\text{if }q(0)\le v\le q(1).
 \end{cases}
\end{equation}
In particular, for every $h\in\R$,
\begin{equation}\label{e.dictionary.value}
 \psi_\circ(\nu_q;h)=q(1)-\E\,u\bigl(0,h+\sqrt{2q(0)}\,Z\bigr).
\end{equation}
\end{lemma}

\begin{proof}
Let $\Phi$ denote the right-hand side of \eqref{e.def.dictionary}. The two expressions agree at $v=q(0)$, so $\Phi$ is continuous on $[0,T]\times\R$, with $\Phi(T,\cdot)=u(1,\cdot)=\log\cosh$, and $|\dr_x\Phi|\le1$ since $|\dr_xu|\le1$. On $[q(0),q(1)]$, the chain rule gives $\dr_v\Phi(v,x)=\dr_ru(q^{-1}(v),x)/\dot q(q^{-1}(v))$, so that \eqref{e.one.spin.branch} reads $-\dr_v\Phi=\dr_x^2\Phi+q^{-1}(v)(\dr_x\Phi)^2=\dr_x^2\Phi+\zeta_\nu(v)(\dr_x\Phi)^2$. On $[0,q(0)]$, the function $\Phi$ solves the backward heat equation $-\dr_v\Phi=\dr_x^2\Phi$, which is \eqref{e.scalar.pde} with $\zeta_\nu=0$. Thus $\Phi$ solves \eqref{e.scalar.pde} classically on each of the finitely many strips delimited by $0$, $q(0)$, the images under $q$ of its breakpoints, and $q(1)$, with bounded derivatives in $x$ on each strip by Lemma~\ref{l.wellposed}. Integrating by parts on each strip, the boundary terms cancel by the continuity of $\Phi$, so $\Phi$ is a weak solution of \eqref{e.def.scalar.pde}, and $\Phi=\Phi_\nu$ by the uniqueness in Lemma~\ref{l.scalar.wellposed}. Finally, \eqref{e.dictionary.value} is the first case of \eqref{e.def.dictionary} at $v=0$, combined with the definition \eqref{e.def.g.m} of $\psi_\circ$ with $T=q(1)$.
\end{proof}

\subsection{A martingale representation, and Legendre duality}
\label{s.scalar.martingale}

We now represent $g_m$ as the value of a stochastic control problem, in which the control is a martingale with prescribed initial value $m$. Fix $T>0$, let $B$ be the standard Brownian motion of \eqref{e.scalar.diffusion}, with its completed natural filtration, and let $\cM_T$ denote the set of martingales $\al=(\al(v))_{v\in[0,T]}$ in this filtration with values in $[-1,1]$. As in the introduction, where the subscript is reserved for the species, we write the time variable of the elements of $\cM_T$ in parentheses, and accordingly $B(v):=B_v$. Since the initial $\sigma$-algebra is trivial, $\al(0)$ is deterministic. For $m\in[-1,1]$, let $\cM_T^m$ be the set of those $\al\in\cM_T$ with $\al(0)=m$. For $\al\in\cM_T$, we set
\begin{equation}\label{e.def.chi.c}
\begin{gathered}
 \chi_\al(x):=\int_0^x\E[\al(v\wedge T)^2]\,\d v\quad(x\ge0),\\
 c_T(\al):=\E\bigl[I(\al(T))-\sqrt2\,\al(T)B(T)\bigr]
                   +\int_0^T\bigl(1-\E[\al(v)^2]\bigr)\,\d v.
\end{gathered}
\end{equation}
We often call functions of the form $\chi_\al$ potentials. 
Since $v\mapsto\E[\al(v)^2]$ is increasing and takes values in $[\al(0)^2,1]$, the function $\chi_\al$ is convex, increasing, $1$-Lipschitz, vanishes at zero, and has slope $\E[\al(T)^2]$ on $[T,\infty)$. For a vector of martingales $\al=(\al_s)_{s\in\sS}$ as in the introduction, the cost there is $\cost_T(\al)=\sum_s\lambda_sc_{T_s}(\al_s)$; the potentials $\chi_{\al_s}$ are combined in the same way in Section~\ref{s.duality.transfer}. For $\nu\in\Prob([0,T])$, we write $\Phi:=\Phi_\nu$ and $\zeta:=\zeta_\nu$ until the end of the proof of Proposition~\ref{p.scalar.martingale}. For $y\in\R$, recall that $X^y$ solves \eqref{e.scalar.diffusion} from $X^y_0=y$, and set $M^y_v:=\dr_x\Phi(v,X^y_v)$, so that $M^y_0=\dr_x\Phi(0,y)$ and $M^y_T=\tanh X^y_T$. By Lemma~\ref{l.scalar.ito}, $M^y$ is a martingale with values in $(-1,1)$, so $M^y\in\cM_T$.

\begin{proposition}[Scalar martingale formula]\label{p.scalar.martingale}
Let $\nu\in\Prob([0,T])$.
\begin{enumerate}
\item For every $y\in\R$,
\begin{equation}\label{e.scalar.sup.martingale}
 \Phi(0,y)=\sup_{\al\in\cM_T}\Bigl\{y\,\al(0)+\E\bigl[\sqrt2\,\al(T)B(T)-I(\al(T))\bigr]
 +\int_0^T\zeta(v)\,\E[\al(v)^2]\,\d v\Bigr\},
\end{equation}
and the supremum is attained exactly at $\al=M^y$.
\item For every $m\in[-1,1]$,
\begin{equation}\label{e.scalar.martingale}
 g_m(\nu)=\inf_{\al\in\cM_T^m}
             \left\{c_T(\al)+\int\chi_\al\,\d\nu\right\}.
\end{equation}
For $|m|<1$, the infimum is attained at $M^{y(\nu,m)}$. For $m=\pm1$, it is attained at the constant martingale $\al\equiv m$, and
\begin{equation}\label{e.g.endpoints}
 g_{\pm1}(\nu)=\log2+\int v\,\nu(\d v).
\end{equation}
\end{enumerate}
\end{proposition}
\begin{proof}
\noindent\emph{Step 1: the unconstrained formula.}\par
In the normalization \eqref{e.parisi.normalization}, with the Brownian motion $W_t:=B_{Tt}/\sqrt T$ on $[0,1]$ and the external field $h=y$, \cite[Lemma~2.2]{mourrat2025uninverting} states that
\begin{equation*}
 \Phi(0,y)=\sup_{\tilde\al}\E\left[\frac12\int_0^1\xi_T''(t)\mu([0,t])\tilde\al_t^2\,\d t
 +\tilde\al_1\int_0^1\sqrt{\xi_T''(t)}\,\d W_t-\phi^*(\tilde\al_1)\right],
\end{equation*}
where the supremum is over the bounded martingales $\tilde\al$ on $[0,1]$ in the completed filtration of $W$, and where $\phi^*(z)=I(z)-yz$ for $z\in[-1,1]$ and $\phi^*(z)=+\infty$ otherwise is the convex conjugate of $x\mapsto\log\cosh(x+y)$; moreover, the supremum is attained exactly at $\tilde\al_t=\dr_x\Phi(Tt,\tilde X_t)$, where $\d\tilde X_t=\xi_T''(t)\mu([0,t])\dr_x\Phi(Tt,\tilde X_t)\,\d t+\sqrt{\xi_T''(t)}\,\d W_t$ and $\tilde X_0=y$. The map $\al\mapsto\tilde\al$, $\tilde\al_t:=\al(Tt)$, is a bijection from $\cM_T$ onto the set of bounded martingales $\tilde\al$ with $|\tilde\al_1|\le1$ almost surely, and the other bounded martingales have objective $-\infty$. Substituting $t=v/T$, the objective becomes $\int_0^T\zeta(v)\al(v)^2\,\d v+\sqrt2\,\al(T)B(T)-I(\al(T))+y\,\al(T)$, and $\E[\al(T)]=\al(0)$. This gives \eqref{e.scalar.sup.martingale}. Finally, $\tilde X_t=X^y_{Tt}$, so the maximizer is $\tilde\al_t=M^y_{Tt}$, that is, $\al=M^y$.

\smallskip
\noindent\emph{Step 2: the lower bound.}\par
Let $\al\in\cM_T$. By Fubini's theorem, since $\zeta(v)=\nu([0,v])$,
\begin{equation}\label{e.fubini.chi}
 \int_0^T\zeta(v)\E[\al(v)^2]\,\d v
 =\int\nu(\d w)\int_w^T\E[\al(v)^2]\,\d v
 =\chi_\al(T)-\int\chi_\al\,\d\nu.
\end{equation}
Since $T-\chi_\al(T)=\int_0^T(1-\E[\al(v)^2])\,\d v$, the objective in \eqref{e.scalar.sup.martingale} is $y\,\al(0)+T-c_T(\al)-\int\chi_\al\,\d\nu$, and \eqref{e.scalar.sup.martingale} reads
\begin{equation}\label{e.scalar.sup.rewritten}
 \Phi(0,y)=T+\sup_{\al\in\cM_T}\Bigl\{y\,\al(0)-c_T(\al)-\int\chi_\al\,\d\nu\Bigr\}.
\end{equation}
For $\al\in\cM_T^m$, this gives $c_T(\al)+\int\chi_\al\,\d\nu\ge T-\Phi(0,y)+ym=\psi_\circ(\nu;y)+ym$ for every $y\in\R$, and taking the supremum over $y$ in \eqref{e.def.g.m}, we obtain $c_T(\al)+\int\chi_\al\,\d\nu\ge g_m(\nu)$.

\smallskip
\noindent\emph{Step 3: attainment for $|m|<1$.}\par
Let $y:=y(\nu,m)$ and $\al:=M^y$, which belongs to $\cM_T^m$ by \eqref{e.def.y.m}. Since $\al$ attains the supremum in \eqref{e.scalar.sup.rewritten}, we have $c_T(\al)+\int\chi_\al\,\d\nu=T-\Phi(0,y)+ym=g_m(\nu)$, by the choice of $y$ in \eqref{e.def.g.m}.

\smallskip
\noindent\emph{Step 4: the endpoints.}\par
Let $m=1$. For the constant martingale $\al\equiv1$, we have $c_T(\al)=I(1)=\log2$ and $\chi_\al(x)=x$, so Step~2 gives $g_1(\nu)\le\log2+\int v\,\nu(\d v)$. For the reverse inequality, let $y\in\R$, $X:=X^y$, $M:=M^y$ and $m_y:=M_0$. By \eqref{e.scalar.optimal.value} and the identity $\log\cosh x=x-\log2+\log(1+e^{-2x})$,
\begin{equation*}
 \Phi(0,y)\le\E[X_T]-\log2+\E\bigl[\log(1+e^{-2X_T})\bigr]
 -\int_0^T\zeta(v)\E[M_v^2]\,\d v.
\end{equation*}
By \eqref{e.scalar.diffusion} and the martingale property, $\E[X_T]=y+2m_y\int_0^T\zeta$, and $2m_y-\E[M_v^2]\le2m_y-m_y^2\le1$ by Jensen's inequality, so that
\begin{equation*}
 \Phi(0,y)\le y-\log2+\int_0^T\zeta+\E\bigl[\log(1+e^{-2X_T})\bigr].
\end{equation*}
Since the drift in \eqref{e.scalar.diffusion} is at least $-2$, we have $X_T\ge y-2T+\sqrt2B_T$, so for $y\ge2T$ the last integrand is bounded by $\log(1+e^{-2\sqrt2B_T})$, which is integrable, and it tends to zero as $y\to\infty$. By dominated convergence, $\sup_y\{y-\Phi(0,y)\}\ge\log2-\int_0^T\zeta$, that is, $g_1(\nu)\ge T+\log2-\int_0^T\zeta=\log2+\int v\,\nu(\d v)$, since $\int_0^T(1-\zeta)=\int v\,\nu(\d v)$. This proves \eqref{e.g.endpoints} for $m=1$, and shows that the infimum is attained at $\al\equiv1$. The case $m=-1$ follows by the symmetry $x\mapsto-x$.
\end{proof}

We next collect the elementary properties of $g_m$ that follow from Legendre duality.

\begin{lemma}[Legendre duality in the magnetization]
\label{l.legendre.m}
Let $\nu\in\Prob([0,T])$. The function $m\mapsto g_m(\nu)$ is finite and continuous on $[-1,1]$, and the function
\begin{equation}\label{e.entropy.corrected}
 m\longmapsto g_m(\nu)-I(m)
\end{equation}
is even and convex on $[-1,1]$, with Lipschitz constant $2\int v\,\nu(\d v)$:
\begin{equation}\label{e.magnetization.lipschitz}
 \bigl|g_m(\nu)-g_{m'}(\nu)-I(m)+I(m')\bigr|
 \le2|m-m'|\int v\,\nu(\d v),\qquad m,m'\in[-1,1].
\end{equation}
In particular, for $|m_0|<1$ and $m\in[-1,1]$,
\begin{equation}\label{e.strong.convexity}
 g_m(\nu)\ge g_{m_0}(\nu)+y(\nu,m_0)(m-m_0)+\tfrac12(m-m_0)^2.
\end{equation}
For every $h\in\R$,
\begin{equation}\label{e.scalar.duality}
 \psi_\circ(\nu;h)=\inf_{m\in[-1,1]}\{g_m(\nu)-hm\},
\end{equation}
where the infimum is attained exactly at $m=\dr_x\Phi_\nu(0,h)\in(-1,1)$. Finally, $g_m(\delta_0)=I(m)$, and $|g_m(\nu)-g_m(\nu')|\le W_1(\nu,\nu')$ for all $\nu,\nu'\in\Prob([0,T])$ and $m\in[-1,1]$. Consequently, $g_m$ extends uniquely to a function on $\Prob_1(\R_+)$ that is $1$-Lipschitz for $W_1$.
\end{lemma}

\begin{proof}
\noindent\emph{Step 1: elementary properties.}\par
For $\nu=\delta_0$, we have $\Phi_\nu(0,x)=\log\cosh x+T$, so $m \mapsto g_m(\delta_0)$ is the convex dual of $\log \cosh$, which is $I$. The Lipschitz bound in $\nu$ stated at the very end follows from \eqref{e.scalar.continuity}, since $\int_0^T|\zeta_\nu-\zeta_{\nu'}|=W_1(\nu,\nu')$, and this bound passes to the supremum over $y$. Since $T$ is arbitrary and $g_m$ does not depend on it, the bound holds for all compactly supported laws, and these are dense in $\Prob_1(\R_+)$ for $W_1$: if $\nu^T$ is the image of $\nu\in\Prob_1(\R_+)$ under $x\mapsto x\wedge T$, then $W_1(\nu,\nu^T)=\int(x-T)_+\,\d\nu\to0$ as $T\to\infty$. The extension follows. The function $m\mapsto g_m(\nu)$ is finite on $[-1,1]$, by \eqref{e.g.endpoints} and \eqref{e.def.y.m}, and convex and lower semicontinuous, as a supremum of affine functions of $m$. It is also even, since $\Phi_\nu(0,\cdot)$ is even.

\smallskip
\noindent\emph{Step 2: convexity of \eqref{e.entropy.corrected}.}\par
The function $m\mapsto g_m(\nu)-T$ is the Legendre transform of the strictly convex function $\Phi_\nu(0,\cdot)$, whose derivative is a continuous increasing bijection onto $(-1,1)$ with inverse $y(\nu,\cdot)$, by Lemma~\ref{l.barriers}. Hence $m\mapsto g_m(\nu)$ is differentiable on $(-1,1)$, with
\begin{equation}\label{e.g.m.derivative}
 \dr_mg_m(\nu)=y(\nu,m),\qquad|m|<1.
\end{equation}
Since $\Phi_\nu(0,\cdot)$ is twice continuously differentiable with $\dr_x^2\Phi_\nu>0$, by Lemmas~\ref{l.scalar.wellposed} and~\ref{l.scalar.ito}, the function $y(\nu,\cdot)$ is continuously differentiable on $(-1,1)$ by the inverse function theorem, and differentiating the identity $\dr_x\Phi_\nu(0,y(\nu,m))=m$ and using \eqref{e.scalar.gradient.range}, we get
\begin{equation*}
 \dr_m^2g_m(\nu)=\dr_my(\nu,m)=\frac1{\dr_x^2\Phi_\nu(0,y(\nu,m))}
 \ge\frac1{1-m^2}=I''(m),\qquad|m|<1.
\end{equation*}
Thus \eqref{e.entropy.corrected} is convex on $(-1,1)$. At the endpoints, lower semicontinuity and convexity give
\begin{equation*}
 \liminf_{m\to1}g_m(\nu)\ge g_1(\nu),\qquad
 g_m(\nu)\le\tfrac{1-m}2\,g_{-1}(\nu)+\tfrac{1+m}2\,g_1(\nu),
\end{equation*}
so that $g_m(\nu)\to g_1(\nu)$ as $m\to1$; the endpoint $-1$ is similar. Hence $m\mapsto g_m(\nu)$ is continuous on $[-1,1]$, and the convexity of \eqref{e.entropy.corrected} extends to $[-1,1]$. Since $I''(m)=(1-m^2)^{-1}\ge1$, the function $m\mapsto g_m(\nu)-m^2/2$ is convex as well, and \eqref{e.strong.convexity} is its supporting inequality at $m_0$, where its derivative is $y(\nu,m_0)-m_0$ by \eqref{e.g.m.derivative}.

\smallskip
\noindent\emph{Step 3: the Lipschitz bound.}\par
Let $\varphi(m):=g_m(\nu)-I(m)$, so that $\varphi'(m)=y(\nu,m)-\operatorname{arctanh}m$ for $|m|<1$, by \eqref{e.g.m.derivative}. Evaluating the barriers \eqref{e.barriers} at $v=0$ and $x=y(\nu,m)$, where $\dr_x\Phi_\nu(0,x)=m$, and inverting the increasing function $\tanh$, we obtain
\begin{equation*}
 y(\nu,m)-2b_\nu(0)\le\operatorname{arctanh}m\le y(\nu,m)+2b_\nu(0),
 \qquad b_\nu(0)=\int v\,\nu(\d v),
\end{equation*}
that is, $|\varphi'|\le2\int v\,\nu(\d v)$ on $(-1,1)$. Since $\varphi$ is continuous on $[-1,1]$ by Step~2, \eqref{e.magnetization.lipschitz} follows.

\smallskip
\noindent\emph{Step 4: duality.}\par
Let $m_0:=\dr_x\Phi_\nu(0,h)$, so that $y(\nu,m_0)=h$ and $g_{m_0}(\nu)-hm_0=T-\Phi_\nu(0,h)$. For every other $m\in[-1,1]$, taking $y=h$ in \eqref{e.def.g.m} gives $g_m(\nu)-hm\ge T-\Phi_\nu(0,h)$, and the minimizer is unique by \eqref{e.strong.convexity}.
\end{proof}

\subsection{Concavity in the law, and supporting functions}
\label{s.scalar.support}

The right-hand side of \eqref{e.scalar.martingale} is an infimum of affine functions of $\nu$, so we can readily deduce the concavity of the map $\nu \mapsto g_m(\nu)$. To describe the supporting functions, for $|m|<1$, we denote by $M^{\nu,m}:=M^{y(\nu,m)}$ (defined with $\Phi=\Phi_\nu$) the minimizing martingale of Proposition~\ref{p.scalar.martingale}, and for $m=\pm1$ we let $M^{\nu,m}\equiv m$. We write $\chi_{\nu,m}:=\chi_{M^{\nu,m}}$ for the associated potential, which is convex, increasing, $1$-Lipschitz, and vanishes at zero. 

\begin{corollary}[Concavity in the law and supporting functions]
\label{c.scalar.support}
Let $m\in[-1,1]$. The function $\nu\mapsto g_m(\nu)$ is concave on $\Prob([0,T])$ for the affine interpolation of measures, and for all $\nu,\nu'\in\Prob([0,T])$,
\begin{equation}\label{e.scalar.support}
 g_m(\nu')\le g_m(\nu)+\int\chi_{\nu,m}\,\d(\nu'-\nu).
\end{equation}
Each $\chi_{\nu,m}$ is continuously differentiable on $\R_+$, with $\chi_{\nu,m}'(v)=\E[(M^{\nu,m}_{v\wedge T})^2]$. If $\zeta_{\nu_n}\to\zeta_\nu$ in $L^1([0,T])$, then $\chi_{\nu_n,m}'\to\chi_{\nu,m}'$ uniformly on $\R_+$; in particular, $\chi_{\nu_n,m}\to\chi_{\nu,m}$ uniformly on $[0,T]$. Consequently, for all $\nu,\nu'\in\Prob([0,T])$, with $\nu_\ep:=(1-\ep)\nu+\ep\nu'$,
\begin{equation}\label{e.scalar.derivative.law}
 \lim_{\ep\downarrow0}\frac{g_m(\nu_\ep)-g_m(\nu)}{\ep}
 =\int\chi_{\nu,m}\,\d(\nu'-\nu).
\end{equation}
\end{corollary}

\begin{proof}
Concavity is clear from \eqref{e.scalar.martingale}. The inequality \eqref{e.scalar.support} follows from applying the inequality in \eqref{e.scalar.martingale} at $\nu'$ with the martingale $\al=M^{\nu,m}$, and the equality at $\nu$ with the same martingale. The function $v\mapsto\E[(M^{\nu,m}_v)^2]$ is continuous on $[0,T]$, since $M^{\nu,m}$ is a bounded continuous process, which gives the differentiability of $\chi_{\nu,m}$. For the convergence statement, the case $m=\pm1$ is trivial. Let $|m|<1$, and let $y_n:=y(\nu_n,m)$ and $y:=y(\nu,m)$. Since $\dr_x\Phi_{\nu_n}(0,\cdot)\to\dr_x\Phi_\nu(0,\cdot)$ uniformly by Lemma~\ref{l.scalar.continuity}, we have $\dr_x\Phi_\nu(0,y_n)\to m=\dr_x\Phi_\nu(0,y)$, and $y_n\to y$ because $\dr_x\Phi_\nu(0,\cdot)$ is a continuous increasing bijection. Let $X^n$ solve \eqref{e.scalar.diffusion} with $\zeta_n:=\zeta_{\nu_n}$ and $\Phi_n:=\Phi_{\nu_n}$ in place of $\zeta$ and $\Phi$, driven by the same Brownian motion as $X:=X^{y}$ and started at $y_n$, and let $M^n_v:=\dr_x\Phi_n(v,X^n_v)$. Since the noise is additive, and $\dr_x\Phi_\nu$ is $1$-Lipschitz in $x$ with $|\dr_x\Phi_n|\le1$,
\begin{equation*}
 |X^n_v-X_v|\le|y_n-y|+2\int_0^T|\zeta_n-\zeta_\nu|
 +2T\sup|\dr_x\Phi_n-\dr_x\Phi_\nu|+2\int_0^v|X^n_w-X_w|\,\d w,
\end{equation*}
and Gronwall's inequality shows that $\sup_v|X^n_v-X_v|\to0$ for every realization of $B$, by Lemma~\ref{l.scalar.continuity}. Since $\dr_x\Phi_n\to\dr_x\Phi_\nu$ uniformly, we also get $\sup_v|M^n_v-M^{\nu,m}_v|\to0$ for every realization. Since these processes are bounded by $1$, bounded convergence gives $\sup_{v\le T}|\E[(M^n_v)^2]-\E[(M^{\nu,m}_v)^2]|\to0$, which is the uniform convergence of the derivatives of the potentials; the potentials themselves converge uniformly on $[0,T]$ since they vanish at $0$. Finally, the supporting inequalities at $\nu$ and at $\nu_\ep$ give
\begin{equation*}
 \int\chi_{\nu_\ep,m}\,\d(\nu'-\nu)
 \le\frac{g_m(\nu_\ep)-g_m(\nu)}{\ep}
 \le\int\chi_{\nu,m}\,\d(\nu'-\nu),
\end{equation*}
and $\zeta_{\nu_\ep}\to\zeta_\nu$ in $L^1$ as $\ep\to0$, so \eqref{e.scalar.derivative.law} follows from the continuity statement.
\end{proof}

\subsection{The susceptibility in magnetization coordinates}
\label{s.scalar.susceptibility}

We now prepare the proof of the path convexity of $g_m$. For illustration, we start with the simplest example, which shows the effect of fixing the magnetization rather than the field. For the constant path $q\equiv v$, that is $\nu=\delta_v$, the recursion \eqref{e.def.scalar.recursion} gives $\psi_\circ(\delta_v;h)=v-\E\log\cosh(h+\sqrt{2v}Z)$, and the heat equation gives
\begin{equation*}
 \frac{\d}{\d v}\psi_\circ(\delta_v;h)=\E\tanh^2\bigl(h+\sqrt{2v}Z\bigr),
 \qquad
 \left.\frac{\d^2}{\d v^2}\psi_\circ(\delta_v;h)\right|_{v=0}
 =2(1-\tanh^2h)(1-3\tanh^2h),
\end{equation*}
which is negative when $\tanh^2h>1/3$: the function $\psi_\circ(\cdot;h)$ is not convex in the path when $|h|$ is large. Let now $|m|<1$. The optimal field $y(\delta_v,m)$ in \eqref{e.def.g.m} solves $\E\tanh(y+\sqrt{2v}Z)=m$, so it equals $\operatorname{arctanh}m$ at $v=0$ and, by implicit differentiation, moves with speed $2m$ there. By the envelope theorem, $\frac{\d}{\d v}g_m(\delta_v)=\E\tanh^2(y(\delta_v,m)+\sqrt{2v}Z)$, and differentiating once more, the motion of the optimal field contributes a second, positive term:
\begin{equation*}
 \left.\frac{\d^2}{\d v^2}g_m(\delta_v)\right|_{v=0}
 =2(1-m^2)(1-3m^2)+2m(1-m^2)\cdot2m=2(1-m^2)^2>0.
\end{equation*}

For a finite-step path with distinct heights $q_0<\cdots<q_\ell$, carried by intervals of lengths $w_0,\ldots,w_\ell$, we will show in Proposition~\ref{p.path.convexity} that $\dr_{q_i}g_m(q)=w_i\E[M_{q_i}^2]$, where $M$ is the magnetization martingale of Lemma~\ref{l.scalar.ito} started at $m$, and convexity will follow from the following two observations. First, raising another height $q_j$, $j\ne i$, lowers the cumulative function $\zeta$ on the interval swept by the atom, which should make the magnetization at the height $q_i$ less spread, so the off-diagonal entries of the Hessian of $g_m$ should be nonpositive. Second, raising all the heights together amounts to starting the magnetization process from a more spread initial point, so the row sums of the Hessian of $g_m$ should be nonnegative. These two properties imply that the Hessian of $g_m$ is positive semidefinite. Lemma~\ref{l.diffusion.comparison} below makes these monotonicity properties precise, in coordinates that we now introduce.

Throughout this subsection and the next, $\nu\in\Prob([0,T))$ has finite support and $T>\max\supp\nu$; we write $\Phi:=\Phi_\nu$ and $\zeta:=\zeta_\nu$, so that $\Phi$ is given by the recursion \eqref{e.def.scalar.recursion}, and is smooth on each strip $[v_i,v_{i+1}]\times\R$ delimited by the atoms of $\nu$. By Lemma~\ref{l.barriers}, for every $v\in[0,T]$ the map $x\mapsto\dr_x\Phi(v,x)$ is an increasing bijection from $\R$ onto $(-1,1)$; we denote by $x(v,z)$ its inverse, which is smooth on each strip $[v_i,v_{i+1}]\times(-1,1)$ by the inverse function theorem, since $\dr_x^2\Phi>0$. We define
\begin{equation}\label{e.def.susceptibility}
 A(v,z):=\dr_x^2\Phi(v,x(v,z)),\qquad
 W(v,z):=\sup_{x\in\R}\{xz-\Phi(v,x)\},\qquad z\in(-1,1).
\end{equation}
Thus $A$ is the susceptibility $\dr_x^2\Phi$, viewed as a function of the magnetization $\dr_x\Phi$, and $W(v,\cdot)$ is the Legendre transform of $\Phi(v,\cdot)$, so that $g_m(\nu)=T+W(0,m)$. Both functions are smooth on each strip $[v_i,v_{i+1}]\times(-1,1)$. The envelope formula gives $\dr_zW(v,z)=x(v,z)$, hence $\dr_z^2W=1/A$, and $\dr_vW=-\dr_v\Phi(v,x(v,z))$, hence by \eqref{e.scalar.pde}
\begin{equation}\label{e.dual.pde}
 \dr_vW=A+\zeta(v)z^2,\qquad \dr_z^2W=\frac1A,\qquad W(T,z)=I(z),
\end{equation}
on each strip. Moreover, by \eqref{e.scalar.gradient.range},
\begin{equation}\label{e.susceptibility.bound}
 0<A(v,z)\le1-z^2,\qquad z\in(-1,1),
\end{equation}
so that $A$ extends continuously to $[0,T]\times[-1,1]$ by $A(v,\pm1):=0$.

The martingale $M$ of Lemma~\ref{l.scalar.ito} is the magnetization process of a single spin, and it is natural to describe its dynamics in its own coordinates. For $s\in[0,T]$ and $z\in(-1,1)$, we define the magnetization process started at time $s$ from $z$ by
\begin{equation}\label{e.def.magnetization.process}
 M^{s,z}_r:=\dr_x\Phi\bigl(r,X^{s,x(s,z)}_r\bigr),\qquad r\in[s,T],
\end{equation}
and for $z=\pm1$ we let $M^{s,\pm1}\equiv\pm1$. We write $\E_z$ for the expectation of functionals of $M^{s,z}$, the starting time $s$ being understood; by the Markov property of \eqref{e.scalar.diffusion}, the process $M^{s,z}$ is a Markov process, and its law depends only on $(s,z)$. For $s=0$, we have $x(0,z)=y(\nu,z)$ by \eqref{e.def.y.m}, so $M^{0,z}$ is the martingale $M^{\nu,z}$ of Proposition~\ref{p.scalar.martingale}. By Lemma~\ref{l.scalar.ito}, $M^{s,z}$ is a martingale started at $z$ with values in $(-1,1)$, and since $\dr_x^2\Phi(r,X_r)=A(r,\dr_x\Phi(r,X_r))$, the identity \eqref{e.magnetization.martingale} reads
\begin{equation}\label{e.magnetization.sde}
 \d M_r=\sqrt2\,A(r,M_r)\,\d B_r.
\end{equation}
In these coordinates, the law $\nu$ enters only through the diffusion coefficient $A$, and the initial condition is the prescribed magnetization. The two comparison statements below say that a larger cumulative function $\zeta$ gives a larger diffusion coefficient, and that a larger diffusion coefficient spreads the magnetization more. We first record the monotone dependence on the starting point. Let $z\le z'$ in $(-1,1)$, so that $x:=x(s,z)\le x':=x(s,z')$, and let $X$ and $X'$ solve \eqref{e.scalar.diffusion} from $(s,x)$ and $(s,x')$, driven by the same Brownian motion. Their difference $D:=X'-X$ is absolutely continuous, with $D'_v=2\zeta(v)[\dr_x\Phi(v,X'_v)-\dr_x\Phi(v,X_v)]$; since $\dr_x\Phi$ is increasing in $x$, the function $D$ is increasing on every interval on which it is positive, so $D_v\ge x'-x\ge0$ for every~$v$. Hence $M^{s,z}_r\le M^{s,z'}_r$ for every $r$, and since both are martingales, for every Lipschitz function $\varphi$ on $[-1,1]$ and every $v\in[s,T]$,
\begin{equation}\label{e.lipschitz.in.z}
 |\E_{z'}\varphi(M_v)-\E_z\varphi(M_v)|
 \le\|\varphi'\|_\infty\E\bigl[M^{s,z'}_v-M^{s,z}_v\bigr]=\|\varphi'\|_\infty(z'-z).
\end{equation}
This remains true when $z'=1$ or $z=-1$, since $\E|1-M^{s,z}_v|=1-z$ by the martingale property. In particular, $z\mapsto\E_z\varphi(M_v)$ is continuous on $[-1,1]$, uniformly in $s$.

\begin{lemma}[Diffusion comparisons]\label{l.diffusion.comparison}
Let $\nu_1,\nu_2\in\Prob([0,T))$ have finite support, with susceptibilities $A_1,A_2$ and magnetization processes $M^1,M^2$, and let $0\le s\le v\le T$. We recall that $\E_z$ denotes the expectation when the process starts from $z\in[-1,1]$, the starting time $s$ being understood.
\begin{enumerate}
\item For every continuous convex function $\varphi$ on $[-1,1]$, the map $z\mapsto\E_z[\varphi(M^1_v)]$ is convex on $[-1,1]$.
\item If $\zeta_{\nu_1}\le\zeta_{\nu_2}$ on $[0,T]$, then $A_1\le A_2$ on $[0,T]\times[-1,1]$.
\item If $A_1\le A_2$, then for every continuous convex function $\varphi$ on $[-1,1]$,
\begin{equation*}
 \E_z[\varphi(M_v^1)]\le\E_z[\varphi(M_v^2)],
 \qquad z\in[-1,1].
\end{equation*}
\end{enumerate}
\end{lemma}

\begin{proof}
\noindent\emph{Step 0: regularity.}\par
We drop the index $1$ in this step and the next. In items~(1) and~(3), we may assume that $\varphi$ is smooth on $[-1,1]$, by extending $\varphi$ outside $[-1,1]$ (e.g.\ via a dilation) and mollifying, which preserves convexity and gives uniform approximations; and the case $s=v$ is trivial, so let $s<v$. For $s'\in[s,v]$ and $x\in\R$, let $F(s',x):=\E[\varphi(\dr_x\Phi(v,X^{s',x}_v))]$. On each strip, the drift of \eqref{e.scalar.diffusion} is smooth in $x$ with bounded derivatives, so by the parabolic theory used in Lemma~\ref{l.wellposed} and It\^o's formula, $F$ is smooth in $x$ with bounded derivatives on each strip, continuous on $[s,v]\times\R$, and solves the backward Kolmogorov equation $-\dr_{s'}F=\dr_x^2F+2\zeta\,\dr_x\Phi\,\dr_xF$ classically on each strip. For $|z|<1$, we have $u(s',z):=\E_z[\varphi(M_v)]=F(s',x(s',z))$, which is therefore smooth on each strip $[v_i,v_{i+1}]\times(-1,1)$ intersected with $[s,v]\times(-1,1)$, with derivatives bounded on $[-1+\eta,1-\eta]$ for every $\eta>0$; the flow of \eqref{e.scalar.diffusion} and the inverse $x(\cdot,\cdot)$ are jointly continuous, so $u$ is continuous on $[s,v]\times(-1,1)$, and on $[s,v]\times[-1,1]$ by \eqref{e.lipschitz.in.z}. Writing $F(s',x)=u(s',\dr_x\Phi(s',x))$, the chain rule and \eqref{e.magnetization.pde} transform the equation for $F$ into
\begin{equation}\label{e.dual.heat}
 -\dr_{s'}u=a\,\dr_z^2u,\qquad a:=A^2,
\end{equation}
on each strip, in the interior.

\smallskip
\noindent\emph{Step 1: convexity in the starting point.}\par
Let $z\in(-1,1)$, $x:=x(s,z)$, $X:=X^{s,x}$ and $M_r:=\dr_x\Phi(r,X_r)$. Since the noise in \eqref{e.scalar.diffusion} is additive and the drift is smooth in $x$ on each strip, for every Brownian path the map $x'\mapsto X^{s,x'}_r$ is continuously differentiable, with derivative $J_r:=\exp(2\int_s^r\zeta(w)\dr_x^2\Phi(w,X_w)\,\d w)\in[1,e^{2T}]$ at $x'=x$. Differentiating under the expectation and using $\dr_zu=\dr_xF/\dr_x^2\Phi$, we obtain
\begin{equation*}
 \dr_zu(s,z)=\E\bigl[\varphi'(M_v)N_v\bigr],
 \qquad N_r:=\frac{\dr_x^2\Phi(r,X_r)J_r}{\dr_x^2\Phi(s,x)},\quad r\in[s,v].
\end{equation*}
By It\^o's formula, the drift of $\dr_x^2\Phi(r,X_r)$ is $-2\zeta(\dr_x^2\Phi)^2$, as in the proof of Lemma~\ref{l.scalar.ito}, while $\d J_r=2\zeta\,\dr_x^2\Phi(r,X_r)J_r\,\d r$, so that
\begin{equation*}
 \d N_r=\sqrt2\,N_r\,\frac{\dr_x^3\Phi}{\dr_x^2\Phi}(r,X_r)\,\d B_r.
\end{equation*}
Thus $N$ is a positive continuous local martingale, started from $1$ at time $s$ and bounded by $e^{2T}/\dr_x^2\Phi(s,x)$, hence a martingale with $\E[N_v]=1$. Let $\tilde\P:=N_v\cdot\P$ on the $\sigma$-algebra at time $v$. By Girsanov's theorem, $\tilde B_r:=B_r-\sqrt2\int_s^r(\dr_x^3\Phi/\dr_x^2\Phi)(w,X_w)\,\d w$ is a Brownian motion under $\tilde\P$ on $[s,v]$, and
\begin{equation*}
 \d X_r=b(r,X_r)\,\d r+\sqrt2\,\d\tilde B_r,
 \qquad b:=2\zeta\,\dr_x\Phi+2\,\dr_x\log\dr_x^2\Phi,
\end{equation*}
while $\dr_zu(s,z)=\tilde\E[\varphi'(\dr_x\Phi(v,X_v))]$. The function $b$ is smooth in $x$ on each strip, since $\dr_x^2\Phi>0$. For a Brownian motion $B'$ and $x'\in\R$, the equation $\d Y_r=b(r,Y_r)\,\d r+\sqrt2\,\d B'_r$ on $[s,v]$, started from $x'$ at time $s$, is, for every Brownian path, an ordinary differential equation for $Y-\sqrt2B'$ with a locally Lipschitz vector field: it has a unique solution up to its explosion time, which is a measurable functional of the Brownian path, and two solutions driven by the same Brownian motion from $x'<x''$ remain ordered as long as both exist. Under $\tilde\P$, the process $X$ is such a solution, driven by $\tilde B$ and started at $x$, and it does not explode before $v$, since it is defined on $[s,v]$ almost surely under $\P$, hence under $\tilde\P$. By uniqueness in law, the solution $Y^{x}$ started at $x$ and driven by any Brownian motion does not explode before $v$ either, and $\dr_zu(s,z)=\E[\varphi'(\dr_x\Phi(v,Y^{x}_v))]$. Now let $z<z'$ in $(-1,1)$, so that $x=x(s,z)<x':=x(s,z')$. Driving $Y^{x}$ and $Y^{x'}$ by the same Brownian motion, we get $Y^{x}_v<Y^{x'}_v$, and since $\varphi'$ and $\dr_x\Phi(v,\cdot)$ are increasing, $\dr_zu(s,z)\le\dr_zu(s,z')$. Hence $u(s,\cdot)$ is convex on $(-1,1)$, and on $[-1,1]$ by continuity. This proves item~(1).

\smallskip
\noindent\emph{Step 2: comparison of the coefficients.}\par
Adding atoms of mass zero, we may assume that $\nu_1$ and $\nu_2$ have the same atoms $v_0<\cdots<v_\ell$, so that $\zeta_1:=\zeta_{\nu_1}$ and $\zeta_2:=\zeta_{\nu_2}$ are constant on the same strips. For $\theta\in[0,1]$, let $\nu_\theta:=(1-\theta)\nu_1+\theta\nu_2$, whose cumulative function is $\zeta_\theta=\zeta_1+\theta\eta$ with $\eta:=\zeta_2-\zeta_1\ge0$, and let $\Phi_\theta$, $x_\theta$, $A_\theta$, $W_\theta$ and $M^\theta$ be the corresponding objects. On each strip, by the recursion \eqref{e.def.scalar.recursion}, $\Phi_\theta(v,x)$ is a finite composition of operators $T_{r,a}$, whose parameters $r=\zeta_\theta(v_i)$ are affine in $\theta$, applied to $\log\cosh$. By the smoothness of $(r,a,x)\mapsto T_{r,a}f(x)$ noted after \eqref{e.def.gaussian.operator}, and by differentiation under the expectation, $\Phi_\theta(v,x)$ is smooth in $(\theta,v,x)$ on $[0,1]\times[v_i,v_{i+1}]\times\R$, with derivatives in $x$ of all orders bounded. By the inverse function theorem, $x_\theta(v,z)$, $A_\theta(v,z)$ and $W_\theta(v,z)$ are then smooth in $(\theta,v,z)$ on $[0,1]\times[v_i,v_{i+1}]\times(-1,1)$. Let $H:=\dr_\theta W_\theta$, which is continuous on $[0,T]\times(-1,1)$ for each $\theta$. Differentiating \eqref{e.dual.pde} in $\theta$, and using $\dr_\theta(1/\dr_z^2W_\theta)=-A_\theta^2\dr_z^2H$, we find that on each strip,
\begin{equation}\label{e.linear.response.pde}
 -\dr_vH=A_\theta^2\,\dr_z^2H-\eta(v)z^2,\qquad H(T,z)=0,
\end{equation}
the terminal condition coming from $W_\theta(T,\cdot)=I$. Moreover, by Lemma~\ref{l.scalar.continuity}, $\sup_x|\Phi_\theta(v,x)-\Phi_{\theta'}(v,x)|\le|\theta-\theta'|\int_v^T\eta$, and the Legendre transform is $1$-Lipschitz for the uniform norm, so
\begin{equation}\label{e.linear.response.bound}
 |H(v,z)|\le\int_v^T\eta,\qquad (v,z)\in[0,T]\times(-1,1).
\end{equation}
Fix $\theta$, and for this step let $s\in[0,T]$ be arbitrary. Let $z\in(-1,1)$, let $M:=M^\theta$ be started at $(s,z)$, and for small $\ep>0$ let $\tau_\ep$ be its exit time from $[-1+\ep,1-\ep]$. On each strip, $H$ is smooth with bounded derivatives on this interval, so It\^o's formula on $[s,T\wedge\tau_\ep]$, together with \eqref{e.magnetization.sde} and \eqref{e.linear.response.pde}, gives
\begin{equation*}
 H(s,z)=\E\bigl[H(T\wedge\tau_\ep,M_{T\wedge\tau_\ep})\bigr]
 -\E\int_s^{T\wedge\tau_\ep}\eta(r)M_r^2\,\d r.
\end{equation*}
Since $M$ takes values in $(-1,1)$ and is continuous, $T\wedge\tau_\ep=T$ for every sufficiently small $\ep$, for every realization. Letting $\ep\to0$, by the continuity of $H$, the bound \eqref{e.linear.response.bound} and bounded convergence,
\begin{equation}\label{e.linear.response}
 H(s,z)=-\int_s^T\eta(r)\,\E_z\bigl[(M^\theta_r)^2\bigr]\,\d r.
\end{equation}
By item~(1) applied to $\nu_\theta$ and $\varphi(z)=z^2$, each function $z\mapsto\E_z[(M^\theta_r)^2]$ is convex, and $\eta\ge0$, so $H(s,\cdot)$ is concave on $(-1,1)$, that is, $\dr_z^2H\le0$. Hence $\dr_\theta A_\theta=\dr_\theta(1/\dr_z^2W_\theta)=-A_\theta^2\dr_z^2H\ge0$: the function $\theta\mapsto A_\theta(s,z)$ is increasing on $[0,1]$, which gives $A_1\le A_2$ on $[0,T]\times(-1,1)$. Both vanish at $z=\pm1$. This proves item~(2).

\smallskip
\noindent\emph{Step 3: comparison of the processes.}\par
Let $u_i(s',z):=\E_z[\varphi(M^i_v)]$ for the process started at time $s'$ from $z$, and $a_i:=A_i^2$. By Step~0 and item~(1), $u_1$ is continuous on $[s,v]\times[-1,1]$, smooth on each strip in the interior, convex in $z$, and solves $-\dr_{s'}u_1=a_1\dr_z^2u_1$. Let $|z|<1$ and $\ep>0$, and let $\tau_\ep$ be the exit time of $M^2$, started at time $s$ from $z$, from $[-1+\ep,1-\ep]$. By It\^o's formula on $[s,v\wedge\tau_\ep]$ and \eqref{e.magnetization.sde}, the drift of $u_1(r,M^2_r)$ is $(a_2-a_1)\dr_z^2u_1\ge0$, so that $u_1(s,z)\le\E[u_1(v\wedge\tau_\ep,M^2_{v\wedge\tau_\ep})]$. Letting $\ep\to0$ and using the continuity of $u_1$, we obtain $u_1(s,z)\le\E[u_1(v,M^2_v)]=u_2(s,z)$ by bounded convergence. The case $|z|=1$ is trivial.
\end{proof}
\subsection{Convexity in the path}
\label{s.scalar.convexity}

We now regard $g_m$ as a function of the path. For a finite-step path
\begin{equation}\label{e.step.path}
 q=\sum_{i=0}^{\ell}q_i\1_{[r_i,r_{i+1})},\qquad
 0=r_0<r_1<\cdots<r_{\ell+1}=1,\qquad 0\le q_0\le q_1\le\cdots\le q_\ell,
\end{equation}
the law $\nu_q$ has atoms $q_i$ with masses $w_i:=r_{i+1}-r_i$, and we write $g_m(q)=g_m(\nu_q)$. In this subsection, the partition $(r_i)$ is fixed and we regard $g_m(q)$ as a function of the heights $(q_0,\ldots,q_\ell)$; the derivative $\dr_{q_i}g_m(q)$ refers to this function. When the heights are distinct, $\zeta_{\nu_q}$ takes the value $r_{i+1}$ on $[q_i,q_{i+1})$. We write $M$ for the magnetization process \eqref{e.def.magnetization.process} associated with $\nu_q$ and started at time $0$, and $\E_z$ for the expectation when it starts from~$z$; in particular, $\E_m[\varphi(M_v)]=\E[\varphi(M^{\nu_q,m}_v)]$.

\begin{proposition}[Convexity at fixed magnetization]\label{p.path.convexity}
For every $m\in[-1,1]$, the function $g_m$ is convex on $\cQ_1$ and satisfies
\begin{equation}\label{e.scalar.lipschitz}
 |g_m(q)-g_m(q')|\le\|q-q'\|_{L^1},\qquad q,q'\in\cQ_1.
\end{equation}
For a finite-step path \eqref{e.step.path} with $0<q_0<\cdots<q_\ell$ and $|m|<1$, its derivative is
\begin{equation}\label{e.scalar.path.derivative}
 \dr_{q_i}g_m(q)=(r_{i+1}-r_i)\,\E_m\bigl[M_{q_i}^2\bigr],
 \qquad i\in\{0,\ldots,\ell\}.
\end{equation}
More generally, for every $a\in\cQ_\infty$, the path
\begin{equation}\label{e.def.z.a}
 z_a(r):=\E\bigl[(M^{\nu_a,m}_{a(r)})^2\bigr]=\chi_{\nu_a,m}'(a(r)),
 \qquad r\in[0,1),
\end{equation}
is increasing with values in $[m^2,1]$, and it supports $g_m$ at $a$: for every $b\in\cQ_\infty$,
\begin{equation}\label{e.scalar.supporting}
 g_m(b)\ge g_m(a)+\int_0^1z_a(r)\bigl(b(r)-a(r)\bigr)\,\d r,
\end{equation}
and the integral in \eqref{e.scalar.supporting} is the right derivative of $\theta\mapsto g_m(a+\theta(b-a))$ at $\theta=0$.
\end{proposition}

As explained at the beginning of Section~\ref{s.scalar.susceptibility}, the proof of convexity consists in showing that the Hessian of $g_m$ in the heights has nonpositive off-diagonal entries and nonnegative row sums. The sign of the off-diagonal entries comes from items~(2) and~(3) of Lemma~\ref{l.diffusion.comparison}, and the sign of the row sums from item~(1). For $m=\pm1$, we have $M\equiv m$ and $z_a\equiv1$, and \eqref{e.scalar.supporting} is an equality by \eqref{e.g.endpoints}.

\begin{proof}[Proof of Proposition~\ref{p.path.convexity}]
\noindent\emph{Step 1: smoothness, and directional derivatives.}\par
Let $|m|<1$ and $0<q_0<\cdots<q_\ell<T$. On the open cone of such heights, $\Phi_{\nu_q}(0,y)$ is a smooth function of $(q_0,\ldots,q_\ell,y)$, as a finite composition of the operators \eqref{e.def.gaussian.operator}, which depend smoothly on their parameter $a$; and $y(\nu_q,m)$ is a smooth function of the heights by the implicit function theorem, since $\dr_x^2\Phi_{\nu_q}>0$. Hence $g_m(q)=T-\Phi_{\nu_q}(0,y(\nu_q,m))+my(\nu_q,m)$ is smooth on this cone.

We now compute directional derivatives. Let $a\in\cQ_\infty$, and let $\ga$ be a bounded measurable function on $[0,1)$ such that $a_\theta:=a+\theta\ga$ belongs to $\cQ_\infty$ for all small $\theta>0$, with uniform norms bounded by some $T'$, and let $\nu_\theta:=\nu_{a_\theta}\in\Prob([0,T'])$. We claim that
\begin{equation}\label{e.directional.derivative}
 \lim_{\theta\downarrow0}\frac{g_m(a_\theta)-g_m(a)}{\theta}
 =\int_0^1\chi'_{\nu_0,m}(a(r))\,\ga(r)\,\d r.
\end{equation}
For every potential $\chi$, we have $\int\chi\,\d(\nu_\theta-\nu_0)=\int_0^1[\chi(a_\theta(r))-\chi(a(r))]\,\d r$, so the supporting inequality \eqref{e.scalar.support}, applied at $\nu_\theta$ and at $\nu_0$, gives
\begin{equation*}
 \int_0^1\bigl[\chi_{\nu_\theta,m}(a_\theta(r))-\chi_{\nu_\theta,m}(a(r))\bigr]\,\d r
 \le g_m(a_\theta)-g_m(a)
 \le\int_0^1\bigl[\chi_{\nu_0,m}(a_\theta(r))-\chi_{\nu_0,m}(a(r))\bigr]\,\d r.
\end{equation*}
Each integrand is $\theta\ga(r)$ times the average of the derivative of the corresponding potential over the segment between $a(r)$ and $a_\theta(r)$. Since $\int_0^{T'}|\zeta_{\nu_\theta}-\zeta_{\nu_0}|=W_1(\nu_\theta,\nu_0) \le\theta\|\ga\|_{L^1}$, Corollary~\ref{c.scalar.support} gives $\chi'_{\nu_\theta,m}\to\chi'_{\nu_0,m}$ uniformly as $\theta\downarrow0$, and $\chi'_{\nu_0,m}$ is continuous. Dividing by $\theta$ and letting $\theta\downarrow0$, both bounds converge to the right-hand side of \eqref{e.directional.derivative}, by bounded convergence.

For the finite-step path $q$ with distinct heights and $\ga=\1_{[r_i,r_{i+1})}$, the paths $a_\theta$ are increasing for small $\theta$, and $\chi'_{\nu_q,m}(q_i)=\E_m[M_{q_i}^2]$ by Corollary~\ref{c.scalar.support}, so \eqref{e.directional.derivative} proves \eqref{e.scalar.path.derivative}. For general $a\in\cQ_\infty$, the function $v\mapsto\E[(M^{\nu_a,m}_v)^2]=\chi'_{\nu_a,m}(v)$ is increasing with values in $[m^2,1]$, by the martingale property and Jensen's inequality, and composing it with the increasing path $a$ gives the properties of $z_a$ stated in \eqref{e.def.z.a}. Taking $\ga=b-a$ with $b\in\cQ_\infty$, for which $a_\theta$ is a convex combination of $a$ and $b$, \eqref{e.directional.derivative} says that $\int_0^1z_a(b-a)$ is the right derivative of $\theta\mapsto g_m(a+\theta(b-a))$ at $\theta=0$.

\smallskip
\noindent\emph{Step 2: off-diagonal entries of the Hessian.}\par
Let $K_{ij}:=\dr_{q_i}\dr_{q_j}g_m(q)$ and $U_i(q,z):=\E_z[M_{q_i}^2]$, so that $K_{ij}=w_i\dr_{q_j}U_i(q,m)$ by Step~1. Let $j\ne i$. Increasing $q_j$ lowers $\zeta_{\nu_q}$ on the interval swept by the atom, and leaves the observation time $q_i$ unchanged. By Lemma~\ref{l.diffusion.comparison}(2), the diffusion coefficient $A$ decreases, and then, by Lemma~\ref{l.diffusion.comparison}(3), $\E_m[\varphi(M_{q_i})]$ decreases for the convex function $\varphi(z)=z^2$. Hence $\dr_{q_j}U_i(q,m)\le0$, that is, $K_{ij}\le0$ for all $j\ne i$.

\smallskip
\noindent\emph{Step 3: row sums of the Hessian.}\par
Let $0<\ep<T-q_\ell$, and consider the path $q+\ep\mathbf1$, with all heights shifted by $\ep$. We write $\tilde\Phi$ and $\tilde A$ for the associated objects, and $\tilde M$ for its magnetization process started at time $0$ from $m$. The cumulative function of $\nu_{q+\ep\mathbf1}$ is $\zeta_{\nu_q}(\cdot-\ep)$ on $[\ep,T]$, and $\Phi(T-\ep,\cdot)=\log\cosh+\ep$ by \eqref{e.scalar.top}, since $T-\ep>q_\ell$. Since the operators $T_{r,a}$ commute with the addition of constants, the recursion \eqref{e.def.scalar.recursion} gives $\tilde\Phi(v,x)=\Phi(v-\ep,x)-\ep$ for $v\in[\ep,T]$, so that $\tilde A(v,z)=A(v-\ep,z)$. By \eqref{e.scalar.diffusion} and \eqref{e.def.magnetization.process}, it follows that for every $z\in(-1,1)$, the processes $(\tilde M^{\ep,z}_{\ep+r})_{r\in[0,T-\ep]}$ and $(M^{0,z}_r)_{r\in[0,T-\ep]}$ have the same law, where $M^{0,z}$ is the magnetization process for $q$. The $i$-th height of $q+\ep\mathbf1$ is $q_i+\ep$, so by the Markov property, the convexity in the starting point of Lemma~\ref{l.diffusion.comparison}(1), and Jensen's inequality for the martingale $\tilde M$,
\begin{equation*}
 U_i(q+\ep\mathbf1,m)
   =\E\bigl[\tilde M_{q_i+\ep}^2\bigr]
   =\E\bigl[U_i(q,\tilde M_\ep)\bigr]
   \ge U_i\bigl(q,\E\tilde M_\ep\bigr)=U_i(q,m).
\end{equation*}
Differentiating at $\ep=0$ gives $\sum_j\dr_{q_j}U_i(q,m)\ge0$, that is, $\sum_jK_{ij}\ge0$.

\smallskip
\noindent\emph{Step 4: conclusion.}\par
For every real vector $\theta=(\theta_0,\ldots,\theta_\ell)$,
\begin{equation*}
 \theta^{\mathsf T}K\theta
 =\sum_i\Bigl(\sum_jK_{ij}\Bigr)\theta_i^2
                +\sum_{i<j}(-K_{ij})(\theta_i-\theta_j)^2\ge0,
\end{equation*}
by Steps~2 and~3. This proves that $g_m$ is convex on the open cone of distinct positive heights. The Lipschitz bound \eqref{e.scalar.lipschitz} is the $W_1$ bound of Lemma~\ref{l.legendre.m}, since $W_1(\nu_q,\nu_{q'})=\|q-q'\|_{L^1}$. Convexity extends by continuity to coincident or zero heights, hence to all finite-step paths, any two of which share a common partition, and then by density and \eqref{e.scalar.lipschitz} to $\cQ_1$. For $m=\pm1$, the function $g_m$ is affine by \eqref{e.g.endpoints}.

\smallskip
\noindent\emph{Step 5: the supporting inequality.}\par
Let $|m|<1$ and $a,b\in\cQ_\infty$. Since $g_m$ is convex on $\cQ_1$, the function $\theta\mapsto g_m(a+\theta(b-a))$ is convex on $[0,1]$, so it lies above its tangent line at $\theta=0$, which by Step~1 is $\theta\mapsto g_m(a)+\theta\int_0^1z_a(b-a)$. Taking $\theta=1$ gives \eqref{e.scalar.supporting}.
\end{proof}

We now return to the multi-species setting. Recall from \eqref{e.def.slice.initial} that, for $m\in[-1,1]^{\sS}$ and $q\in\cQ_1^{\sS}$, we have $\psi_m(q)=\sum_s\lambda_sg_{m_s}(q_s)$. We write $\cH$ for $L^2([0,1);\R^{\sS})$, and $\la\cdot,\cdot\ra_{\cH}$ for its inner product. Let $(\cQ_2^{\sS})^*$ denote the dual cone of $\cQ_2^{\sS}$ in $\cH$, that is, the set of $\kappa\in\cH$ with $\la\kappa,q\ra_\cH\ge0$ for every $q\in\cQ_2^{\sS}$. It has the concrete description \cite[Lemma~3.4]{chen2022hamilton}: for $\kappa\in\cH$,
\begin{equation}\label{e.dual.cone.characterization}
 \kappa\in(\cQ_2^{\sS})^*
 \quad\Longleftrightarrow\quad
 \int_r^1\kappa_s(u)\,\d u\ge0\quad\text{for every }s\in\sS
 \text{ and }r\in[0,1).
\end{equation}
Indeed, the condition is necessary by testing against the paths $\1_{[r,1)}e_s\in\cQ_2^{\sS}$, and it is sufficient because every finite-step path in $\cQ_2^{\sS}$ is a combination of such paths with nonnegative coefficients, and these are dense. In particular, $\cQ_2^{\sS}\subset(\cQ_2^{\sS})^*$. We say that a function $u$ on $\cQ_2^{\sS}$ is $(\cQ_2^{\sS})^*$-increasing if $u(q')\ge u(q)$ whenever $q'-q\in(\cQ_2^{\sS})^*$. For every $a\ge0$, we write $\cQ_{\infty,\le a}$ for the set of paths in $\cQ_\infty$ that remain bounded by $a$, and we set
\begin{equation}\label{e.def.Q.lambda}
 \cQ_{\infty,\le\lambda}^{\sS}:=\prod_{s\in\sS}\cQ_{\infty,\le\lambda_s},
\end{equation}
that is, the set of $p\in\cQ_\infty^{\sS}$ such that $0\le p_s\le\lambda_s$ for every $s\in\sS$.

\begin{corollary}[The slice initial condition]\label{c.psi.m}
For every $m\in[-1,1]^{\sS}$, the function $\psi_m$ is convex and $(\cQ_2^{\sS})^*$-increasing on $\cQ_2^{\sS}$, and
\begin{equation}\label{e.initial.lipschitz}
 |\psi_m(q)-\psi_m(q')|\le\sum_s\lambda_s\|q_s-q_s'\|_{L^1}
 \le\|q-q'\|_{L^1},\qquad q,q'\in\cQ_1^{\sS}.
\end{equation}
For every $q\in\cQ_\infty^{\sS}$, the path $p:=(\lambda_sz_{q_s})_{s\in\sS}$, where $z_{q_s}$ is the path \eqref{e.def.z.a} for the magnetization $m_s$, belongs to $\cQ_{\infty,\le\lambda}^{\sS}$ and supports $\psi_m$ at $q$:
\begin{equation}\label{e.vector.supporting}
 \psi_m(q')\ge\psi_m(q)+\la q'-q,p\ra_\cH,
 \qquad q'\in\cQ_\infty^{\sS}.
\end{equation}
Finally, let $\mcl I_\lambda(m):=\sum_s\lambda_sI(m_s)$. For every $q\in\cQ_1^{\sS}$, the function $m\mapsto\psi_m(q)-\mcl I_\lambda(m)$ is convex on $[-1,1]^{\sS}$, and
\begin{equation}\label{e.initial.magnetization.lipschitz}
 \bigl|\psi_m(q)-\psi_{m'}(q)-\mcl I_\lambda(m)+\mcl I_\lambda(m')\bigr|
 \le2\sum_s\lambda_s\|q_s\|_{L^1}|m_s-m_s'|,
 \qquad m,m'\in[-1,1]^{\sS};
\end{equation}
in particular, $m\mapsto\psi_m(q)$ is continuous on $[-1,1]^{\sS}$.
\end{corollary}

\begin{proof}
Convexity and the Lipschitz bound follow from Proposition~\ref{p.path.convexity}, and so does \eqref{e.vector.supporting}, by summing the supporting inequalities \eqref{e.scalar.supporting} with the weights $\lambda_s$. For the monotonicity, let $q,q'\in\cQ_\infty^{\sS}$ with $q'-q\in(\cQ_2^{\sS})^*$. Since $p$ belongs to $\cQ_{\infty,\le\lambda}^{\sS}\subset\cQ_2^{\sS}$, we have $\la q'-q,p\ra_\cH\ge0$, and \eqref{e.vector.supporting} gives $\psi_m(q')\ge\psi_m(q)$. For general $q,q'\in\cQ_2^{\sS}$ with $q'-q\in(\cQ_2^{\sS})^*$, we approximate $q$ and $q'$ in $L^1$ by their averages over the dyadic intervals of length $2^{-n}$, which are bounded paths; this preserves the relation $q'-q\in(\cQ_2^{\sS})^*$, by \eqref{e.dual.cone.characterization}, and we conclude by \eqref{e.initial.lipschitz}. The statements in $m$ follow from Lemma~\ref{l.legendre.m} for bounded paths, since $\int v\,\nu_{q_s}(\d v)=\|q_s\|_{L^1}$, and they extend to $\cQ_1^{\sS}$ by \eqref{e.initial.lipschitz}, since the right-hand side of \eqref{e.initial.magnetization.lipschitz} is continuous in $q$ for the $L^1$ norm.
\end{proof}
%
%
%
%
\section{The initial condition on the slices}
\label{s.initial}

In this section, we identify the initial condition of the constrained models: we show that $\bar F_{N,k}(0,\cdot)$ is close to $\psi_{m_{N,k}}$ as $N$ tends to infinity, uniformly over the slices. We first show that, for every $t\ge0$, the free energy of the full model is, up to a vanishing error, the minimum of the slice free energies, each corrected by the contribution of the external field. This relation will be used again in Section~\ref{s.conclusion} to identify the limit free energy of the full model; here, we only need it at $t=0$. At $t=0$, the free energies of the full model and of the slices factorize over the species, which reduces the problem to a single species, and the free energy of the full model with external field $h\in\R$ is then exactly $\psi_\circ(\cdot;h)$. The relation above thus states that $h\mapsto-\psi_\circ(\cdot;h)$ is, up to a small error, the convex conjugate of the slice free energies, viewed as functions of the magnetization; by Lemma~\ref{l.legendre.m}, it is also the convex conjugate of $m\mapsto g_m$. A priori, a convex conjugate only determines the convex envelope of a function, but here the strong convexity (and continuity) of $g_m$ in $m$ allows us to conclude.

\subsection{Gluing the slices}
\label{s.initial.gluing}

Recall that $\bar F_N$ denotes the free energy of the full cube with the weights $w(\si)=2^{-N}e^{h^{(N)}\cdot\si}$, where $h^{(N)}_i=h_s$ for $i\in I_{N,s}$, and that, by Lemma~\ref{l.zero}, its value at $q=0$ is exactly \eqref{e.def.free.energy}. Let $U^{\Sigma_N}$ and $U^{\Omega_{N,k}}$ denote the solutions of \eqref{e.branch.pde} for the full cube with these weights, and for the slices with the weights $2^{-N}$. Since the field $h^{(N)}\cdot\si$ is constant on each slice, with value $N\lambda_N\cdot(hm_{N,k})$, where $hm:=(h_sm_s)_s$, the terminal function of the full cube is a logarithmic sum of the slice terminal functions:
\begin{equation}\label{e.terminal.gluing}
 U^{\Sigma_N}(1,x;g)
 =\log\sum_k\exp\bigl(N\lambda_N\cdot(hm_{N,k})+U^{\Omega_{N,k}}(1,x;g)\bigr).
\end{equation}
Throughout this section, we write $n_{N,s}:=|I_{N,s}|$ for the number of spins of species $s$. The number of slices is $M_N:=\prod_s(n_{N,s}+1)\le(N+1)^{|\sS|}$, and we set
\begin{equation}\label{e.def.aN}
 a_N:=\frac{|\sS|\log(N+1)}N\ge\frac{\log M_N}N.
\end{equation}

\begin{lemma}[Gluing the slices]\label{l.gluing}
For every $N$, $t\ge0$, $h\in\R^{\sS}$, and every $q\in\cQ_1^{\sS}$,
\begin{equation}\label{e.gluing}
 0\le\min_k\bigl\{\bar F_{N,k}(t,q)-\lambda_N\cdot(hm_{N,k})\bigr\}
 -\bar F_N(t,q)
 \le a_N+2\sqrt{a_N\bigl(t\xi(\mathbf1)+\|q\|_{L^1}\bigr)}.
\end{equation}
\end{lemma}

This lemma is perhaps most easily understood through the cascade representation of Appendix~\ref{s.cascade}, where the quenched partition function of the full model is clearly the sum of the partition functions over the slices. This would easily yield a statement in the spirit of \eqref{e.gluing} at the level of the quenched free energies, but then a concentration argument, uniform over the slices, would need to be supplied in order to conclude. This is delicate, since the randomness of the cascade includes its weights, which are not Gaussian, and since the total variance of the field at a leaf is $2N\lambda_N\cdot q(1)$, which is not controlled by the $L^1$ norm of the path.

The proof given here avoids any concentration argument over the cascade, as we work directly with the branch equation. The terminal datum of the full cube is the logarithmic sum of those of the slices, and the lower bound is the comparison principle. For the upper bound, we would like the logarithmic sum of the slice solutions to stay above the solution of the full cube for every $r$. But a logarithmic sum of solutions is not a solution, and the computation \eqref{e.logsum} shows that it is a subsolution, which lies below the solution of the full cube and only reproduces the lower bound. We therefore use instead the tempered logarithmic sum $\rho^{-1}\log\sum_ke^{\rho V_k}$ with $\rho<1$, which is a supersolution on $[\rho,1]$, and on $[0,\rho)$ we add a correction which involves only the values of the path on $[0,\rho]$. The concentration inequality is then used only at the root, for $g$ and $X_0$, where it is Gaussian, to compare the tempered logarithmic sum with the maximum; the value at the root of the correction cancels the contribution of the initial field to this estimate, leaving an error controlled by the $L^1$ norm of the path, and the optimization over $\rho$ produces the term of order $\sqrt{a_N}$ in \eqref{e.gluing}.

\begin{proof}[Proof of Lemma~\ref{l.gluing}]
We first assume that $q$ is regular, and we denote by $\msc G_N(q)$ the middle term in~\eqref{e.gluing}. Let
\begin{equation*}
 V_k(r,x):=N\lambda_N\cdot(hm_{N,k})+U^{\Omega_{N,k}}(r,x;g),
 \qquad
 V(r,x):=U^{\Sigma_N}(r,x;g).
\end{equation*}
By Lemma~\ref{l.symmetries}(1), the functions $V_k$ solve the branch equation \eqref{e.branch.pde} for the path $q$ on $[0,1]\times\R^N$, and so does $V$; by \eqref{e.terminal.gluing}, their terminal data satisfy $V(1,\cdot)=\log\sum_ke^{V_k(1,\cdot)}$. With $X_0=\sqrt{2\diag(q(0))}\,Z$, the identity \eqref{e.free.energy.V} gives
\begin{align*}
 \bar F_{N,k}(t,q)-\lambda_N\cdot(hm_{N,k})&=t\xi(\lambda_N)+\lambda_N\cdot q(1)-\frac1N\E\,V_k(0,X_0),\\
 \bar F_N(t,q)&=t\xi(\lambda_N)+\lambda_N\cdot q(1)-\frac1N\E\,V(0,X_0),
\end{align*}
so that
\begin{equation}\label{e.D.N}
 \msc G_N(q)=\frac1N\Bigl(\E\,V(0,X_0)-\max_k\E\,V_k(0,X_0)\Bigr).
\end{equation}

\smallskip
\noindent\emph{Step 1: a logarithmic sum of solutions.}\par
Fix $\rho\in(0,1]$, and let $W:=\rho^{-1}\log\sum_ke^{\rho V_k}$ and $\pi_k:=e^{\rho V_k}/\sum_je^{\rho V_j}$, so that $\sum_k\pi_k=1$. Then $\dr_rW=\sum_k\pi_k\dr_rV_k$, $\nabla W=\sum_k\pi_k\nabla V_k$, and
\begin{equation*}
 \nabla^2W=\sum_k\pi_k\nabla^2V_k
 +\rho\sum_k\pi_k(\nabla V_k-\nabla W)(\nabla V_k-\nabla W)^{\mathsf T}.
\end{equation*}
Since each $V_k$ solves the branch equation, we have $\mcl L[V_k]=0$, where $\mcl L$ is the operator defined in \eqref{e.def.branch.operator}. Using this and the identity $\sum_k\pi_k\nabla V_k\cdot D\nabla V_k-\nabla W\cdot D\nabla W=\sum_k\pi_k(\nabla V_k-\nabla W)\cdot D(\nabla V_k-\nabla W)$, we find
\begin{equation}\label{e.logsum}
 \mcl L[W]=(\rho-r)\,\mcl V,
 \qquad
 \mcl V:=\sum_k\pi_k(\nabla V_k-\nabla W)\cdot D(\nabla V_k-\nabla W).
\end{equation}
The matrix $D(r)$ is diagonal with entries $\dot q_s(r)$ on $I_{N,s}$, and $|\dr_{x_i}V_k|\le1$ for every $i$ and $k$, so that the variance of $\dr_{x_i}V_k$ under the weights $\pi_k$ is at most $1$, and
\begin{equation*}
 0\le\mcl V\le\tr D(r)=N\lambda_N\cdot\dot q(r).
\end{equation*}
The sign of $\mcl L[W]$ is favorable only for $r\ge\rho$. For $r<\rho$, we compensate it by the function $c(r):=N\int_r^1(\rho-u)_+\,\lambda_N\cdot\dot q(u)\,\d u$, which satisfies $\mcl L[W+c]=(\rho-r)\mcl V-(\rho-r)_+N\lambda_N\cdot\dot q(r)\le0$ on $[0,1)\times\R^N$. The function $W+c$ is continuous with bounded gradient and classical on each piece of $q$. At $r=1$, we have $W(1,\cdot)=\rho^{-1}\log\sum_ke^{\rho V_k(1,\cdot)}\ge\log\sum_ke^{V_k(1,\cdot)}=V(1,\cdot)$, because $p\mapsto(\sum_ka_k^p)^{1/p}$ is decreasing for positive $a_k$ and $\rho\le1$. Lemma~\ref{l.comparison}, applied to $W_1=W+c$ and $W_2=V$, gives $V\le W+c$; and since $V_k(1,\cdot)\le V(1,\cdot)$, it also gives $V_k\le V$ for every $k$. Hence, for every $x\in\R^N$,
\begin{equation}\label{e.root.comparison}
 \max_kV_k(0,x)\le V(0,x)\le\frac1\rho\log\sum_ke^{\rho V_k(0,x)}+c(0),
\end{equation}
where an integration by parts gives $c(0)=N\bigl(\int_0^\rho\lambda_N\cdot q(u)\,\d u-\rho\,\lambda_N\cdot q(0)\bigr)$.

\smallskip
\noindent\emph{Step 2: Gaussian concentration at the root.}\par
We evaluate \eqref{e.root.comparison} at $x=X_0$ and take expectations. The lower bound and \eqref{e.D.N} give $\msc G_N(q)\ge0$. For the upper bound, each $Y_k:=V_k(0,X_0)$ is a function of the standard Gaussian vector $(g,Z)$. By \eqref{e.response.first}, applied to the branch of the slice $\Omega_{N,k}$ started from $(0,x)$, the gradient of $U^{\Omega_{N,k}}(0,x;g)$ in $g$ is $\sqrt{2t}$ times a conditional mean of $v_\si$, hence of norm at most $\sqrt{2tN\xi(\lambda_N)}$ by \eqref{e.gram}; and the gradient of $V_k(0,X_0)$ in $Z$ is $\sqrt{2\diag(q(0))}\,\nabla_xV_k$, of squared norm at most $2\sum_sn_{N,s}q_s(0)=2N\lambda_N\cdot q(0)$, since $|\dr_{x_i}V_k|\le1$. Thus $Y_k$ is a Lipschitz function of $(g,Z)$ with constant $\Lambda:=\sqrt{2N(t\xi(\lambda_N)+\lambda_N\cdot q(0))}$, and the Gaussian concentration inequality, in its exponential form $\E\,e^{\rho(Y_k-\E Y_k)}\le e^{\rho^2\Lambda^2/2}$, together with Jensen's inequality, gives
\begin{equation*}
 \E\,\frac1\rho\log\sum_ke^{\rho Y_k}
 \le\frac1\rho\log\sum_k\E\,e^{\rho Y_k}
 \le\max_k\E\,Y_k+\frac{\log M_N}\rho+\rho N\bigl(t\xi(\lambda_N)+\lambda_N\cdot q(0)\bigr).
\end{equation*}
We insert this into \eqref{e.root.comparison} and divide by $N$. The term $\rho\,\lambda_N\cdot q(0)$ cancels against $c(0)/N$, and \eqref{e.D.N} and \eqref{e.def.aN} yield
\begin{equation}\label{e.gluing.rho}
 \msc G_N(q)\le\frac{a_N}\rho+\rho\,t\xi(\lambda_N)+\int_0^\rho\lambda_N\cdot q(u)\,\d u
 \le\frac{a_N}\rho+\rho\bigl(t\xi(\mathbf1)+\|q\|_{L^1}\bigr),
\end{equation}
where we used that $\xi(\lambda_N)\le\xi(\mathbf1)$, and that the average of the increasing function $\lambda_N\cdot q$ over $[0,\rho]$ is at most its average over $[0,1]$, which is $\sum_s\lambda_{N,s}\|q_s\|_{L^1}\le\|q\|_{L^1}$.

\smallskip
\noindent\emph{Step 3: choice of $\rho$, and general paths.}\par
Let $C:=t\xi(\mathbf1)+\|q\|_{L^1}$. If $C\le a_N$, we take $\rho=1$ in \eqref{e.gluing.rho}; otherwise, we take $\rho=\sqrt{a_N/C}$. In both cases, $\msc G_N(q)\le a_N+2\sqrt{a_NC}$, which is \eqref{e.gluing} for regular paths. For $q\in\cQ_1^{\sS}$, let $q^n$ be regular paths with $\|q^n-q\|_{L^1}\to0$, which exist by the density noted in the proof of Proposition~\ref{p.extension}. All the free energies in \eqref{e.gluing} satisfy the Lipschitz bound \eqref{e.finite.lipschitz}, and so does their minimum over $k$, so that $|\msc G_N(q)-\msc G_N(q^n)|\le2\|q^n-q\|_{L^1}$; since $\|q^n\|_{L^1}\to\|q\|_{L^1}$, \eqref{e.gluing} for $q^n$ gives \eqref{e.gluing} for $q$.
\end{proof}

\subsection{Factorization at $t=0$}
\label{s.initial.factorization}

At $t=0$, the branch equation separates over blocks of coordinates. We make this precise, and we use it to express the initial conditions of the full cube and of the slices through one-species quantities. For $n\ge1$, $0\le k\le n$ and $a\in\cQ_1$, let
\begin{equation}\label{e.def.G}
 G_{n,k}(a):=\bar F_{n,k}(0,a)
 \qquad\text{for the system with a single species and $n$ sites},
\end{equation}
that is, the quantity \eqref{e.def.branch.free.energy} at $t=0$, extended to $\cQ_1$ by Proposition~\ref{p.extension}, with $\sS$ a singleton, $N=n$, $\lambda_N=1$, and $\Omega=\Omega_{n,k}=\{\si\in\{-1,1\}^n:|\{i:\si_i=1\}|=k\}$ with the weights $2^{-n}$. We write $m_k:=2k/n-1$ for the magnetization of this slice. Recall that $\psi_\circ$ is defined in \eqref{e.def.g.m}. 

\begin{lemma}[Factorization at $t=0$]\label{l.factorization}
For every $N$, $q\in\cQ_1^{\sS}$, $h\in\R^{\sS}$ and $k$, we have
\begin{align}
 \bar F_N(0,q)&=\sum_s\lambda_{N,s}\,\psi_\circ(q_s;h_s),
                                            \label{e.factorization}\\
 \bar F_{N,k}(0,q)&=\sum_s\lambda_{N,s}\,G_{n_{N,s},k_s}(q_s).
                                            \label{e.slice.factorization}
\end{align}
Moreover, for every $n\ge1$ and $a\in\cQ_1$,
\begin{equation}\label{e.single.configuration}
 G_{n,0}(a)=G_{n,n}(a)=\log2+\int_0^1a(r)\,\d r.
\end{equation}
\end{lemma}

\begin{proof}
By the Lipschitz bound \eqref{e.finite.lipschitz}, by Lemma~\ref{l.scalar.continuity} and by the density of regular paths, it suffices to consider a regular path $q$. At $t=0$, the slice $\Omega_{N,k}$ is the product over $s$ of the one-species slices $\Omega_{n_{N,s},k_s}$ in the coordinates $I_{N,s}$, and the terminal function \eqref{e.def.terminal} is
\begin{equation*}
 U(1,x)=\sum_s\log\sum_{\si_s\in\Omega_{n_{N,s},k_s}}
            2^{-n_{N,s}}\exp(x_{I_{N,s}}\cdot\si_s),
\end{equation*}
the sum over the blocks $I_{N,s}$ of the terminal functions of the one-species slices. By Lemma~\ref{l.symmetries}(2), $U(0,x)=\sum_sU_s(0,x_{I_{N,s}})$, where $U_s$ is the solution of the branch equation for the one-species slice $\Omega_{n_{N,s},k_s}$ and the path $q_s$; moreover the blocks $X_{0,I_{N,s}}=\sqrt{2q_s(0)}\,Z_{I_{N,s}}$ of the initial field are independent, each with the law of the initial field of the corresponding one-species system. Hence $\E\,U(0,X_0)$ is the sum over $s$ of the corresponding one-species quantities, and \eqref{e.slice.factorization} follows from \eqref{e.free.energy.V}, since $\lambda_N\cdot q(1)=\sum_s\lambda_{N,s}q_s(1)$ and $N^{-1}=\lambda_{N,s}n_{N,s}^{-1}$.

For the full cube with the weights $2^{-N}e^{h^{(N)}\cdot\si}$, we take the sites as blocks. The terminal function is $U^{\Sigma_N}(1,x)=\sum_s\sum_{i\in I_{N,s}}\log\cosh(x_i+h_s)$. By Lemma~\ref{l.symmetries}(1) and (2), $U^{\Sigma_N}(0,x)=\sum_s\sum_{i\in I_{N,s}}u_s(0,x_i+h_s)$, where $u_s$ is the solution of the one-spin branch equation \eqref{e.one.spin.branch} for the path $q_s$, and the coordinates of $X_0$ are independent with $X_{0,i}=\sqrt{2q_s(0)}\,Z_i$ for $i\in I_{N,s}$. By Lemma~\ref{l.dictionary}, $\E\,u_s(0,h_s+\sqrt{2q_s(0)}\,Z)=q_s(1)-\psi_\circ(q_s;h_s)$, hence
\begin{equation*}
 \E\,U^{\Sigma_N}(0,X_0)
 =\sum_sn_{N,s}\bigl[q_s(1)-\psi_\circ(q_s;h_s)\bigr],
\end{equation*}
and \eqref{e.factorization} follows from \eqref{e.free.energy.V}.

Finally, the slice $\Omega_{n,n}$ consists of the single configuration $\mathbf1$, so that the overlap $R_n(\si^1,\si^2)=\frac1n\si^1\cdot\si^2$ of two replicas is identically $1$. For a regular path $a$ and $\theta>0$, the path $\theta a$ is regular, and $\theta\mapsto\theta a$ is an admissible variation with $h=a$, so Theorem~\ref{t.derivatives} gives
\begin{equation*}
 \frac{\d}{\d\theta}G_{n,n}(\theta a)=\int_0^1a(r)\,\E\langle R_n\rangle_r\,\d r=\int_0^1a(r)\,\d r
 \qquad(\theta>0).
\end{equation*}
The function $\theta\mapsto G_{n,n}(\theta a)$ is continuous on $[0,\infty)$ by \eqref{e.finite.lipschitz}, and $G_{n,n}(0)=-\frac1n\log2^{-n}=\log2$ by Lemma~\ref{l.zero}. Integrating over $[0,1]$ proves \eqref{e.single.configuration} for $k=n$ and regular $a$; the slice $\Omega_{n,0}=\{-\mathbf1\}$ is treated identically, and the general case follows by \eqref{e.finite.lipschitz}.
\end{proof}

\subsection{Continuity in the magnetization after subtracting the entropy}
\label{s.initial.continuity}

We next compare two slices of the one-species system. Near $m=\pm1$, the slice free energies depend on the magnetization in a singular way, through the entropy of the slice, and the comparison is made after subtracting this entropy. We will use a coupling argument between the two slices that incorporates the appropriate entropy correction, and proceed by interpolation. The interpolation goes through a system with $2n$ coordinates in which, at the endpoints, the path of one of the two blocks vanishes; we first identify the free energy of such a system.

\begin{lemma}[A block with zero path]\label{l.block}
Let $t=0$, let $\sS=\sS_0\cup\sS_1$ be a partition of the species, and let $q=(q^{(0)},0)$, where $q^{(0)}$ is a regular path indexed by $\sS_0$ and $q_s\equiv0$ for $s\in\sS_1$. Let $N_0$ be the number of coordinates of the species in $\sS_0$, write $\si_{(0)}$ for the restriction of $\si\in\Omega$ to these coordinates, let $\Omega_{(0)}:=\{\si_{(0)}:\si\in\Omega\}$, and let $w_{(0)}(\tau):=\sum_{\si\in\Omega:\si_{(0)}=\tau}w(\si)$. Then
\begin{equation}\label{e.block}
 N\bar F_{N,\Omega}(0,q)=N_0\bar F_{N_0,\Omega_{(0)}}(0,q^{(0)}),
\end{equation}
where the right-hand side is the free energy of the system with the $N_0$ coordinates of the species in $\sS_0$, the configuration set $\Omega_{(0)}$ and the weights $w_{(0)}$.
\end{lemma}

\begin{proof}
For $\ep>0$, let $q^\ep:=(q^{(0)},\ep(1+r)\mathbf1)$, which is regular and converges to $q$ in $L^1$, so that $\bar F_{N,\Omega}(0,q^\ep)\to\bar F_{N,\Omega}(0,q)$ by \eqref{e.finite.lipschitz}. Let $N_1:=N-N_0$, write $x=(x_{(0)},x_{(1)})$ according to the two groups of species, let $U^\ep$ be the solution of \eqref{e.branch.pde} for $q^\ep$ in dimension $N$, and recall from \eqref{e.hessian.bounds} that $0\le\tr\nabla^2U^\ep\le N$ and $|\nabla U^\ep|^2\le N$. Let $V$ be the solution of the branch equation in dimension $N_0$, for the path $q^{(0)}$, with the terminal datum $x_{(0)}\mapsto U^\ep(1,x_{(0)},0)=\log\sum_{\tau\in\Omega_{(0)}}w_{(0)}(\tau)e^{x_{(0)}\cdot\tau}$, which is the terminal function of the reduced system. The function $W(r,x_{(0)}):=U^\ep(r,x_{(0)},0)$ solves the equation of the reduced system up to the source $-\ep(\tr\nabla_{(1)}^2U^\ep+r|\nabla_{(1)}U^\ep|^2)(r,x_{(0)},0)$, which is bounded in absolute value by $2\ep N$, so Lemma~\ref{l.comparison} gives $|W-V|\le2\ep N$. Moreover $X^\ep_0=(X_{0,(0)},\sqrt{2\ep}Z_{(1)})$, with $X_{0,(0)}\sim\mcl N(0,2\diag(q^{(0)}(0)))$ the initial field of the reduced system and $Z_{(1)}$ a standard Gaussian vector in $\R^{N_1}$, so that \eqref{e.heat.interpolation}, applied in the second group of coordinates, gives $0\le\E\,U^\ep(0,X^\ep_0;g)-\E\,W(0,X_{0,(0)})\le N\ep$. By \eqref{e.free.energy.V}, since $N\lambda_{N,s}=n_{N,s}$,
\begin{align*}
 N\bar F_{N,\Omega}(0,q^\ep)&=\sum_{s\in\sS_0}n_{N,s}q_s(1)+2\ep N_1-\E\,U^\ep(0,X^\ep_0;g),\\
 N_0\bar F_{N_0,\Omega_{(0)}}(0,q^{(0)})&=\sum_{s\in\sS_0}n_{N,s}q_s(1)-\E\,V(0,X_{0,(0)}),
\end{align*}
and letting $\ep\to0$ proves \eqref{e.block}.
\end{proof}

\begin{lemma}[Entropy-corrected continuity in the magnetization]\label{l.composition}
For every $n\ge1$, $0\le k,k'\le n$, and $a\in\cQ_1$, we have
\begin{equation}\label{e.composition}
 \bigl|(G_{n,k}(a)-I_{n,k})-(G_{n,k'}(a)-I_{n,k'})\bigr|
 \le2\|a\|_{L^1}\,|m_k-m_{k'}|,
\end{equation}
where $I_{n,k}:=\log2-\frac1n\log\binom nk$ satisfies $|I_{n,k}-I(m_k)|\le\frac{\log(n+1)}n$.
\end{lemma}

\begin{proof}
We first prove the entropy estimate. Let $\mathsf h(p):=-p\log p-(1-p)\log(1-p)$ for $p\in[0,1]$, with $0\log0=0$. Writing $p=\frac{1+z}2$ in the definition of $I$, we get $I(z)=p\log(2p)+(1-p)\log(2(1-p))=\log2-\mathsf h(\frac{1+z}2)$, and since $\frac{1+m_k}2=\frac kn$, we have
\begin{equation*}
 I_{n,k}-I(m_k)=\mathsf h\bigl(\tfrac kn\bigr)-\frac1n\log\binom nk.
\end{equation*}
Let $p:=k/n$ and $b_j:=\binom nj p^j(1-p)^{n-j}$ for $0\le j\le n$, so that $\sum_jb_j=1$ and $b_k=\binom nk e^{-n\mathsf h(p)}$. In particular, $b_k\le1$. If $k\in\{0,n\}$, then $b_k=1$. Otherwise, for $0\le j<n$, the ratio $b_{j+1}/b_j=\frac{(n-j)k}{(j+1)(n-k)}$ is at least $1$ for $j<k$ and less than $1$ for $j\ge k$, since $(n-j)k-(j+1)(n-k)=n(k-1-j)+k$. So $b_k$ is the largest of the $n+1$ numbers $b_j$, and $b_k\ge(n+1)^{-1}$. Thus $(n+1)^{-1}e^{n\mathsf h(k/n)}\le\binom nk\le e^{n\mathsf h(k/n)}$, that is, $0\le I_{n,k}-I(m_k)\le\frac{\log(n+1)}n$.

We now turn to the proof of \eqref{e.composition}. By~\eqref{e.finite.lipschitz}, both sides of this inequality are continuous in $a$ for the $L^1$ norm, so it suffices to consider a regular path $a$.

\smallskip
\noindent\emph{Step 1: uniform weights, and a coupling of the two slices.}\par
Since $2^{-n}=e^{-nI_{n,k}}\binom nk^{-1}$, and since adding a constant to the terminal datum adds it to the solution of the branch equation, by Lemma~\ref{l.symmetries}(1), the quantity $G_{n,k}(a)-I_{n,k}$ is the free energy at $t=0$ of the slice $\Omega_{n,k}$ with the uniform weights $\binom nk^{-1}$. Let $U_1,\dots,U_n$ be independent random variables, uniformly distributed on $[0,1]$. Let $\si\in\Omega_{n,k}$ be the configuration whose plus spins are the $k$ sites with the smallest values of $U_i$, and define $\si'\in\Omega_{n,k'}$ in the same way with $k'$. Then $\si$ and $\si'$ are each uniformly distributed on their slices, and they differ at exactly $d:=|k-k'|$ sites. Let $\tilde\Omega\subset\{-1,1\}^{2n}$ be the set of pairs $(\si,\si')$ obtained in this way, with the weights $w(\si,\si')$ given by their probabilities. Since $\sum_{\si'}w(\si,\si')=\binom nk^{-1}$ for every $\si\in\Omega_{n,k}$, the slice $\Omega_{n,k}$ with the uniform weights is the projection of $(\tilde\Omega,w)$ on the first block of $n$ coordinates, in the sense of Lemma~\ref{l.block}, and $\Omega_{n,k'}$ with the uniform weights is its projection on the second block.

\smallskip
\noindent\emph{Step 2: a doubled system.}\par
For $u\in[0,1]$, consider the system with $2n$ coordinates, partitioned into two species, the first and the second block of $n$ coordinates, with the configuration set $\tilde\Omega$ and the weights $w$, at $t=0$, with the path $q^u:=((1-u)a,ua)$. Its free energy, normalized by $2n$ as in Definition~\ref{d.free.energy}, is $\bar F_{2n,\tilde\Omega}(0,q^u)$, and we set $G(u):=2\bar F_{2n,\tilde\Omega}(0,q^u)$. The auxiliary field of this system has two blocks, driven by independent Brownian motions, with covariance paths $(1-u)a$ and $ua$; the terminal weights couple the two blocks. Since $\lambda=(\frac12,\frac12)$ for this system, the self-overlap correction $2\lambda\cdot q^u(1)=a(1)$ does not depend on $u$. The function $G$ is continuous on $[0,1]$ by \eqref{e.finite.lipschitz}. At $u=0$, the second path vanishes, and Lemma~\ref{l.block} together with Step~1 gives $2n\bar F_{2n,\tilde\Omega}(0,(a,0))=n\bigl[G_{n,k}(a)-I_{n,k}\bigr]$, that is, $G(0)=G_{n,k}(a)-I_{n,k}$; similarly, $G(1)=G_{n,k'}(a)-I_{n,k'}$. For $u\in(0,1)$, the path $q^u$ is regular, and $u\mapsto q^u$ is an admissible variation with $h=(-a,a)$, so by Theorem~\ref{t.derivatives}, with the overlaps $\frac1{2n}\si^1\cdot\si^2$ and $\frac1{2n}\si'^1\cdot\si'^2$ of the two blocks, the function $G$ is differentiable on $(0,1)$, with
\begin{equation*}
 G'(u)=\frac1n\int_0^1a(r)\,
 \E\bigl\langle\si'^1\cdot\si'^2-\si^1\cdot\si^2\bigr\rangle_r\,\d r,
\end{equation*}
where $(\si^1,\si'^1)$ and $(\si^2,\si'^2)$ are two replicas at level $r$. Since $\si^b$ and $\si'^b$ differ at $d$ sites, at most $2d$ terms of the sum defining $\si'^1\cdot\si'^2-\si^1\cdot\si^2$ are nonzero, each bounded by $2$. Hence $|G'(u)|\le\frac{4d}n\|a\|_{L^1}=2|m_k-m_{k'}|\,\|a\|_{L^1}$ on $(0,1)$, and \eqref{e.composition} follows by the continuity of $G$ on $[0,1]$.
\end{proof}

\subsection{Identification of the initial condition}
\label{s.initial.identification}

We now identify the limit of the slice initial conditions $\bar F_{N,k}(0,\cdot)$ as $N$ tends to infinity. By the factorization \eqref{e.slice.factorization}, this is a statement about a single species, which we prove first.

\begin{proposition}[Initial data on one-species slices]
\label{p.scalar.initial}
For every $R<\infty$, we have
\begin{equation}\label{e.scalar.initial.limit}
 \lim_{n\to\infty}\ \sup_{0\le k\le n}\ \sup_{\substack{a\in\cQ_1\\\|a\|_{L^1}\le R}}
 \bigl|G_{n,k}(a)-g_{m_k}(a)\bigr|=0.
\end{equation}
More precisely, with $\ep_n(R):=\log(n+1)/n+2\sqrt{R\log(n+1)/n}$, for every $0\le k\le n$ and every $a\in\cQ_1$ with $\|a\|_{L^1}\le R$,
\begin{equation}\label{e.scalar.initial.bound}
 0\le G_{n,k}(a)-g_{m_k}(a)\le E_n(R):=\ep_n(R)+4R\sqrt{2\ep_n(R)}+\frac{2\log(n+1)}n.
\end{equation}
\end{proposition}

\begin{proof}
We first assume that $a$ is bounded, so that $\nu_a\in\Prob([0,T])$ for some $T>0$ and the results of Section~\ref{s.scalar} apply to it; note that $\int v\,\nu_a(\d v)=\|a\|_{L^1}\le R$. We write $d_k:=G_{n,k}(a)-g_{m_k}(a)$, and we begin with two preliminary bounds.

Lemma~\ref{l.gluing} at $t=0$, applied to the system with a single species and $n$ sites, for which $a_N$ becomes $\log(n+1)/n$, together with \eqref{e.factorization} for this system, which reads $\bar F_n(0,a)=\psi_\circ(a;h)$, gives
\begin{equation}\label{e.scalar.gluing}
 0\le\min_{0\le j\le n}\{G_{n,j}(a)-hm_j\}-\psi_\circ(a;h)\le\ep_n(R)
 \qquad\text{for every }h\in\R.
\end{equation}
The lower bound gives $G_{n,k}(a)\ge\psi_\circ(a;h)+hm_k$ for every $k$ and every $h\in\R$; taking the supremum over $h$ and recalling the definition \eqref{e.def.g.m} of $g_m$, we obtain $d_k\ge0$ for every $k$. Next, by Lemma~\ref{l.composition}, including its entropy estimate, and by \eqref{e.magnetization.lipschitz}, for all $j,k$,
\begin{multline}\label{e.d.continuity}
 |d_k-d_j|\le\bigl|(G_{n,k}(a)-I_{n,k})-(G_{n,j}(a)-I_{n,j})\bigr|
 +\frac{2\log(n+1)}n\\
 +\bigl|(g_{m_k}(a)-I(m_k))-(g_{m_j}(a)-I(m_j))\bigr|
 \le4R|m_k-m_j|+\frac{2\log(n+1)}n.
\end{multline}

By \eqref{e.scalar.duality}, the upper bound in \eqref{e.scalar.gluing} says that the Legendre transforms in the magnetization of $k\mapsto G_{n,k}(a)$ and of $m\mapsto g_m(a)$ nearly coincide. A Legendre transform identifies in general only a convex envelope; here, the strong convexity of $g_m$ in $m$ localizes the minimizing slice. Let $0<k<n$, so that $|m_k|<1$, let $h:=y(\nu_a,m_k)$, so that $\psi_\circ(a;h)=g_{m_k}(a)-hm_k$ by \eqref{e.scalar.duality}, and let $j$ attain the minimum in \eqref{e.scalar.gluing}. The upper bound there reads
\begin{equation*}
 d_j+\bigl[g_{m_j}(a)-g_{m_k}(a)-h(m_j-m_k)\bigr]\le\ep_n(R),
\end{equation*}
and by the strong convexity \eqref{e.strong.convexity} at $m_0=m_k$, where $y(\nu_a,m_k)=h$, the bracket is at least $\frac12(m_j-m_k)^2$. Since $d_j\ge0$, we obtain
\begin{equation}\label{e.localization}
 d_j+\tfrac12(m_j-m_k)^2\le\ep_n(R):
\end{equation}
the minimizing slice has both a small error and a magnetization close to $m_k$. Transporting its error to the slice $k$ by \eqref{e.d.continuity},
\begin{equation*}
 d_k\le d_j+4R|m_j-m_k|+\frac{2\log(n+1)}n
 \le\ep_n(R)+4R\sqrt{2\ep_n(R)}+\frac{2\log(n+1)}n=E_n(R).
\end{equation*}
For $k\in\{0,n\}$, we have $m_k=\pm1$, and \eqref{e.single.configuration} and \eqref{e.g.endpoints} give $d_0=d_n=0$. This proves \eqref{e.scalar.initial.bound} for bounded $a$.

Let finally $a\in\cQ_1$ with $\|a\|_{L^1}\le R$, and let $a^M:=a\wedge M$ for $M>0$. Then $a^M$ is bounded, $\|a^M\|_{L^1}\le R$, and $\|a^M-a\|_{L^1}\to0$ as $M\to\infty$. Both $G_{n,k}$ and $g_{m_k}$ are $1$-Lipschitz for the $L^1$ norm, by \eqref{e.finite.lipschitz} and \eqref{e.scalar.lipschitz}, so $d_k$ at $a$ is the limit of $d_k$ at $a^M$, which lies in $[0,E_n(R)]$. Since $E_n(R)\to0$ as $n\to\infty$, this also proves \eqref{e.scalar.initial.limit}.
\end{proof}

We now return to the system with several species, and compare $\bar F_{N,k}(0,\cdot)$ with the function $\psi_m(q)=\sum_s\lambda_sg_{m_s}(q_s)$ of \eqref{e.def.slice.initial}, evaluated at the magnetization $m=m_{N,k}$ of the slice.

\begin{corollary}[Initial data on slices]\label{c.initial.limit}
For every $R<\infty$, we have
\begin{equation}\label{e.initial.limit}
 \lim_{N\to\infty}\ \sup_k\ \sup_{\substack{q\in\cQ_1^{\sS}\\\|q\|_{L^1}\le R}}
 \bigl|\bar F_{N,k}(0,q)-\psi_{m_{N,k}}(q)\bigr|=0.
\end{equation}
\end{corollary}

\begin{proof}
By \eqref{e.slice.factorization} and the definition \eqref{e.def.slice.initial} of $\psi_m$,
\begin{equation*}
 \bar F_{N,k}(0,q)-\psi_{m_{N,k}}(q)
 =\sum_s\lambda_{N,s}\bigl[G_{n_{N,s},k_s}(q_s)-g_{m_{N,k,s}}(q_s)\bigr]
 +\sum_s(\lambda_{N,s}-\lambda_s)\,g_{m_{N,k,s}}(q_s).
\end{equation*}
Since $\|q_s\|_{L^1}\le\|q\|_{L^1}\le R$ for every $s$, the first sum is bounded in absolute value by $\sum_s\lambda_{N,s}E_{n_{N,s}}(R)$, by Proposition~\ref{p.scalar.initial}, and this tends to zero because $n_{N,s}=N\lambda_{N,s}\to\infty$ for every $s$. For the second, the Lipschitz bounds of Lemma~\ref{l.legendre.m} and Proposition~\ref{p.path.convexity} give $|g_m(q_s)|\le g_m(0)+\|q_s\|_{L^1}\le\log2+R$, so the second sum is at most $(\log2+R)\sum_s|\lambda_{N,s}-\lambda_s|$, which tends to zero because $\lambda_N\to\lambda$.
\end{proof}
%
%
%
%
\section{Convergence of the slice free energies}
\label{s.hj}

In this section, we prove assertions \eqref{e.slice.convergence}--\eqref{e.slice.hopf} of Theorem~\ref{t.slices}: the slice free energies converge, uniformly over the slices, to the viscosity solution $f_m$ of the Hamilton--Jacobi equation \eqref{e.slice.hj} with initial condition $\psi_m$, and $f_m$ is given by the Hopf formula. By the remainder estimate \eqref{e.small.defect}, the enriched free energy is an approximate solution of the equation, and by Corollary~\ref{c.initial.limit}, its initial condition is close to $\psi_{m_{N,k}}$ on bounded sets of paths. A comparison principle with error terms, which also expresses the finite speed of propagation of the equation, turns these two facts into an explicit convergence estimate. We start by recalling the analytic framework.

\subsection{The equation on the closed cone}
\label{s.hj.framework}

We view $\cQ_2^{\sS}$ as a closed convex cone in the Hilbert space $L^2 = L^2([0,1);\R^{\sS})$, with dual cone $(\cQ_2^{\sS})^*$ described in \eqref{e.dual.cone.characterization}. This cone has empty interior in $\cH$, which is the source of the technical difficulties in the definition of viscosity solutions; we only recall the outcome of \cite{chen2022hamilton}, in the form used in \cite[Section~7]{chen2026ising}. The solution on $\cQ_2^{\sS}$ is constructed there as a limit of solutions on finite-dimensional cones, which we also use. Let $L$ be a power of two, let $J_1,\ldots,J_L$ be the dyadic intervals of length $1/L$ partitioning $[0,1)$, and let
\begin{equation}\label{e.def.finite.cone}
 \msf C_L:=\{x\in(\R_+^{\sS})^L:x_1\le\cdots\le x_L\},
 \qquad
 \la x,y\ra_L:=\frac1L\sum_{j=1}^Lx_j\cdot y_j,
 \qquad |x|_L^2:=\la x,x\ra_L,
\end{equation}
which is a closed convex cone in $(\R^{\sS})^L$, with dual cone $\msf C_L^*$. For $x\in(\R^{\sS})^L$ and $p\in\cH$, let
\begin{equation}\label{e.def.projections}
 \iota_Lx:=\sum_jx_j\1_{J_j},
 \qquad
 \pi_Lp:=\Bigl(L\int_{J_j}p(r)\,\d r\Bigr)_{j\le L}
\end{equation}
be the step path with values $x_j$, and the vector of cell averages. These maps are adjoint to each other, $\la\iota_Lx,p\ra_\cH=\la x,\pi_Lp\ra_L$; the map $\iota_L$ is an isometry, and it sends $\msf C_L$ into $\cQ_\infty^{\sS}$; the map $\pi_L$ sends $\cQ_2^{\sS}$ into $\msf C_L$, with $|\pi_Lp|_L\le\|p\|_\cH$ by Jensen's inequality; and $\iota_L$ sends $\msf C_L^*$ into $(\cQ_2^{\sS})^*$, since $\la\iota_Lz,p\ra_\cH=\la z,\pi_Lp\ra_L\ge0$ for $z\in\msf C_L^*$ and $p\in\cQ_2^{\sS}$.

The function $\xi$ is proper on $\R_+^{\sS}$, in the sense that it is coordinatewise increasing and so is $x\mapsto\xi(x+y)-\xi(x)$ for every $y\in\R_+^{\sS}$; this follows from the nonnegativity of its power-series coefficients. It is locally but not globally Lipschitz. By \cite[Lemma~4.4]{chen2022hamilton}, there exists a Lipschitz, proper function $\bar\xi:\R_+^{\sS}\to\R$ that agrees with $\xi$ on the positive part of the closed unit ball; we fix such a regularization and define, for every $\kappa\in\cH$,
\begin{equation}\label{e.def.extended.hamiltonian}
 \msf H(\kappa):=\inf\left\{\int_0^1\bar\xi(p(r))\,\d r:
                   p\in\cQ_2^{\sS}\cap\Ll(\kappa+\Ll(\cQ_2^{\sS}\Rr)^*\Rr)\right\}.
\end{equation}
By \cite[Lemma~4.6]{chen2022hamilton}, $\msf H$ is Lipschitz and $(\cQ_2^{\sS})^*$-increasing on $\cH$, and $\msf H(p)=\int_0^1\bar\xi(p(r))\,\d r$ for $p\in\cQ_2^{\sS}$; in particular, $\msf H(p)=\mcl H(p)$ for $p\in\cQ_2^{\sS}$ with $|p(r)|\le1$ for every $r$, where we recall that $\mcl H(p)=\int_0^1\xi(p(r))\,\d r$ is the nonlinearity in \eqref{e.def.remainder}. We denote by $L_\msf H$ the Lipschitz constant of $\msf H$ for the norm of $\cH$, and we set $\msf H_L(x):=\msf H(\iota_Lx)$ for $x\in(\R^{\sS})^L$, which is Lipschitz with the constant $L_\msf H$ for the norm $|\cdot|_L$, since $\iota_L$ is an isometry, and increasing for the dual cone $\msf C_L^*$, since $\iota_L$ maps $\msf C_L^*$ into $(\cQ_2^{\sS})^*$.

The following definitions apply to both settings: the cone $C=\cQ_2^{\sS}$ in $E=\cH$, with the nonlinearity $H=\msf H$; and the cone $C=\msf C_L$ in $E=(\R^{\sS})^L$, with the inner product $\la\cdot,\cdot\ra_L$ and the nonlinearity $H=\msf H_L$. We write $\|\cdot\|$ for the norm of $E$. We say that a function $\phi:(0,\infty)\times C\to\R$ is smooth if there is a unique $(\dr_t\phi,\nabla\phi)(t,x)\in\R\times E$ such that $\phi(t',x')-\phi(t,x)=\dr_t\phi(t,x)(t'-t)+\la\nabla\phi(t,x),x'-x\ra+O(|t'-t|^2+\|x'-x\|^2)$ as $(t',x')\to(t,x)$ in $(0,\infty)\times C$, and if this differential is continuous; this is \cite[Definition~1.2]{chen2022hamilton}. A continuous function $u$ on $\R_+\times C$ is a viscosity subsolution of
\begin{equation}\label{e.hj.H}
 \dr_tu-H(\nabla u)=0\qquad\text{on }\R_+\times C
\end{equation}
if, whenever a smooth $\phi$ is such that $u-\phi$ has a local maximum at some $(t,x)\in(0,\infty)\times C$, we have $\dr_t\phi(t,x)-H(\nabla\phi(t,x))\le0$; it is a viscosity supersolution if the reverse inequality holds at local minima; and it is a viscosity solution if it is both. These are \cite[Definitions~1.4 and~4.3]{chen2022hamilton}. No boundary condition is imposed at the boundary of the cone: the contacts with test functions are required at all points of $C$, and the monotonicity of $H$ for the dual-cone order is what makes this consistent. A function on $\R_+\times C$ is Lipschitz if it is so for the norm of $\R\times E$. In Theorem~\ref{t.slices}, a viscosity solution of \eqref{e.slice.hj} means a viscosity solution of \eqref{e.hj.H} on $\cQ_2^{\sS}$; by \cite[Remark~4.10]{chen2022hamilton}, the solutions below do not depend on the choice of the regularization $\bar\xi$. We use the following results from \cite{chen2022hamilton}.

\begin{theorem}[Well-posedness, comparison, and Hopf formula \cite{chen2022hamilton}]\label{t.cx}
Let $\psi:\cQ_2^{\sS}\to\R$ be $(\cQ_2^{\sS})^*$-increasing and $1$-Lipschitz for the $L^1$ norm.
\begin{enumerate}
\item There exists a unique Lipschitz viscosity solution $f$ of \eqref{e.hj.H} on $\cQ_2^{\sS}$ with $f(0,\cdot)=\psi$.
\item If $u$ and $v$ are respectively a Lipschitz viscosity subsolution and a Lipschitz viscosity supersolution of \eqref{e.hj.H} on $\cQ_2^{\sS}$, and $u(0,\cdot)\le v(0,\cdot)$ on $\cQ_2^{\sS}$, then $u\le v$ on $\R_+\times\cQ_2^{\sS}$.
\item If $\psi$ is convex, then for every $(t,q)\in\R_+\times\cQ_2^{\sS}$,
\begin{equation}\label{e.hopf.general}
 f(t,q)=\sup_{p\in\cQ_\infty^{\sS}}\inf_{q'\in\cQ_\infty^{\sS}}
 \left\{\psi(q')+\la q-q',p\ra_\cH+t\int_0^1\xi(p(r))\,\d r\right\}.
\end{equation}
\item For every $L$ that is a power of $2$, the equation \eqref{e.hj.H} on $\msf C_L$, with initial condition $x\mapsto\psi(\iota_Lx)$, has a unique Lipschitz viscosity solution $f_L$, the comparison principle of item~(2) holds on $\msf C_L$, and $f_L(t,\pi_Lq)\to f(t,q)$ as $L\to\infty$, for every $(t,q)\in\R_+\times\cQ_2^{\sS}$.
\end{enumerate}
\end{theorem}

Items (1) to (3) are \cite[Theorem~4.7]{chen2022hamilton}, together with \cite[Proposition~3.8 and Remark~4.8]{chen2022hamilton} for the comparison principle, specialized to the present $\R^{\sS}$-valued paths as explained in \cite[Section~7]{chen2026ising}. In item~(4), the existence of a unique Lipschitz viscosity solution on the closed cone $\msf C_L$ and its comparison principle are \cite[Theorem~2.2 and Propositions~2.3 and~2.5]{chen2022hamilton}, applied to the Lipschitz and monotone nonlinearity $\msf H_L$. The dyadic partitions form a good collection of partitions in the sense of \cite[Section~3.1]{chen2022hamilton}, and the convergence of $f_L(t,\pi_Lq)$ to $f(t,q)$ is \cite[Proposition~3.10]{chen2022hamilton}; see also \cite[Theorem~4.7(1) and Remark~4.11]{chen2022hamilton}, where the solution on the cone is constructed as the limit of these projections.

By Corollary~\ref{c.psi.m}, the initial condition $\psi_m$ is $(\cQ_2^{\sS})^*$-increasing, $1$-Lipschitz for the $L^1$ norm, and convex, for every $m\in[-1,1]^{\sS}$. Theorem~\ref{t.cx} therefore yields a unique Lipschitz viscosity solution of \eqref{e.slice.hj}, which we take as the definition of $f_m$, together with its Hopf representation \eqref{e.slice.hopf}. The second assertion of Theorem~\ref{t.slices} amounts to the existence and uniqueness of this solution, and the third to the Hopf formula; both are now established, and it remains to prove the convergence \eqref{e.slice.convergence}, the continuity in $m$, and the representation \eqref{e.slice.martingale}.

We will compare the free energy with the finite-dimensional solutions $f_{m,L}$ of item~(4) of Theorem~\ref{t.cx}, although the former solves the equation only approximately, and although the two initial conditions are close only on a bounded set. The comparison principle handles both errors at once in the following form, which also expresses the finite speed of propagation of the equation.

\begin{lemma}[Comparison with error terms]\label{l.stability}
In either of the two settings above, let $v$ be a Lipschitz viscosity solution of \eqref{e.hj.H}, let $\ep\ge0$, and let $u$ be a Lipschitz function on $\R_+\times C$ such that $u-\ep t$ is a viscosity subsolution and $u+\ep t$ a viscosity supersolution of \eqref{e.hj.H}. Let $\omega:\R_+\to\R_+$ be increasing and Lipschitz, and suppose that $|u(0,x)-v(0,x)|\le\omega(\|x\|)$ for every $x\in C$. Then
\begin{equation}\label{e.stability}
 |u(t,x)-v(t,x)|\le\ep t+\omega\bigl(\|x\|+L_\msf Ht\bigr)
 \qquad\text{for every }(t,x)\in\R_+\times C.
\end{equation}
\end{lemma}

\begin{proof}
We first assume that $\omega$ is smooth with bounded derivative. For $a>0$, let $n_a(x):=(\|x\|^2+a^2)^{1/2}$, which is smooth with $\|\nabla n_a\|\le1$ and $n_a\ge\|x\|$, and let
\begin{equation*}
 \Gamma(t,x):=\omega\bigl(n_a(x)+L_\msf Ht\bigr).
\end{equation*}
The function $\Gamma$ is smooth and Lipschitz, and $\dr_t\Gamma=L_\msf H\,\omega'\ge L_\msf H\|\nabla\Gamma\|$, the derivative of $\omega$ being evaluated at $n_a(x)+L_\msf Ht$. We claim that $v+\Gamma$ is a Lipschitz viscosity supersolution. Let $\phi$ be smooth and such that $(v+\Gamma)-\phi$ has a local minimum at $(t,x)$. Then $\phi-\Gamma$ is smooth and $v-(\phi-\Gamma)$ has a local minimum at $(t,x)$, so $\dr_t\phi-\dr_t\Gamma-H(\nabla\phi-\nabla\Gamma)\ge0$, and the Lipschitz continuity of $H$ gives
\begin{equation*}
 \dr_t\phi-H(\nabla\phi)
 \ge\dr_t\Gamma-\bigl|H(\nabla\phi)-H(\nabla\phi-\nabla\Gamma)\bigr|
 \ge\dr_t\Gamma-L_\msf H\|\nabla\Gamma\|\ge0.
\end{equation*}
In the same way, $v-\Gamma$ is a Lipschitz viscosity subsolution. At $t=0$, since $\omega$ is increasing and $n_a\ge\|x\|$, the assumption gives $u(0,\cdot)\le v(0,\cdot)+\Gamma(0,\cdot)$ and $v(0,\cdot)-\Gamma(0,\cdot)\le u(0,\cdot)$. The comparison principle of Theorem~\ref{t.cx}, applied to the pairs $(u-\ep t,v+\Gamma)$ and $(v-\Gamma,u+\ep t)$, gives $|u-v|\le\ep t+\Gamma$, and letting $a\to0$ yields \eqref{e.stability}. For a general $\omega$, we extend it by $\omega(0)$ to the negative half-line, and we apply the previous case to $\omega_b(s):=\int\omega(s+b-y)\eta_b(y)\,\d y$, where $b>0$ and $\eta_b$ is a smooth probability density supported in $[-b,b]$. Each $\omega_b$ is smooth and increasing, with derivative bounded by the Lipschitz constant of $\omega$, it satisfies $\omega_b\ge\omega$ since $\omega$ is increasing, and $\omega_b\to\omega$ pointwise as $b\to0$; the conclusion for $\omega_b$ passes to the limit.
\end{proof}

\subsection{The free energy as an approximate solution in finite dimensions}
\label{s.hj.finite.comparison}

The remainder estimate \eqref{e.small.defect} applies to regular paths with slopes bounded below. We now show that the enriched free energy, evaluated at regular interpolations of the vectors of $\msf C_L$, is an approximate viscosity solution of the finite-dimensional equation, including at the boundary of the cone. Besides the maps \eqref{e.def.projections}, we use a continuous interpolation. Let $c_j:=(j-\frac12)/L$, for $j\in\{1,\dots,L\}$, be the midpoints of the dyadic intervals. For $x\in(\R^{\sS})^L$, let $\mcl T_Lx$ be the continuous path which is piecewise affine with the values $x_j$ at the midpoints $c_j$, and constant on $[0,c_1]$ and on $[c_L,1]$. Thus
\begin{equation}\label{e.def.interpolation}
 \mcl T_Lx(r)=\sum_jw_j(r)x_j,\qquad w_j\ge0,\qquad\sum_jw_j=1,\qquad\int_0^1w_j(r)\,\d r=\frac1L,
\end{equation}
where $w_j$ is supported in the union of $J_j$ and its neighbors, and the probability densities $Lw_j$ are stochastically ordered from left to right. For $p\in\cH$, let
\begin{equation}\label{e.def.A}
 A_Lp:=\Bigl(L\int_0^1w_j(r)p(r)\,\d r\Bigr)_{j\le L}
\end{equation}
be the vector of averages of $p$ against these densities. The interpolation $\mcl T_L$ produces, after the addition of a positive slope, the regular paths to which the derivative identities and the remainder estimate apply, and $A_L$, which is its adjoint, gives the gradient of the interpolated free energy. For $x\in\msf C_L$, the path $\mcl T_Lx$ belongs to $\cQ_\infty^{\sS}$, and
\begin{equation}\label{e.interpolation.error}
 \|\mcl T_Lx-\iota_Lx\|_{L^1}
 \le\frac1{4L}\sum_{j=1}^{L-1}|x_{j+1}-x_j|
 \le\frac{\sqrt{|\sS|}}{4\sqrt L}\,|x|_L
 \qquad\text{and}\qquad \|\mcl T_Lx\|_{\cH}\le|x|_L,
\end{equation}
the last bound being valid for every $x\in(\R^{\sS})^L$. Indeed, on $[c_j,c_{j+1}]$, the two paths differ by two triangles of total area $|x_{j+1}-x_j|/(4L)$, and $\sum_j|x_{j+1}-x_j|\le\sum_s(x_{L,s}-x_{1,s})\le\sqrt{|\sS|}\,|x_L|\le\sqrt{|\sS|L}\,|x|_L$; and $\mcl T_Lx(r)$ is an average of the $x_j$ with the weights $w_j(r)$, so Jensen's inequality gives $\|\mcl T_Lx\|_\cH^2\le\int_0^1\sum_jw_j(r)|x_j|^2\,\d r=|x|_L^2$. Moreover, $\mcl T_L\pi_Lq\to q$ in $L^1$ as $L\to\infty$, for every $q\in\cQ_2^{\sS}$, since the dyadic averages $\iota_L\pi_Lq$ of an integrable function converge to it in $L^1$, and $|\pi_Lq|_L\le\|q\|_\cH$ is bounded.

\begin{lemma}[The interpolation and its adjoint]\label{l.cone.preservation}
The map $A_L$ sends $\cQ_2^{\sS}$ into $\msf C_L$, and
\begin{equation}\label{e.adjoint}
 \la\mcl T_Lz,p\ra_\cH=\la z,A_Lp\ra_L
 \qquad\text{for all }z\in(\R^{\sS})^L\text{ and }p\in\cH.
\end{equation}
\end{lemma}

\begin{proof}
For $p\in\cQ_2^{\sS}$, we have $(A_Lp)_j=\int p\,\d\mu_j$, where the probability measures $\mu_j:=Lw_j(r)\,\d r$ are stochastically ordered from left to right, so $(A_Lp)_j$ is nonnegative and increasing in $j$. The identity \eqref{e.adjoint} is immediate from \eqref{e.def.interpolation}.
\end{proof}

\begin{proposition}[Approximate solution in finite dimensions]\label{p.approximate}
Fix $N$, a slice $\Omega_{N,k}$, $0<\kappa\le1$, and a power of two $L$. For $x\in\msf C_L$, let
\begin{equation}\label{e.shifted.approximation}
 q^\kappa_x(r):=\mcl T_Lx(r)+\kappa(1+r)\mathbf1,
 \qquad u(t,x):=\bar F_{N,k}(t,q^\kappa_x).
\end{equation}
The function $u$ is Lipschitz on $\R_+\times\msf C_L$, and the functions $u-\ep t$ and $u+\ep t$ are respectively a viscosity subsolution and a viscosity supersolution of \eqref{e.hj.H} on $\msf C_L$, where
\begin{equation}\label{e.finite.residual}
 \ep:=15M_\xi(\kappa N)^{-1/3}+\frac{3L_\xi}{L},
 \qquad \text{with } \ L_\xi:=\sup_{[0,1]^{\sS}}|\nabla\xi|,
\end{equation}
and where we recall that $M_\xi=\sup_{[-1,1]^{\sS}}\|\nabla^2\xi\|$.
\end{proposition}

\begin{proof}
\noindent\emph{Step 1: the residual.}\par
The shift by $\kappa(1+r)\mathbf1$ in \eqref{e.shifted.approximation} has two purposes: it supplies the positive path slope required by Theorem~\ref{t.remainder}, and it makes $u$ differentiable in a neighborhood of the closed cone. Indeed, for $x\in\msf C_L$, every coordinate of $q^\kappa_x$ starts at $x_{1,s}+\kappa\ge\kappa$ and has slope at least $\kappa$. Moreover, $q^\kappa_x$ remains a regular path, with positive initial value and positive slopes, on the open neighborhood of $\msf C_L$ where
\begin{equation*}
 x_{1,s}>-\kappa,\qquad
 x_{j+1,s}-x_{j,s}>-\kappa/L .
\end{equation*}
On this neighborhood, the variations $x\mapsto q^\kappa_x$ are admissible in the sense of Section~\ref{s.derivatives}, and the derivatives of Theorem~\ref{t.derivatives} depend continuously on $(t,x)$ by the parameter dependence noted after Lemma~\ref{l.comparison}. Consequently $u$ is continuously differentiable on this neighborhood for $t>0$, with a right derivative in $t$ at $t=0$. We use the scalar product $\la\cdot,\cdot\ra_L$ to define gradients. By Theorem~\ref{t.derivatives} and \eqref{e.adjoint}, the projected gradient is still an average of expected Gibbs overlaps,
\begin{equation}\label{e.projected.gradient}
 \nabla_Lu=A_Lp_N,\qquad p_N(r)=\E\langle R_N\rangle_r,
 \qquad
 \dr_tu=\int_0^1\E\langle\xi(R_N)\rangle_r\,\d r ,
\end{equation}
where the Gibbs averages are those of the path $q^\kappa_x$. By Lemmas~\ref{l.energy} and \ref{l.cone.preservation}, $A_Lp_N\in\msf C_L$; also $0\le(A_Lp_N)_{j,s}\le\lambda_{N,s}$, so that each $(A_Lp_N)_j$ lies in the positive part of the closed unit ball, and $\msf H_L(A_Lp_N)=L^{-1}\sum_j\xi((A_Lp_N)_j)$. Since $w_j$ is supported in the cell $J_j$ and its neighbors, the $L^1$ distance between $p_N$ and $\iota_LA_Lp_N$ is at most $3/L$ times the sum over $s$ of the total variations of $p_{N,s}$, which is at most $\sum_s\lambda_{N,s}=1$. This gives
\begin{equation}\label{e.quadrature}
 \Ll|\mcl H(p_N)-\msf H_L(\nabla_Lu)\Rr|
 =\Ll|\int_0^1\xi(p_N(r))\,\d r-\frac1L\sum_j\xi\bigl((A_Lp_N)_j\bigr)\Rr|\le\frac{3L_\xi}{L}.
\end{equation}
Combining \eqref{e.defect.identity}, Theorem~\ref{t.remainder} and \eqref{e.quadrature}, we obtain, for every $x\in\msf C_L$ and $t>0$,
\begin{equation}\label{e.residual.bound}
 |\dr_tu-\msf H_L(\nabla_Lu)|\le\ep .
\end{equation}
The function $u$ is Lipschitz in $x$ by \eqref{e.finite.lipschitz} and \eqref{e.interpolation.error}, since $\|\mcl T_L(x-y)\|_{L^1}\le\|\mcl T_L(x-y)\|_\cH\le|x-y|_L$, and in $t$ since $|\dr_t\bar F_{N,k}|\le\xi(\mathbf1)$.

\smallskip
\noindent\emph{Step 2: contacts at the boundary.}\par
Suppose that a smooth function $\phi$ is such that $u-\phi$ has a local maximum at $(t,x)\in(0,\infty)\times\msf C_L$, relative to the closed cone. The extension of $u$ to an open neighborhood gives $\dr_t\phi=\dr_tu$ and, for every $y\in\msf C_L$, the point $x+\theta y$ belongs to $\msf C_L$ for $\theta>0$, so $\la\nabla_L\phi-\nabla_Lu,y\ra_L\ge0$; that is, $\nabla_L\phi-\nabla_Lu\in\msf C_L^*$. Since $\msf H_L$ is increasing for the dual-cone order, \eqref{e.residual.bound} gives $\dr_t\phi-\msf H_L(\nabla_L\phi)\le\dr_tu-\msf H_L(\nabla_Lu)\le\ep$. At a local minimum, the dual-cone relation is reversed and gives $\dr_t\phi-\msf H_L(\nabla_L\phi)\ge-\ep$. Hence $u-\ep t$ is a viscosity subsolution and $u+\ep t$ a viscosity supersolution of the equation on $\msf C_L$. 
\end{proof}

\subsection{Convergence, and continuity in the magnetization}
\label{s.hj.convergence}

For $R>0$, let
\begin{equation}\label{e.def.deltaN}
 \de_N(R):=\sup_k\ \sup_{\substack{q\in\cQ_2^{\sS}\\\|q\|_\cH\le R}}
            |\bar F_{N,k}(0,q)-\psi_{m_{N,k}}(q)|,
\end{equation}
which tends to zero as $N\to\infty$ by Corollary~\ref{c.initial.limit}, since $\|q\|_{L^1}\le\|q\|_\cH$.

\begin{proposition}[Convergence estimate]\label{p.convergence}
For every $N$, every slice $\Omega_{N,k}$, every $t\ge0$, every $q\in\cQ_2^{\sS}$, every $0<\kappa\le1$, and every $R>\|q\|_\cH+L_\msf Ht$,
\begin{equation}\label{e.convergence}
 \bigl|\bar F_{N,k}(t,q)-f_{m_{N,k}}(t,q)\bigr|
 \le\de_N(R)+3\kappa+15M_\xi\,t\,(\kappa N)^{-1/3}.
\end{equation}
\end{proposition}

\begin{proof}
Write $m:=m_{N,k}$, fix a power of two $L$, let $u$ be the function \eqref{e.shifted.approximation}, and let $f_{m,L}$ be the solution of Theorem~\ref{t.cx}(4) with initial condition $x\mapsto\psi_m(\iota_Lx)$. By \eqref{e.finite.lipschitz}, for every $q\in\cQ_1^{\sS}$ and $t\ge0$,
\begin{equation}\label{e.shift.cost}
 \bigl|\bar F_{N,k}\bigl(t,q+\kappa(1+r)\mathbf1\bigr)-\bar F_{N,k}(t,q)\bigr|
 \le\sum_s\lambda_{N,s}\int_0^1\kappa(1+r)\,\d r=\frac32\kappa ,
\end{equation}
which is the cost of the shift in \eqref{e.shifted.approximation}.

\smallskip
\noindent\emph{Step 1: the initial data.}\par
Let $x\in\msf C_L$. By \eqref{e.finite.lipschitz}, \eqref{e.shift.cost} and \eqref{e.interpolation.error},
\begin{equation*}
 |u(0,x)-\bar F_{N,k}(0,\iota_Lx)|
 \le\frac32\kappa+\|\mcl T_Lx-\iota_Lx\|_{L^1}
 \le\frac32\kappa+\frac{\sqrt{|\sS|}}{4\sqrt L}\,|x|_L.
\end{equation*}
If $|x|_L\le R$, then $\|\iota_Lx\|_\cH=|x|_L\le R$, and $|\bar F_{N,k}(0,\iota_Lx)-\psi_m(\iota_Lx)|\le\de_N(R)$ by definition. Otherwise, let $x':=(R/|x|_L)x$, which belongs to $\msf C_L$ and satisfies $\|\iota_L(x-x')\|_{L^1}\le\|\iota_L(x-x')\|_\cH=|x|_L-R$; since $\bar F_{N,k}(0,\cdot)$ and $\psi_m$ are $1$-Lipschitz for the $L^1$ norm, by Proposition~\ref{p.extension} and \eqref{e.initial.lipschitz}, we get $|\bar F_{N,k}(0,\iota_Lx)-\psi_m(\iota_Lx)|\le\de_N(R)+2(|x|_L-R)$. Altogether,
\begin{equation}\label{e.initial.error}
 |u(0,x)-\psi_m(\iota_Lx)|\le\omega(|x|_L),
 \qquad
\text{ where } \  \omega(s):=\de_N(R)+\frac32\kappa+\frac{\sqrt{|\sS|}}{4\sqrt L}\,s+2(s-R)_+ .
\end{equation}

\smallskip
\noindent\emph{Step 2: comparison.}\par
The function $\omega$ is increasing and Lipschitz. By Proposition~\ref{p.approximate} and Lemma~\ref{l.stability}, applied on $\msf C_L$ with $v=f_{m,L}$, we have
\begin{equation}\label{e.comparison.error}
 |u(t,x)-f_{m,L}(t,x)|\le\ep t+\omega\bigl(|x|_L+L_\msf Ht\bigr)
 \qquad\text{for every }(t,x)\in\R_+\times\msf C_L,
\end{equation}
with $\ep$ as in \eqref{e.finite.residual}.

\smallskip
\noindent\emph{Step 3: passage to the limit.}\par
Let $x:=\pi_Lq$. Since $|x|_L\le\|q\|_\cH$ and $R>\|q\|_\cH+L_\msf Ht$, the last term of $\omega$ vanishes in \eqref{e.comparison.error}, which gives
\begin{equation*}
 |u(t,\pi_Lq)-f_{m,L}(t,\pi_Lq)|
 \le\ep t+\de_N(R)+\frac32\kappa+\frac{\sqrt{|\sS|}}{4\sqrt L}\bigl(\|q\|_\cH+L_\msf Ht\bigr).
\end{equation*}
By \eqref{e.finite.lipschitz} and \eqref{e.shift.cost}, $|\bar F_{N,k}(t,q)-u(t,\pi_Lq)|\le\|q-\mcl T_L\pi_Lq\|_{L^1}+\frac32\kappa$. We now let $L\to\infty$ along the powers of two: $\|q-\mcl T_L\pi_Lq\|_{L^1}\to0$, the term $3L_\xi t/L$ in $\ep t$ and the last term of the previous display vanish, and $f_{m,L}(t,\pi_Lq)\to f_m(t,q)$ by Theorem~\ref{t.cx}(4). This proves \eqref{e.convergence}.
\end{proof}

\begin{proof}[Proof of \eqref{e.slice.convergence} in Theorem~\ref{t.slices}, and continuity in $m$]
Let $t\ge0$ and $q\in\cQ_2^{\sS}$, and choose $R>\|q\|_\cH+L_\msf Ht$. Proposition~\ref{p.convergence} with $\kappa=N^{-1/4}$ gives
\begin{equation*}
 \sup_k\bigl|\bar F_{N,k}(t,q)-f_{m_{N,k}}(t,q)\bigr|
 \le\de_N(R)+(3+15M_\xi t)N^{-1/4}\xrightarrow[N\to\infty]{}0,
\end{equation*}
which is \eqref{e.slice.convergence}.

It remains to show that $m\mapsto f_m(t,q)$ is continuous. Let $m,m'\in[-1,1]^{\sS}$, and let $c:=\mcl I_\lambda(m)-\mcl I_\lambda(m')$, with $\mcl I_\lambda$ as in Corollary~\ref{c.psi.m}. The function $\psi_{m'}+c$ satisfies the assumptions of Theorem~\ref{t.cx}, and its solution is $f_{m'}+c$. By \eqref{e.initial.magnetization.lipschitz}, for every $q'\in\cQ_2^{\sS}$,
\begin{equation*}
 |\psi_m(q')-(\psi_{m'}(q')+c)|
 \le2\sum_s\lambda_s\|q'_s\|_{L^1}|m_s-m'_s|\le2\|q'\|_\cH\,|m-m'|_\infty,
\end{equation*}
where $|m-m'|_\infty:=\max_s|m_s-m'_s|$. Lemma~\ref{l.stability} on $\cQ_2^{\sS}$, with $u=f_m$, $v=f_{m'}+c$, $\ep=0$ and $\omega(s)=2s|m-m'|_\infty$, gives
\begin{equation}\label{e.f.magnetization.lipschitz}
 \bigl|f_m(t,q)-f_{m'}(t,q)-\mcl I_\lambda(m)+\mcl I_\lambda(m')\bigr|
 \le2\bigl(\|q\|_\cH+L_\msf Ht\bigr)|m-m'|_\infty .
\end{equation}
Since $\mcl I_\lambda$ is continuous on $[-1,1]^{\sS}$, so is $m\mapsto f_m(t,q)$.
\end{proof}

The estimate \eqref{e.f.magnetization.lipschitz} is more precise than the continuity asserted in Theorem~\ref{t.slices}: the entropy $\mcl I_\lambda(m)$ carries all the singular dependence of the slice free energy on the magnetization, and the remainder is Lipschitz in $m$, uniformly over bounded sets of $(t,q)$.
%
%
%
%
\section{From the Hopf formula to martingales}\label{s.duality}

In this section, we prove the martingale representation \eqref{e.slice.martingale} of $f_m(t,q)$, starting from the Hopf formula \eqref{e.slice.hopf}. The two convexity properties of $g_m$ established in Section~\ref{s.scalar} play complementary roles. The convexity in the path gives the Hopf formula, and we first rewrite this formula in conjugate form, in terms of the path conjugate of $g_m$ (Section~\ref{s.duality.conjugate}). The concavity in the law, in the form of the martingale representation of Proposition~\ref{p.scalar.martingale}, expresses $g_m$ as an infimum of affine functionals of the law, parametrized by the potentials $\chi_\al$. We then show that the Hopf formula can be written as an infimum over such potentials (Sections~\ref{s.duality.potentials} to~\ref{s.duality.minimax}), and we replace the potentials by martingales (Section~\ref{s.duality.transfer}). The passage from the Hopf formula to the potentials requires exchanging a supremum and an infimum; this is done with Sion's minimax theorem, after enlarging the set of paths to a compact set of couplings. The enlargement does not change the value, by a rearrangement inequality that uses the nonnegativity of the coefficients of $\xi$.

Three different conjugacies appear in this section, and we fix the notation for them now. For increasing paths, we use the path conjugates $\psi_m^*$ and $J_{m}$ of Section~\ref{s.duality.conjugate}, which pair a path with an increasing slope. We give the path conjugate $J_m$ of $g_m$ its own letter, rather than an upper star, because it is conjugated once more in Section~\ref{s.duality.scalar}, for a different pairing. For laws, we pair a probability measure $\nu$ with a potential $\chi$ through $\int\chi\,\d\nu$, and we represent the concave functional $g_m$ by its supporting potentials; the corresponding concave conjugate is denoted by $g_*$ below, with a lower star. Finally, for functions of a finite-dimensional variable, such as a potential $\chi$ on $\R_+^{\sS}$ or a function $v$ on $[0,1]$, we write $\chi^*$ and $v^*$ for the usual Legendre transforms. Upper stars always denote convex conjugates. The two scalar identities of Section~\ref{s.duality.scalar} relate the first two conjugacies through the third one, and they also show that the path conjugate $J_m$ is convex for the affine interpolation of the laws of slopes.

\subsection{The conjugate form of the Hopf formula}
\label{s.duality.conjugate}

It is convenient to rewrite the Hopf formula \eqref{e.slice.hopf} in conjugate form, and to record the duality in the path variable that comes with the supporting paths of Proposition~\ref{p.path.convexity}. For $m\in[-1,1]$ and $z\in\cQ_\infty$ with $0\le z\le1$, we define the path conjugate
\begin{equation}\label{e.def.J}
 J_m(z):=\sup_{a\in\cQ_\infty}
 \left\{\int_0^1z(r)a(r)\,\d r-g_m(a)\right\}\in[-g_m(0),+\infty].
\end{equation}

\begin{lemma}[Path biconjugacy]\label{l.path.bidual}
For every $m\in[-1,1]$ and $a\in\cQ_\infty$, the path $z_a$ of \eqref{e.def.z.a} satisfies $J_m(z_a)=\int_0^1z_aa-g_m(a)$, and
\begin{equation}\label{e.scalar.path.bidual}
 g_m(a)=\sup_{\substack{z\in\cQ_\infty\\0\le z\le1}}
       \left\{\int_0^1z(r)a(r)\,\d r-J_m(z)\right\},
\end{equation}
the supremum being attained at $z=z_a$.
\end{lemma}

\begin{proof}
By \eqref{e.scalar.supporting}, $\int_0^1z_ab-g_m(b)\le\int_0^1z_aa-g_m(a)$ for every $b\in\cQ_\infty$, with equality for $b=a$; this is the first claim. The inequality $\ge$ in \eqref{e.scalar.path.bidual} holds for every admissible $z$ by the definition of $J_m$, and the first claim gives equality at $z=z_a$.
\end{proof}

In words, $g_m$ is the conjugate of its own conjugate, and the dual variable may be restricted to increasing paths with values in $[0,1]$. 
For $p\in\cQ_\infty^{\sS}$, let
\begin{equation}\label{e.def.path.conjugate}
 \psi_m^*(p):=\sup_{q'\in\cQ_\infty^{\sS}}
                 \{\la p,q'\ra_{\cH}-\psi_m(q')\}\in[-\psi_m(0),+\infty].
\end{equation}
Recall the set $\cQ_{\infty,\le\lambda}^{\sS}$ of paths with $0\le p_s\le\lambda_s$ from \eqref{e.def.Q.lambda}.

\begin{lemma}[Conjugate form of the Hopf formula]\label{l.hopf.conjugate}
For every $m\in[-1,1]^{\sS}$, $t\ge0$ and $q\in\cQ_2^{\sS}$, we have
\begin{equation}\label{e.hopf.conjugate}
 f_m(t,q)=\sup_{p\in\cQ_{\infty,\le\lambda}^{\sS}}
 \left\{\la q,p\ra_{\cH}+t\int_0^1\xi(p(r))\,\d r-\psi_m^*(p)\right\}.
\end{equation}
Moreover, for every $p\in\cQ_{\infty,\le\lambda}^{\sS}$, we have
\begin{equation}\label{e.conjugate.separable}
 \psi_m^*(p)=\sum_s\lambda_s\,J_{m_s}\Bigl(\frac{p_s}{\lambda_s}\Bigr),
\end{equation}
where $J_{m_s}$ is the path conjugate \eqref{e.def.J}.
\end{lemma}

\begin{proof}
The Hopf formula \eqref{e.slice.hopf} reads $f_m(t,q)=\sup_{p\in\cQ_\infty^{\sS}}\{\la q,p\ra_\cH+t\int_0^1\xi(p) -\psi_m^*(p)\}$. If $p\in\cQ_\infty^{\sS}\setminus\cQ_{\infty,\le\lambda}^{\sS}$, there are $s\in\sS$ and $r_0\in[0,1)$ with $p_s(r_0)>\lambda_s$, by right-continuity, so $\int_{r_0}^1(p_s-\lambda_s)>0$ by monotonicity. Testing \eqref{e.def.path.conjugate} with $q'=a\1_{[r_0,1)}e_s$ and using \eqref{e.initial.lipschitz}, we get $\psi_m^*(p)\ge a\int_{r_0}^1(p_s-\lambda_s)-\psi_m(0)\to+\infty$ as $a\to\infty$. Such $p$ do not contribute to the supremum. For $p\in\cQ_{\infty,\le\lambda}^{\sS}$, the separable form \eqref{e.conjugate.separable} follows from \eqref{e.def.slice.initial}, since the components of $q'$ vary independently; the lower bound on $\psi_m^*$ is obtained by taking $q'=0$.
\end{proof}

At $t=0$, where $f_m(0,\cdot)=\psi_m$, formula \eqref{e.hopf.conjugate} reduces to the biconjugacy of Lemma~\ref{l.path.bidual}, summed over the species.

\subsection{Potentials, and localization to the horizon}
\label{s.duality.potentials}

Throughout the rest of this section, we fix $m\in[-1,1]^{\sS}$, $t\ge0$ and $q\in\cQ_\infty^{\sS}$, and horizons $T=(T_s)_{s\in\sS}$ such that
\begin{equation}\label{e.horizon}
 T_s>0,\qquad T_s\ge\|q_s\|_\infty+t\,\dr_s\xi(\lambda)
 \qquad\text{for every }s\in\sS.
\end{equation}
We write $[0,T]:=\prod_s[0,T_s]$ and $\Lambda_\lambda:=\prod_s[0,\lambda_s]$.

Let $\cX_1$ be the set of convex, increasing, $1$-Lipschitz functions $\chi:\R_+\to\R$ with $\chi(0)=0$, and let
\begin{equation}
\label{e.def.cX}
 \cX:=\Bigl\{x\mapsto\sum_s\lambda_s\chi_s(x_s):\chi_s\in\cX_1\Bigr\}
\end{equation}
be the corresponding separable functions on $\R_+^{\sS}$. We often call elements of $\cX$ potentials. For $\chi\in\cX$, we define its conjugate and its evolution at time $t$ by
\begin{equation}\label{e.def.potential.evolution}
 \chi^*(p):=\sup_{x\in\R_+^{\sS}}\{p\cdot x-\chi(x)\},\qquad
 \msf S_t\chi(x):=\sup_{p\in\R_+^{\sS}}
                        \{p\cdot x+t\xi(p)-\chi^*(p)\}.
\end{equation}
If $p_s>\lambda_s$ for some $s$, then $\chi^*(p)=+\infty$, since $\chi$ has slope at most $\lambda_s$ in the direction $e_s$; so the supremum defining $\msf S_t\chi$ is over $\Lambda_\lambda$. For a functional $g(\nu)=\sum_s\lambda_sg_s(\nu_s)$ of the marginals of $\nu\in\Prob_1(\R_+^{\sS})$, we define its conjugate for the pairing with potentials, and the value associated with the pair $(g,q)$, by
\begin{align}
 g_*(\chi)&:=\inf_{\nu\in\Prob_1(\R_+^{\sS})}
                   \left\{\int\chi\,\d\nu-g(\nu)\right\},
                                                    \label{e.def.g.conjugate}\\
 \mcl U_t[g](q)&:=\inf_{\chi\in\cX}
           \left\{\int_0^1 \msf S_t\chi(q(r))\,\d r-g_*(\chi)\right\}.
                                                    \label{e.def.support.value}
\end{align}
Since $\chi$ and $g$ are separable, and since every family of laws in $\Prob_1(\R_+)$ is the family of marginals of some $\nu\in\Prob_1(\R_+^{\sS})$, the infimum in \eqref{e.def.g.conjugate} decouples over the species:
\begin{equation}\label{e.g.conjugate.species}
 g_*(\chi)=\sum_s\lambda_s\,(g_s)_*(\chi_s),\qquad
 (g_s)_*(\chi_s):=\inf_{\nu_s\in\Prob_1(\R_+)}
                   \left\{\int\chi_s\,\d\nu_s-g_s(\nu_s)\right\}.
\end{equation}
Our goal is to prove that $f_m(t,q)=\mcl U_t[\psi_m](q)$, where $\psi_m$ is regarded as the functional $\nu\mapsto\sum_s\lambda_sg_{m_s}(\nu_s)$ of laws.
The next lemma shows that, thanks to the horizon condition \eqref{e.horizon}, the problem \eqref{e.def.support.value} can be localized to the box $[0,T]$.

\begin{lemma}[Localization]\label{l.localization}
Let $g(\nu)=\sum_s\lambda_sg_s(\nu_s)$ be a functional of the marginals, as above.
\begin{enumerate}
\item For every $\chi\in\cX$ and every $x\in\R_+^{\sS}$ with $x_s+t\dr_s\xi(\lambda)\le T_s$ for all $s$, we have
\begin{equation}\label{e.localized.evolution}
 \msf S_t\chi(x)=\sup_{p\in\Lambda_\lambda}
 \left\{p\cdot x+t\xi(p)
              -\sup_{y\in [0,T]}(p\cdot y-\chi(y))\right\}.
\end{equation}
In particular, $\msf S_t\chi(x)$ depends only on the restriction of $\chi$ to $[0,T]$; and if $\chi\le\chi'+c$ on $[0,T]$ for a constant $c$, then $\msf S_t\chi(x)\le \msf S_t\chi'(x)+c$.
\item Suppose that each $g_s$ is $1$-Lipschitz for the $1$-Wasserstein distance. Then the infimum in \eqref{e.def.support.value} may be restricted to the potentials $\chi=\sum_s\lambda_s\chi_s$ such that $\chi_s$ has slope $1$ on $[T_s,\infty)$ for every $s$, and for these potentials, we have
\begin{equation}\label{e.localized.measure.conjugate}
 g_*(\chi)=\inf_{\nu\in\Prob([0,T])}
                         \left\{\int\chi\,\d\nu-g(\nu)\right\}.
\end{equation}
\end{enumerate}
\end{lemma}

\begin{proof}
\noindent\emph{Step 1: the evolution.}\par
Write $\chi^*_T(p):=\sup_{y\in [0,T]}(p\cdot y-\chi(y))$, which is a finite convex function of $p$, separable in the coordinates. Since $\chi^*_T\le\chi^*$, the right-hand side of \eqref{e.localized.evolution} is at least $\msf S_t\chi(x)$. For the converse, let $c_s:=\lambda_s\chi_{s}'(T_s-)$ be the left derivative of $\lambda_s\chi_s$ at $T_s$. If $p_s\ge c_s$, the concave function $y_s\mapsto p_sy_s-\lambda_s\chi_s(y_s)$ is increasing on $[0,T_s]$, so its supremum over $[0,T_s]$ is attained at $y_s=T_s$, and $\chi^*_T$ is affine in $p_s$ with slope $T_s$ on this range. There, the derivative in $p_s$ of the objective in \eqref{e.localized.evolution} is $x_s+t\dr_s\xi(p)-T_s\le x_s+t\dr_s\xi(\lambda)-T_s\le0$, because $\dr_s\xi$ is coordinatewise increasing on $\R_+^{\sS}$. Hence decreasing $p_s$ to $c_s$ does not decrease the objective. After doing this in every coordinate, we obtain a point $p$ with $p_s\le c_s$ for all $s$; for such $p$, the function $y_s\mapsto p_sy_s-\lambda_s\chi_s(y_s)$ is decreasing on $[T_s,\infty)$, so $\chi^*(p)=\chi^*_T(p)$. The right-hand side of \eqref{e.localized.evolution} is therefore at most $\msf S_t\chi(x)$. The supremum is attained by continuity and compactness. The last two assertions follow, since $\chi^*_T$ only involves $\chi$ on $[0,T]$, and $\chi\le\chi'+c$ on $[0,T]$ implies $\chi^*_T\ge\chi'^*_T-c$.

\smallskip
\noindent\emph{Step 2: the measure conjugate.}\par
Let $\chi\in\cX$, and let $\hat\chi_s$ coincide with $\chi_s$ on $[0,T_s]$ and have slope $1$ on $[T_s,\infty)$. Since $\chi_s$ is $1$-Lipschitz, $\hat\chi_s\ge\chi_s$, so $g_*(\hat\chi)\ge g_*(\chi)$, while $\msf S_t\hat\chi(q(r))=\msf S_t\chi(q(r))$ for every $r$ by Step~1 and \eqref{e.horizon}. Thus replacing $\chi$ by $\hat\chi$ does not increase the objective in \eqref{e.def.support.value}. Now let $\chi$ have slope $1$ beyond the horizons, let $\nu\in\Prob_1(\R_+^{\sS})$, and let $\nu^T$ be its image under the coordinatewise truncation $x\mapsto x\wedge T$. For each $s$, $\int\chi_s\,\d\nu_s-\int\chi_s\,\d\nu^T_s=\int(x_s-T_s)_+\,\d\nu_s =W_1(\nu_s,\nu_s^T)$, and $g_s(\nu_s)\le g_s(\nu_s^T)+W_1(\nu_s,\nu_s^T)$. Hence $\int\chi\,\d\nu^T-g(\nu^T)\le\int\chi\,\d\nu-g(\nu)$. Finally, $\Prob([0,T])$ is compact for the weak topology, and $\nu\mapsto\int\chi\,\d\nu-g(\nu)$ is continuous on it, since $\chi$ is continuous on $[0,T]$ and each $g_s$ is continuous for the weak topology on $\Prob([0,T_s])$. This proves \eqref{e.localized.measure.conjugate}.
\end{proof}

By Lemma~\ref{l.legendre.m}, the functions $g_m$ are $1$-Lipschitz for $W_1$, so the lemma applies to $\psi_m$. The supporting inequality \eqref{e.scalar.support} also yields a representation of $g_m$ by potentials. Here $(g_m)_*(\chi)=\inf_{\nu\in\Prob_1(\R_+)}\{\int\chi\,\d\nu-g_m(\nu)\}$ is the scalar conjugate appearing in \eqref{e.g.conjugate.species}.

\begin{lemma}[Dual representation of $g_m$]\label{l.scalar.bidual}
For every $m\in[-1,1]$, $T>0$ and $\nu\in\Prob([0,T])$,
\begin{equation}\label{e.scalar.measure.bidual}
 g_m(\nu)=\inf_{\chi\in\cX_1}
                         \left\{\int\chi\,\d\nu-(g_m)_*(\chi)\right\},
\end{equation}
and the infimum is attained at the potential $\hat\chi_{\nu,m}\in\cX_1$ that coincides with the supporting potential $\chi_{\nu,m}$ of Corollary~\ref{c.scalar.support} on $[0,T]$ and has slope $1$ on $[T,\infty)$.
\end{lemma}


\begin{proof}
The inequality $\le$ in \eqref{e.scalar.measure.bidual} holds for every $\chi\in\cX_1$, by the definition of $(g_m)_*$. Conversely, let $\nu'\in\Prob_1(\R_+)$, and let $\nu'^T$ be its image under $x\mapsto x\wedge T$. By Lemma~\ref{l.legendre.m}, $g_m(\nu')\le g_m(\nu'^T)+W_1(\nu',\nu'^T)=g_m(\nu'^T)+\int(x-T)_+\,\d\nu'$, and by \eqref{e.scalar.support} at $\nu'^T\in\Prob([0,T])$, $g_m(\nu'^T)\le g_m(\nu)+\int\chi_{\nu,m}\,\d(\nu'^T-\nu)$. Since $\hat\chi_{\nu,m}(x)=\chi_{\nu,m}(x\wedge T)+(x-T)_+$ and $\nu$ is supported in $[0,T]$, the two bounds combine into $g_m(\nu')\le g_m(\nu)+\int\hat\chi_{\nu,m}\,\d\nu'-\int\hat\chi_{\nu,m}\,\d\nu$. Taking the infimum over $\nu'$, we get $(g_m)_*(\hat\chi_{\nu,m})\ge\int\hat\chi_{\nu,m}\,\d\nu-g_m(\nu)$, which is the reverse inequality at $\chi=\hat\chi_{\nu,m}$.
\end{proof}

\subsection{Two scalar conjugacies}
\label{s.duality.scalar}

We fix a scalar magnetization $m_0\in[-1,1]$ and abbreviate $g:=g_{m_0}$ and $J:=J_{m_0}$, where $J$ is the path conjugate \eqref{e.def.J}. For $\rho\in\Prob([0,1])$, we write $J(\rho):=J(p_\rho)$, where $p_\rho$ is the quantile function of $\rho$; recall that the quantile functions are exactly the right-continuous increasing functions on $[0,1)$ with values in $[0,1]$, so that $J(\rho)$ is defined for every $\rho$ and every such $z$ is $p_\rho$ for the law $\rho$ of $z(U)$. The function $g$ is convex in the path, with conjugate $J$; and $g$ is concave in the law, with concave conjugate $g_*$ for the pairing with potentials. The following lemma relates the two conjugates. Let $\cV$ be the convex set of convex, increasing, Lipschitz functions $v:[0,1]\to\R$ with $v(0)=0$. For $v\in\cV$, its conjugate $v^*(x):=\sup_{0\le z\le1}(xz-v(z))$, $x\ge0$, is convex, increasing, $1$-Lipschitz, vanishes at $0$, and has slope $1$ beyond the Lipschitz constant of $v$; conversely, if $\chi\in\cX_1$ has slope $1$ beyond some point, then $\chi^*|_{[0,1]}\in\cV$ and $(\chi^*)^*=\chi$. Thus $v\mapsto v^*$ is a bijection from $\cV$ onto the potentials in $\cX_1$ with eventual slope $1$, and $v^{**}=v$ on $[0,1]$. We also regard $J$ as a function on $\Prob([0,1])$, and we define its convex conjugate for the pairing $\int v\,\d\rho$ between functions and measures on $[0,1]$:
\begin{equation}\label{e.def.J.star}
 J^*(v):=\sup_{\rho\in\Prob([0,1])}
       \left\{\int v\,\d\rho-J(\rho)\right\},\qquad v\in\cV.
\end{equation}
The following lemma identifies $J^*$ with the concave conjugate $g_*$, and shows that $J$ is its own biconjugate. 

\begin{lemma}[Scalar dual identities]\label{l.scalar.conjugacies}
For every $v\in\cV$ and every $\rho\in\Prob([0,1])$,
\begin{align}
 J^*(v)&=-g_*(v^*),                                  \label{e.J.star.g}\\
 J(\rho)&=\sup_{v\in\cV}\left\{\int v\,\d\rho-J^*(v)\right\}.
                                                     \label{e.J.bidual}
\end{align}
In particular, $J^*$ is finite, convex and $1$-Lipschitz for the uniform norm on $\cV$.
\end{lemma}

\begin{proof}
\noindent\emph{Step 1: a pointwise conjugacy.}\par
Let $\chi\in\cX_1$ have slope $1$ beyond $x_1\ge0$, and let $z\in\cQ_\infty$ with $0\le z\le1$. We claim that
\begin{equation}\label{e.pointwise.conjugacy}
 \sup_{a\in\cQ_\infty}\int_0^1[z(r)a(r)-\chi(a(r))]\,\d r
 =\int_0^1\chi^*(z(r))\,\d r.
\end{equation}
The inequality $\le$ is clear. For the converse, let $a(r)$ be the smallest maximizer of $x\mapsto z(r)x-\chi(x)$ on $\R_+$; it exists and belongs to $[0,x_1]$, since this function is concave and decreasing beyond $x_1$, and it is increasing in $r$ because the function has increasing differences in $(x,z)$ and $z$ is increasing. Replacing $a$ by its right-continuous version, which differs from it on a countable set only, we obtain an element of $\cQ_\infty$ attaining the supremum.

\smallskip
\noindent\emph{Step 2: a first formula for $J$.}\par
Let $z:=p_\rho$ be the quantile function of $\rho$. Inserting \eqref{e.scalar.measure.bidual} into the definition of $J$, where by Lemma~\ref{l.scalar.bidual} the infimum may be restricted to potentials with eventual slope $1$, and exchanging the two suprema,
\begin{equation*}
 J(\rho)=\sup_{a\in\cQ_\infty}\sup_\chi
 \left\{\int_0^1[z(r)a(r)-\chi(a(r))]\,\d r+g_*(\chi)\right\}
 =\sup_\chi\left\{\int_0^1\chi^*(z(r))\,\d r+g_*(\chi)\right\},
\end{equation*}
by \eqref{e.pointwise.conjugacy}. Writing $v:=\chi^*|_{[0,1]}\in\cV$, we have $\chi=v^*$ and $\int_0^1v(z(r))\,\d r=\int v\,\d\rho$, so that
\begin{equation}\label{e.J.first}
 J(\rho)=\sup_{v\in\cV}\left\{\int v\,\d\rho+g_*(v^*)\right\}.
\end{equation}

\smallskip
\noindent\emph{Step 3: proof of \eqref{e.J.star.g} and \eqref{e.J.bidual}.}\par
Since $v^*$ has slope $1$ beyond the Lipschitz constant of $v$, the truncation argument of Lemma~\ref{l.localization} shows that the infimum defining $g_*(v^*)$ may be taken over compactly supported laws, that is, $-g_*(v^*)=\sup_{a\in\cQ_\infty}\{g(a)-\int_0^1v^*(a(r))\,\d r\}$. Inserting the biconjugacy \eqref{e.scalar.path.bidual} of Lemma~\ref{l.path.bidual} and exchanging the suprema,
\begin{multline*}
 -g_*(v^*)=\sup_{\substack{z\in\cQ_\infty\\0\le z\le1}}
 \left\{\sup_{a\in\cQ_\infty}\int_0^1[z(r)a(r)-v^*(a(r))]\,\d r-J(z)\right\}\\
 =\sup_z\left\{\int_0^1v(z(r))\,\d r-J(z)\right\},
\end{multline*}
by \eqref{e.pointwise.conjugacy} applied to $\chi=v^*$, whose conjugate is $v$. Since the paths $z$ in the supremum are exactly the quantile functions of the laws $\rho\in\Prob([0,1])$, the right-hand side is $J^*(v)$, which proves \eqref{e.J.star.g}. Inserting \eqref{e.J.star.g} into \eqref{e.J.first} gives \eqref{e.J.bidual}. Finally, $J^*$ is finite by \eqref{e.J.star.g}, since $g(\nu)\le g(\delta_0)+\int x\,\d\nu$ by the Lipschitz bound of Lemma~\ref{l.legendre.m} and $x-v(1)\le v^*(x)\le x$, so that $-v(1)-g(\delta_0)\le g_*(v^*)\le-g(\delta_0)$; and \eqref{e.def.J.star} represents $J^*$ as a supremum of affine functions of $v$, each $1$-Lipschitz for the uniform norm.
\end{proof}

The identity \eqref{e.J.bidual} will be used to express the Hopf value through potentials, and the convexity of $J^*$ will be used when appealing to the minimax theorem. The lemma also shows that $J$, which was defined as a conjugate in the path variable, is convex for the affine interpolation of the laws $\rho$, as a supremum of affine functions of $\rho$.

\subsection{Minimax and a common increasing rearrangement}
\label{s.duality.minimax}

We now express the Hopf value through potentials.

\begin{proposition}[The Hopf value through potentials]
\label{p.hopf.support}
For every $m\in[-1,1]^{\sS}$, $t\ge0$ and $q\in\cQ_\infty^{\sS}$, we have
\begin{equation}\label{e.hopf.support}
 f_m(t,q)=\mcl U_t[\psi_m](q).
\end{equation}
\end{proposition}

\begin{proof}
\noindent\emph{Step 1: from increasing paths to couplings.}\par
For each $s\in\sS$, let $J_s$ and $J_s^*$ be the functions of Section~\ref{s.duality.scalar} with magnetization $m_s$, and for $z\in[0,1]^{\sS}$, let $\lambda z:=(\lambda_sz_s)_s$. By Lemma~\ref{l.hopf.conjugate}, writing $p=\lambda z$ with $z$ an increasing path with values in $[0,1]^{\sS}$,
\begin{equation*}
 f_m(t,q)=\sup_z\left\{\int_0^1\bigl[q(r)\cdot\lambda z(r)
                            +t\xi(\lambda z(r))\bigr]\,\d r
                      -\sum_s\lambda_sJ_s(z_s)\right\}.
\end{equation*}
Let $\mu_q$ be the law of $q(U)$ on $\R_+^{\sS}$, and let $\Pi_q$ be the set of probability measures $\pi$ on $\R_+^{\sS}\times[0,1]^{\sS}$ whose first marginal is $\mu_q$. Since $q$ is bounded, every $\pi\in\Pi_q$ is supported in the compact set $K\times[0,1]^{\sS}$, where $K:=\prod_s[0,\|q_s\|_\infty]$; hence $\Pi_q$ is a convex, weakly closed subset of $\Prob(K\times[0,1]^{\sS})$, and thus weakly compact, and the integrands below are bounded and continuous on $K\times[0,1]^{\sS}$. We write $(x,z)$ for a generic point of $\R_+^{\sS}\times[0,1]^{\sS}$, with $x=(x_s)_{s\in\sS}$ and $z=(z_s)_{s\in\sS}$, and for $\pi\in\Pi_q$ and $s\in\sS$, we write $\pi_s\in\Prob([0,1])$ for the image of $\pi$ under the map $(x,z)\mapsto z_s$. For an increasing path $z$ as in the last display, let $\pi^z\in\Pi_q$ be the law of $(q(U),z(U))$; we call these couplings monotone. The marginal $\pi^z_s$ is the law of $z_s(U)$, so that $J_s(\pi^z_s)=J_s(z_s)$, and the integral above is $\int[x\cdot\lambda z+t\xi(\lambda z)]\,\d\pi^z$. We claim that
\begin{equation}\label{e.coupling.formula}
 f_m(t,q)=\sup_{\pi\in\Pi_q}
 \left\{\int[x\cdot\lambda z+t\xi(\lambda z)]\,\d\pi
                             -\sum_s\lambda_sJ_s(\pi_s)\right\},
\end{equation}
that is, that the supremum over all couplings is already attained on the monotone ones. Let $\pi\in\Pi_q$, with marginals $\mu_q$ and $(\rho_s)_s$ for the coordinates $z_s$, and let $\pi^\uparrow:=\pi^z$ for $z:=(p_{\rho_s})_s$, where $p_{\rho_s}$ is the quantile function of $\rho_s$. This monotone coupling has the same marginals as $\pi$, so the costs $J_s$ are unchanged, and it suffices to show that the first integral in \eqref{e.coupling.formula} does not decrease when $\pi$ is replaced by $\pi^\uparrow$. For nonnegative random variables $Y_1,\ldots,Y_n$ with prescribed marginal laws, the layer-cake formula gives
\begin{equation*}
 \E\Bigl[\prod_iY_i\Bigr]
 =\int_{\R_+^n}\P(Y_1>u_1,\ldots,Y_n>u_n)\,\d u
 \le\int_{\R_+^n}\min_i\P(Y_i>u_i)\,\d u,
\end{equation*}
and the monotone coupling $Y_i=p_{\Law(Y_i)}(U)$ achieves the bound, since for it the joint tail probability equals the smallest marginal tail probability. We apply this to each product $x_sz_s$ and to each monomial $\prod_sz_s^{a_s}$ in the power series of $\xi(\lambda z)$. Since the coefficients of this series are nonnegative, and since it converges absolutely, summing the inequalities proves the claim.

\smallskip
\noindent\emph{Step 2: exchange of the supremum and the infimum.}\par
By \eqref{e.J.bidual}, $J_s(\pi_s)=\sup_{v_s\in\cV}\{\int v_s\,\d\pi_s-J_s^*(v_s)\}$, so that \eqref{e.coupling.formula} becomes
\begin{equation}\label{e.minimax}
 f_m(t,q)=\sup_{\pi\in\Pi_q}\ \inf_{v\in\cV^{\sS}}
 \Bigl\{\int[x\cdot\lambda z+t\xi(\lambda z)]\,\d\pi
          -\sum_s\lambda_s\int v_s\,\d\pi_s
                          +\sum_s\lambda_sJ_s^*(v_s)\Bigr\}.
\end{equation}
The expression in braces is affine in $\pi$, and weakly continuous since its integrands are bounded and continuous on $K\times[0,1]^{\sS}$; it is convex and continuous for the uniform norm in $v$, by Lemma~\ref{l.scalar.conjugacies}. Since $\Pi_q$ is convex and compact, Sion's minimax theorem \cite{sion1958general} allows us to exchange the supremum and the infimum.

\smallskip
\noindent\emph{Step 3: the supremum over couplings at a fixed potential.}\par
Let $v=(v_s)_{s\in\sS}\in\cV^{\sS}$, and let $\chi_v:=\sum_s\lambda_sv_s^*\in\cX$. Then, for $z\in[0,1]^{\sS}$,
\begin{equation}\label{e.chi.v.conjugate}
 \chi_v^*(\lambda z)=\sum_s\lambda_sv_s(z_s),
 \qquad
 -(\psi_m)_*(\chi_v)=\sum_s\lambda_sJ_s^*(v_s).
\end{equation}
The first identity follows from $v_s^{**}=v_s$, and the second from \eqref{e.g.conjugate.species} and \eqref{e.J.star.g}. The supremum over $\pi\in\Pi_q$ of the integral of a continuous function of $(x,z)$ is the integral against $\mu_q$ of its pointwise supremum over $z\in[0,1]^{\sS}$. By the first identity in \eqref{e.chi.v.conjugate}, and since $\chi_v^*=+\infty$ outside $\Lambda_\lambda$, the definition \eqref{e.def.potential.evolution} gives
\begin{equation*}
 \sup_{\pi\in\Pi_q}\int\Bigl[x\cdot\lambda z+t\xi(\lambda z)
                  -\sum_s\lambda_sv_s(z_s)\Bigr]\,\d\pi(x,z)
 =\int_0^1\msf S_t\chi_v(q(r))\,\d r.
\end{equation*}
Together with the second identity in \eqref{e.chi.v.conjugate}, the exchanged form of \eqref{e.minimax} therefore reads
\begin{equation*}
 f_m(t,q)=\inf_{v\in\cV^{\sS}}
 \left\{\int_0^1\msf S_t\chi_v(q(r))\,\d r-(\psi_m)_*(\chi_v)\right\}.
\end{equation*}
By Lemma~\ref{l.localization}, the infimum defining $\mcl U_t[\psi_m](q)$ may be restricted to the potentials with slope $1$ beyond the horizons, and these are of the form $\chi_v$ with $v\in\cV^{\sS}$. The infimum over all $v\in\cV^{\sS}$ lies between this restricted infimum and the infimum over all of $\cX$, so it equals $\mcl U_t[\psi_m](q)$. This proves \eqref{e.hopf.support}.
\end{proof}

\subsection{Replacing potentials by martingales}
\label{s.duality.transfer}

Recall from the introduction that $\cM_T$ is the set of families $\al=(\al_s)_{s\in\sS}$ of $[-1,1]$-valued martingales on $[0,T_s]$, $\al_s$ being a martingale for the filtration $\cF^s$ of the Brownian motion $B_s$, and that $\cM_T^m$ is the subset of those with $\al_s(0)=m_s$ for every $s$. For $\al\in\cM_T$, we define the potential
\begin{equation}\label{e.def.chi.alpha}
 \chi_\al(x):=\sum_s\lambda_s\chi_{\al_s}(x_s)
 =\sum_s\lambda_s\int_0^{x_s}\E[\al_s(v\wedge T_s)^2]\,\d v,
 \qquad x\in\R_+^{\sS},
\end{equation}
with $\chi_{\al_s}$ as in \eqref{e.def.chi.c}. It belongs to $\cX$. The cost $\cost_T(\al)$ of the introduction is $\sum_s\lambda_sc_{T_s}(\al_s)$, with $c_T$ as in \eqref{e.def.chi.c}.

\begin{proposition}[Transfer and attainment]\label{p.transfer}
For every $m$, $t$, $q$ and $T$ as in \eqref{e.horizon},
\begin{equation}\label{e.transfer}
 \mcl U_t[\psi_m](q)=\inf_{\al\in\cM_T^m}
       \left\{\cost_T(\al)+\int_0^1\msf S_t\chi_\al(q(r))\,\d r\right\},
\end{equation}
and the infimum is attained.
\end{proposition}

\begin{proof}
The compactness used for the attainment is that of the potentials, not of the martingales: we first find an optimal potential, and we then construct a martingale from a law attaining the concave conjugate at this potential.

\smallskip
\noindent\emph{Step 1: an optimal potential.}\par
By Lemma~\ref{l.localization}(2), the infimum defining $\mcl U_t[\psi_m](q)$ may be restricted to the potentials $\chi=\sum_s\lambda_s\chi_s$ with $\chi_s\in\cX_1$ of slope $1$ beyond $T_s$. Such a potential is determined by its restriction to $[0,T]$, and these restrictions form a compact set for the uniform norm, by the Arzel\`a--Ascoli theorem, since they are $1$-Lipschitz, vanish at $0$, and form a closed class. Both terms of the objective in \eqref{e.def.support.value} are $1$-Lipschitz for the uniform norm on $[0,T]$, by \eqref{e.localized.evolution} and \eqref{e.localized.measure.conjugate}. An optimal potential $\chi$ therefore exists; we fix it for the rest of the proof.

\smallskip
\noindent\emph{Step 2: from the optimal potential to a martingale.}\par
By \eqref{e.g.conjugate.species} and Lemma~\ref{l.localization}(2), $(\psi_m)_*(\chi)=\sum_s\lambda_s(g_{m_s})_*(\chi_s)$, and each $(g_{m_s})_*(\chi_s)$ is an infimum over $\Prob([0,T_s])$, attained at some $\nu_s$. Let $\al_s:=M^{\nu_s,m_s}$ be the minimizing martingale of Proposition~\ref{p.scalar.martingale} on $[0,T_s]$, so that $\al\in\cM_T^m$ and $g_{m_s}(\nu_s)=c_{T_s}(\al_s)+\int\chi_{\al_s}\,\d\nu_s$. We claim that
\begin{equation}\label{e.support.optimality}
 \chi_s(x)-\chi_{\al_s}(x)\ge\vartheta_s
 \qquad\text{for every }x\in[0,T_s],\qquad
 \vartheta_s:=\int(\chi_s-\chi_{\al_s})\,\d\nu_s.
\end{equation}
This is the first-order optimality condition for $\nu_s$. Indeed, let $x\in[0,T_s]$, $\ep\in(0,1)$, and $\nu':=(1-\ep)\nu_s+\ep\delta_x\in\Prob([0,T_s])$. By the minimality of $\nu_s$,
\begin{equation*}
 \frac{g_{m_s}(\nu')-g_{m_s}(\nu_s)}{\ep}
 \le\frac1\ep\int\chi_s\,\d(\nu'-\nu_s)
 =\chi_s(x)-\int\chi_s\,\d\nu_s,
\end{equation*}
and letting $\ep\to0$ with \eqref{e.scalar.derivative.law} turns the left-hand side into $\chi_{\al_s}(x)-\int\chi_{\al_s}\,\d\nu_s$. In words, $\chi_s$ touches $\chi_{\al_s}+\vartheta_s$ from above on $[0,T_s]$, and the two integrate to the same value against~$\nu_s$, so $\nu_s$ is supported on the contact set. Consequently,
\begin{equation*}
 (\psi_m)_*(\chi)=\sum_s\lambda_s\Bigl[\int\chi_s\,\d\nu_s-g_{m_s}(\nu_s)\Bigr]
 =\sum_s\lambda_s\vartheta_s-\cost_T(\al),
\end{equation*}
and, since $\chi\ge\chi_\al+\sum_s\lambda_s\vartheta_s$ on $[0,T]$, the ordering property of Lemma~\ref{l.localization}(1) gives $\msf S_t\chi(q(r))\ge \msf S_t\chi_\al(q(r))+\sum_s\lambda_s\vartheta_s$ for every $r$. Therefore
\begin{equation}\label{e.transfer.lower}
 \mcl U_t[\psi_m](q)=\int_0^1\msf S_t\chi(q(r))\,\d r-(\psi_m)_*(\chi)
 \ge\cost_T(\al)+\int_0^1\msf S_t\chi_\al(q(r))\,\d r.
\end{equation}

\smallskip
\noindent\emph{Step 3: every martingale gives a potential.}\par
Let $\al'\in\cM_T^m$ be arbitrary, and let $\hat\chi_{\al'}:=\sum_s\lambda_s\hat\chi_{\al'_s}$, where $\hat\chi_{\al'_s}(x):=\chi_{\al'_s}(x\wedge T_s)+(x-T_s)_+$ coincides with $\chi_{\al'_s}$ on $[0,T_s]$ and has slope $1$ beyond $T_s$. Let $\nu\in\Prob_1(\R_+^{\sS})$, and let $\nu^T$ be its image under the coordinatewise truncation $x\mapsto x\wedge T$. By Lemma~\ref{l.legendre.m}, $g_{m_s}(\nu_s)\le g_{m_s}(\nu^T_s)+\int(x_s-T_s)_+\,\d\nu_s$, and by the inequality in \eqref{e.scalar.martingale} at $\nu^T_s\in\Prob([0,T_s])$, $g_{m_s}(\nu^T_s)\le c_{T_s}(\al'_s)+\int\chi_{\al'_s}\,\d\nu^T_s$. Since $\int\chi_{\al'_s}\,\d\nu^T_s+\int(x_s-T_s)_+\,\d\nu_s=\int\hat\chi_{\al'_s}\,\d\nu_s$, this gives $g_{m_s}(\nu_s)\le c_{T_s}(\al'_s)+\int\hat\chi_{\al'_s}\,\d\nu_s$, hence $\psi_m(\nu)\le\cost_T(\al')+\int\hat\chi_{\al'}\,\d\nu$ and $(\psi_m)_*(\hat\chi_{\al'})\ge-\cost_T(\al')$. Moreover, $\msf S_t\hat\chi_{\al'}(q(r))=\msf S_t\chi_{\al'}(q(r))$ by Lemma~\ref{l.localization}(1). Thus
\begin{equation*}
 \mcl U_t[\psi_m](q)\le\int_0^1\msf S_t\hat\chi_{\al'}(q(r))\,\d r-(\psi_m)_*(\hat\chi_{\al'})
 \le\cost_T(\al')+\int_0^1\msf S_t\chi_{\al'}(q(r))\,\d r.
\end{equation*}
Together with \eqref{e.transfer.lower}, this proves \eqref{e.transfer}, with the infimum attained at the martingale $\al$ of Step~2.
\end{proof}

The minimizer constructed in Step~2 has a special form: each $\al_s$ is the martingale~$M^{\nu_s,m_s}$ of Proposition~\ref{p.scalar.martingale}, for some law $\nu_s\in\Prob([0,T_s])$. The infimum in \eqref{e.transfer}, and hence in \eqref{e.slice.martingale}, may therefore always be sought among these martingales.

\begin{proof}[Proof of \eqref{e.slice.martingale} in Theorem~\ref{t.slices}]
Let $q\in\cQ_\infty^{\sS}$ and let $T$ satisfy \eqref{e.horizon}. We claim that $\msf S_t\chi_\al(q(r))=\mcl K^{T}_{t,\al}(q(r))$ for every $\al\in\cM_T$ and $r\in[0,1]$. Since $q_s(r)\le\|q_s\|_\infty$, the point $x=q(r)$ satisfies the condition of Lemma~\ref{l.localization}(1), so $\msf S_t\chi_\al(q(r))$ is given by \eqref{e.localized.evolution}, which only involves the conjugate of $\chi_\al$ truncated to $[0,T]$. This truncated conjugate is explicit: for an increasing function $\varrho:[0,T_s]\to[0,1]$ and $z\in[0,1]$,
\begin{equation}\label{e.explicit.conjugate}
 \sup_{0\le y\le T_s}\left\{zy-\int_0^y\varrho(v)\,\d v\right\}
                             =\int_0^{T_s}(z-\varrho(v))_+\,\d v,
\end{equation}
since the function of $y$ on the left increases while $\varrho<z$ and decreases afterwards. Applying this in each coordinate of \eqref{e.localized.evolution}, with $p=\lambda z$ and $\varrho(v)=\E[\al_s(v)^2]$, which is increasing since $\al_s$ is a martingale, and comparing with \eqref{e.def.K}, proves the claim. Combining it with Propositions~\ref{p.hopf.support} and~\ref{p.transfer} gives
\begin{equation*}
 f_m(t,q)=\inf_{\al\in\cM_T^m}
 \left\{\cost_T(\al)+\int_0^1\mcl K^{T}_{t,\al}(q(r))\,\d r\right\},
\end{equation*}
with the infimum attained. This completes the proof of Theorem~\ref{t.slices}.
\end{proof}

%
%
%
%
\section{Restoring the external fields}\label{s.conclusion}

We now sum over the slices and prove Theorem~\ref{t.main}. Recall that $\bar F_N(t,q)$ denotes the enriched free energy of the full cube with the weights $2^{-N}e^{h^{(N)}\cdot\si}$, which coincides with \eqref{e.def.free.energy} at $q=0$ by Lemma~\ref{l.zero}. The gluing lemma expresses it through the slice free energies, and Theorem~\ref{t.slices} identifies their limits.

\begin{theorem}[Enriched free energy with external fields]\label{t.enriched}
Let $t\ge0$ and $h\in\R^{\sS}$. For every $q\in\cQ_2^{\sS}$, the identity \eqref{e.full.slices} holds, and the infimum there is attained. For every $q\in\cQ_\infty^{\sS}$ and every $T$ satisfying \eqref{e.horizon}, we have
\begin{equation}\label{e.main.enriched}
 \lim_{N\to\infty}\bar F_N(t,q)
 =\inf_{\al\in\cM_T}
 \left\{\cost_T(\al)-\sum_s\lambda_sh_s\al_s(0)
                       +\int_0^1\mcl K^{T}_{t,\al}(q(r))\,\d r\right\},
\end{equation}
and the infimum is attained.
\end{theorem}

In the case $q=0$, the condition \eqref{e.horizon} reads $T_s>0$ and $T_s\ge t\dr_s\xi(\lambda)$, and \eqref{e.main.enriched} reduces to \eqref{e.main.martingale}. Theorem~\ref{t.main} is thus the special case $q=0$ of Theorem~\ref{t.enriched}.

\begin{proof}[Proof of Theorem~\ref{t.enriched}]
We start with the proof of \eqref{e.full.slices}. 
Let $t\ge0$, $h\in\R^{\sS}$ and $q\in\cQ_2^{\sS}$. By Lemma~\ref{l.gluing},
\begin{equation*}
 \bar F_N(t,q)
 -\min_k\bigl\{\bar F_{N,k}(t,q)-\lambda_N\cdot(hm_{N,k})\bigr\}
 \longrightarrow0.
\end{equation*}
By Theorem~\ref{t.slices}, $\sup_k|\bar F_{N,k}(t,q)-f_{m_{N,k}}(t,q)|\to0$, the function $m\mapsto f_m(t,q)$ is continuous on $[-1,1]^{\sS}$, the set of compositions $\{m_{N,k}\}_k$ becomes dense in $[-1,1]^{\sS}$ as $N\to\infty$, and $\lambda_N\to\lambda$. Hence
\begin{equation*}
 \lim_{N\to\infty}\min_k\bigl\{\bar F_{N,k}(t,q)-\lambda_N\cdot(hm_{N,k})\bigr\}
 =\inf_{m\in[-1,1]^{\sS}}\bigl\{f_m(t,q)-\lambda\cdot(hm)\bigr\},
\end{equation*}
where the infimum on the right is attained by continuity and compactness. This shows~\eqref{e.full.slices}.

Turning now to the proof of \eqref{e.main.enriched}, we let $t\ge0$, $h\in\R^{\sS}$ and $q\in\cQ_\infty^{\sS}$, and let $T$ satisfy~\eqref{e.horizon}. By \eqref{e.slice.martingale}, we have for every $m\in[-1,1]^{\sS}$ that
\begin{equation*}
 f_m(t,q)-\lambda\cdot(hm)
 =\inf_{\al\in\cM_T^m}\left\{
 \cost_T(\al)-\sum_s\lambda_sh_s\al_s(0)
                       +\int_0^1\mcl K^{T}_{t,\al}(q(r))\,\d r\right\},
\end{equation*}
since $\al_s(0)=m_s$ on $\cM_T^m$. The sets $\cM_T^m$, $m\in[-1,1]^{\sS}$, form a partition of $\cM_T$, because the initial value of a martingale in $\cM_T$ is deterministic. Minimizing over $m$ in \eqref{e.full.slices} therefore gives \eqref{e.main.enriched}, and the infimum is attained: choose a minimizing $m$ in \eqref{e.full.slices}, and then a minimizing martingale in $\cM_T^m$.
\end{proof}

We conclude with a few comments.

At $h=0$, the functions $m\mapsto g_m(\nu)$ are even and convex, so $\psi_0\le\psi_m$ for every $m$, and the comparison principle gives $f_0\le f_m$. Thus $m=0$ is a minimizing slice in \eqref{e.full.slices}, and the formula reduces to the result of \cite{chen2026ising} for centered spins: the free energy of the centered model is the slice free energy at zero magnetization, which is the solution of the Hamilton--Jacobi equation with the centered initial condition $\psi_0$, and this is consistent with $g_0=\psi_\circ(\cdot;0)$.

For general $h$, the limit free energy of the unconstrained system also admits a representation by potentials, as conjectured in \cite{issa2024hopflike} and shown in \cite[Theorem~1.4]{chen2025convex} for convex models. To see this, we set
\begin{equation*}
 \psi^h(q):=\sum_s\lambda_s\psi_\circ(q_s;h_s),
 \qquad q\in\cQ_1^{\sS},
\end{equation*}
which is the limit of $\bar F_N(0,q)$ by \eqref{e.factorization}. We regard $\psi^h$ as a functional of laws, as in Section~\ref{s.duality.potentials}, and we recall the definition of the lower-star operator from \eqref{e.def.g.conjugate}, of $\cX$ from \eqref{e.def.cX}, and of $\msf S_t$ from \eqref{e.def.potential.evolution}.

\begin{theorem}[Representation by potentials]\label{t.potential}
Let $t\ge0$ and $h\in\R^{\sS}$. For every $q\in\cQ_\infty^{\sS}$, 
\begin{equation}\label{e.potential.formula}
 \lim_{N\to\infty}\bar F_N(t,q)=\mcl U_t[\psi^h](q)
 =\inf_{\chi\in\cX}\left\{
             \int_0^1\msf S_t\chi(q(r))\,\d r-(\psi^h)_*(\chi)\right\}.
\end{equation}
\end{theorem}

\begin{proof}
Let $t\ge0$, $h\in\R^{\sS}$ and $q\in\cQ_\infty^{\sS}$. By Lemma~\ref{l.legendre.m} applied in each coordinate,
\begin{equation*}
 \psi^h(\nu)=\inf_{m\in[-1,1]^{\sS}}
                       \{\psi_m(\nu)-\lambda\cdot(hm)\}
\end{equation*}
for every law $\nu$. Each $\psi_\circ(\cdot;h_s)$ is $1$-Lipschitz for $W_1$ by Lemma~\ref{l.scalar.continuity}, so Lemma~\ref{l.localization} applies to $\psi^h$: the infimum in \eqref{e.potential.formula} may be restricted to potentials $\chi$ with slope $1$ beyond the horizons $T$ of \eqref{e.horizon}, and for these potentials,
\begin{equation*}
 -(\psi^h)_*(\chi)
 =\sup_{\nu\in\Prob([0,T])}\ \inf_{m\in[-1,1]^{\sS}}
 \Bigl\{\psi_m(\nu)-\lambda\cdot(hm)-\int\chi\,\d\nu\Bigr\}.
\end{equation*}
The function in braces is continuous on the compact convex set $\Prob([0,T])\times[-1,1]^{\sS}$, by the Lipschitz bounds in $\nu$ and in $m$ of Lemma~\ref{l.legendre.m}; it is concave in $\nu$ by Corollary~\ref{c.scalar.support}, and convex in $m$ by Lemma~\ref{l.legendre.m}. By Sion's minimax theorem, the supremum and the infimum may be exchanged, and Lemma~\ref{l.localization} applied to $\psi_m$ then gives
\begin{equation*}
 -(\psi^h)_*(\chi)=\inf_{m\in[-1,1]^{\sS}}
                       \{-(\psi_m)_*(\chi)-\lambda\cdot(hm)\}.
\end{equation*}
Taking the infimum over $\chi$ and exchanging the two infima, we obtain, by Proposition~\ref{p.hopf.support},
\begin{equation*}
 \mcl U_t[\psi^h](q)=\inf_{m\in[-1,1]^{\sS}}
                    \{\mcl U_t[\psi_m](q)-\lambda\cdot(hm)\}
 =\inf_{m\in[-1,1]^{\sS}}\{f_m(t,q)-\lambda\cdot(hm)\},
\end{equation*}
and \eqref{e.potential.formula} follows from \eqref{e.full.slices}.
\end{proof}

%
%
%
%
\appendix
\section{The cascade interpretation of the branch}\label{s.cascade}

This appendix relates the enriched free energy of Definition~\ref{d.free.energy} to its usual description through a Poisson--Dirichlet cascade; its results are not used elsewhere in the paper. We prove three statements. First, $\bar F_{N,\Omega}$ agrees with the cascade free energy on all of $\cQ_1^{\sS}$ (Proposition~\ref{p.cascade}). Second, in a cascade with finitely many levels, suppose that we sample a leaf and a configuration from the Gibbs measure; as the number of levels tends to infinity, the joint law of the fields along this leaf and of this configuration converges in total variation to that of the values of the branch at the levels and of its sampled configuration (Proposition~\ref{p.selected.limit}). Third, concerning the Gibbs measures at level $r$ of Section~\ref{s.gibbs}, the law of two configurations under $\E\langle\cdot\rangle_r$ is the limit of the law of two replicas of the cascade Gibbs measure, conditioned on their leaves separating at level $r$ (Lemma~\ref{l.selected.pair} and Proposition~\ref{p.cascade.pair}). The agreement of the free energies does not by itself imply the other claims, and it is the second statement that justifies calling $X$ a branch of the cascade. The third statement makes explicit the sense in which the derivative identities of Theorem~\ref{t.derivatives} are the familiar cascade computation, disintegrated according to the level at which two replicas separate.

Throughout, $(\Omega,w)$ is as in Section~\ref{s.branch}, $t\ge0$, and $q$ is a regular path. We write $\ga_q$ for the law of the Gaussian path $Z_r=X_0+\int_0^r\sqrt{2D(u)}\,\d B_u$ of Section~\ref{s.concentration.law}, so that the measure $\ga$ of Lemma~\ref{l.joint.density} is the product of the law of $g$ and of $\ga_q$.

\subsection{The frozen-coefficient problem}
\label{s.cascade.frozen}

Let $n\ge1$, $r_j:=j/n$, and let $\rho_n(r):=r_j$ for $r\in[r_j,r_{j+1})$, with $\rho_n(1):=1$. Let $U^{(n)}$ be the solution of \eqref{e.branch.pde} with $\rho_n(r)$ in place of the coefficient $r$, and with the terminal condition \eqref{e.def.terminal}. Lemmas~\ref{l.wellposed} and~\ref{l.comparison} apply to this equation, with the same proofs, the pieces being now delimited by the breakpoints of $q$ and the mesh points; in particular, $|\dr_{x_i}U^{(n)}|\le1$ for every $i$. For $r\in[0,1]$ and $a\in\R_+^{\sS}$, we define the Gaussian operator $T_{r,a}$ on functions of $x\in\R^N$ with at most linear growth by
\begin{equation}\label{e.def.vector.operator}
\begin{gathered}
 T_{r,a}f(x):=\frac1r\log\E_z\exp\bigl(rf(x+\sqrt{2\diag(a)}\,z)\bigr)
 \quad(r>0),\\
 T_{0,a}f(x):=\E_zf(x+\sqrt{2\diag(a)}\,z),
\end{gathered}
\end{equation}
where $z$ is a standard Gaussian vector in $\R^N$. This is the $N$-dimensional version of the scalar operator \eqref{e.def.gaussian.operator}, with the same semigroup property $T_{r,a}T_{r,a'}=T_{r,a+a'}$.

\begin{lemma}[Cole--Hopf recursion]\label{l.recursion}
For every $x\in\R^N$ and $0\le j\le n-1$,
\begin{equation}\label{e.recursion}
 U^{(n)}(r_j,x;g)=T_{r_j,\,q(r_{j+1})-q(r_j)}U^{(n)}(r_{j+1},\cdot\,;g)(x),
\end{equation}
with $U^{(n)}(1,\cdot\,;g)$ the terminal function \eqref{e.def.terminal}. Moreover, let $X^{(n)}$ solve \eqref{e.branch.sde} with $\rho_n(r)$ in place of the coefficient $r$ and $U^{(n)}$ in place of $U$, driven by $X_0$ and $B$. For every realization of $g$, the law of $X^{(n)}$ has the density
\begin{equation}\label{e.frozen.density}
\begin{split}
 \Lambda_n(x;g)
 &:=\prod_{j=1}^{n-1}\exp\Bigl(r_j\bigl[U^{(n)}(r_{j+1},x_{r_{j+1}};g)-U^{(n)}(r_j,x_{r_j};g)\bigr]\Bigr)\\
 &\;=\exp\Bigl(\Bigl(1-\frac1n\Bigr)U(1,x_1;g)-\frac1n\sum_{j=1}^{n-1}U^{(n)}(r_j,x_{r_j};g)\Bigr)
\end{split}
\end{equation}
at the path $x$ with respect to $\ga_q$. In particular, the values $(X^{(n)}_{r_1},\ldots,X^{(n)}_{r_n})$ form a Markov chain, whose initial value $X^{(n)}_{r_1}$ is a centered Gaussian vector with covariance $2\diag(q(r_1))$, and in which, for $1\le j\le n-1$ and given $X^{(n)}_{r_j}=x$, the increment $X^{(n)}_{r_{j+1}}-x$ has density
\begin{equation}\label{e.selected.transition}
 \exp\Bigl(r_j\bigl[U^{(n)}(r_{j+1},x+z;g)-U^{(n)}(r_j,x;g)\bigr]\Bigr)
\end{equation}
at $z$ with respect to the centered Gaussian law with covariance $2\diag(q(r_{j+1})-q(r_j))$. Conditionally on these values, $X^{(n)}$ has the conditional law of the Gaussian path $Z$ given the same values.
\end{lemma}

\begin{proof}
On $[r_j,r_{j+1}]$ the coefficient is the constant $r_j$, and for $j\ge1$ the function $\exp(r_jU^{(n)})$ solves the backward heat equation $\dr_rw+\tr(D\nabla^2w)=0$, whose solution is the Gaussian convolution with covariance $2\int_{r_j}^{r_{j+1}}D=2\diag(q(r_{j+1})-q(r_j))$. For $j=0$ the equation for $U^{(n)}$ itself is this heat equation. This proves \eqref{e.recursion}.

For the density, we repeat the Girsanov computation of the proof of Lemma~\ref{l.joint.density} with the coefficient $\rho_n$. Given $g$, the frozen branch has the bounded drift $2\rho_nD\nabla U^{(n)}$ relative to the Gaussian path, so its law has density
\begin{equation*}
 \exp\Ll(\int_0^1\rho_n(r)\,\nabla U^{(n)}\cdot\sqrt{2D}\,\d\tilde B_r
          -\int_0^1\rho_n(r)^2\,\nabla U^{(n)}\cdot D\nabla U^{(n)}\,\d r\Rr)
\end{equation*}
with respect to $\ga_q$, where $\d x_r=\sqrt{2D(r)}\,\d\tilde B_r$ under $\ga_q$, and where $U^{(n)}$ is evaluated at $(r,x_r;g)$. Under $\ga_q$, It\^o's formula and the frozen equation give $\d U^{(n)}(r,x_r;g)=-\rho_n(r)\,\nabla U^{(n)}\cdot D\nabla U^{(n)}\,\d r+\nabla U^{(n)}\cdot\sqrt{2D}\,\d\tilde B_r$ on each interval between consecutive mesh points and breakpoints of $q$, and $U^{(n)}$ is continuous across these points. Since $\rho_n=r_j$ on $[r_j,r_{j+1})$ and $r_0=0$, the exponent equals
\begin{equation*}
 \sum_{j=1}^{n-1}r_j\bigl[U^{(n)}(r_{j+1},x_{r_{j+1}};g)-U^{(n)}(r_j,x_{r_j};g)\bigr],
\end{equation*}
which is the first expression in \eqref{e.frozen.density}; the second follows by summation by parts, using $r_j-r_{j-1}=1/n$ and $U^{(n)}(1,\cdot\,;g)=U(1,\cdot\,;g)$. Under $\ga_q$, the value $x_{r_1}=X_0+(Z_{r_1}-Z_0)$ is a centered Gaussian vector with covariance $2\diag(q(r_1))$, and the increments $x_{r_{j+1}}-x_{r_j}$, $1\le j\le n-1$, are independent of it and of each other, centered Gaussian with covariances $2\diag(q(r_{j+1})-q(r_j))$. The density $\Lambda_n$ depends on the path only through its values at the mesh points, and it is a product of functions of two consecutive values, each of which integrates to one in the later value by \eqref{e.recursion}. The two remaining assertions follow.
\end{proof}

The recursion \eqref{e.recursion} is the one that computes the free energy of a cascade with $n-1$ levels, and the transition \eqref{e.selected.transition} is the law of the increment of the field along a leaf of this cascade selected according to the Gibbs weights; this is the content of Lemma~\ref{l.selected.leaf} below. We first record the analytic estimate that controls the effect of freezing the coefficient $r$.

\begin{lemma}[Convergence of the frozen problem]
\label{l.frozen.convergence}
Let $\|\dot q\|_\infty$ denote the largest slope of $q$ over all pieces and species. For every $n\ge1$ and every realization of $g$,
\begin{equation}\label{e.frozen.convergence}
 \sup_{[0,1]\times\R^N}\bigl|U^{(n)}-U\bigr|\le\frac{N\|\dot q\|_\infty}{n}.
\end{equation}
\end{lemma}

\begin{proof}
Since $U^{(n)}$ solves the frozen equation, the operator \eqref{e.def.branch.operator} gives
\begin{equation*}
 \mcl L[U^{(n)}]=\bigl(r-\rho_n(r)\bigr)\,\nabla U^{(n)}\cdot D(r)\nabla U^{(n)}
 \in\Bigl[0,\frac{N\|\dot q\|_\infty}{n}\Bigr],
\end{equation*}
because $0\le r-\rho_n(r)\le1/n$ and $|\nabla U^{(n)}|^2\le N$, while $\mcl L[U]=0$ and the two functions have the same terminal datum. The estimate \eqref{e.comparison} of Lemma~\ref{l.comparison}, applied on the pieces delimited by the breakpoints of $q$ and the mesh points, gives \eqref{e.frozen.convergence}.
\end{proof}

\subsection{Finite-level cascades and the selected leaf}
\label{s.cascade.leaf}

Consider a step path $\bar q$ taking the value $q_j\in\R^{\sS}_+$ on $[r_{j-1},r_j)$, $1\le j\le n$, where $0=r_0<r_1<\cdots<r_n=1$ is an arbitrary partition and $q_1\le\cdots\le q_n$ coordinatewise. Let $(\msf v_\al)_{\al\in\N^{n-1}}$ be a Poisson--Dirichlet cascade with parameters $r_1<\cdots<r_{n-1}$, independent of $g$: for $n=1$, it has a single leaf $\al=\emptyset$ of weight one, and in general the weights are the normalized products, along the branches of the tree, of the points of independent Poisson--Dirichlet point processes, with parameter $r_j$ at depth $j$; see \cite[Section~2.3]{pan} or \cite[Section~5.6]{HJbook}. For a leaf $\al=(\al_1,\ldots,\al_{n-1})$ and $0\le j\le n-1$, let $\al_{|j}:=(\al_1,\ldots,\al_j)$ denote its ancestor at depth $j$, with $\al_{|0}=\emptyset$ the root. We attach to every vertex $\be$ of the tree $\bigcup_{j=0}^{n-1}\N^j$ an independent standard Gaussian vector~$z_\be$ of $\R^N$, independent of the weights and of $g$ and which we call a mark, and we define the field at a vertex $\be$ of depth $j$ by
\begin{equation}\label{e.cascade.fields}
 x_\be:=\sqrt{2\diag(q_1)}\,z_\emptyset
      +\sum_{\ell=1}^{j}\sqrt{2\diag(q_{\ell+1}-q_\ell)}\,z_{\be_{|\ell}}.
\end{equation}
Thus the fields $(x_{\al_{|0}},\ldots,x_{\al_{|n-1}})$ along the branch leading to the leaf $\al$ have independent centered Gaussian increments, $x_{\al_{|j}}$ has covariance $2\diag(q_{j+1})$, and the fields $x_\al,x_{\al'}$ at two leaves whose most recent common ancestor has depth $j$ have covariance $2\diag(q_{j+1})=2\diag(\bar q(r_j))$. The cascade free energy at $\bar q$ is
\begin{equation}\label{e.cascade.definition}
\begin{split}
 \bar F^{\mathrm{cas}}_{N,\Omega}(t,\bar q)
 &:=t\xi(\lambda_N)+\lambda_N\cdot q_n
  -\frac1N\E\log\sum_\al\msf v_\al\sum_{\si\in\Omega}w(\si)
                 e^{\sqrt{2t}\,g\cdot v_\si+x_\al\cdot\si}\\
 &\;=t\xi(\lambda_N)+\lambda_N\cdot q_n
  -\frac1N\E\log\sum_\al\msf v_\al e^{U(1,x_\al;g)}.
\end{split}
\end{equation}
For the full cube without external field, this is the enriched free energy of \cite{chen2026ising} at the step path $\bar q$. There, the Poisson--Dirichlet cascade has a continuum of levels, with uniform overlap distribution, and the field is $\sqrt2$ times a Gaussian field indexed by the cascade with covariance $N\bar q(\al\wedge\al')\cdot R_N(\si,\si')$; for a step path, this continuous cascade may be replaced by the cascade with $n-1$ levels just described without changing the free energy, by \cite[Proposition~4.6]{chen2025free}. For two step paths $\bar q,\bar q'$, we have
\begin{equation}\label{e.cascade.lipschitz}
 \bigl|\bar F^{\mathrm{cas}}_{N,\Omega}(t,\bar q)-\bar F^{\mathrm{cas}}_{N,\Omega}(t,\bar q')\bigr|
 \le\sum_{s\in\sS}\lambda_{N,s}\|\bar q_s-\bar q'_s\|_{L^1},
\end{equation}
by Gaussian integration by parts in the cascade, exactly as in \cite[Proposition~5.1]{chen2025free}: the derivative of the free energy in the direction of a path is the expected Gibbs average of the overlap of two replicas, evaluated at the level where their leaves separate, and each species overlap lies in $[-\lambda_{N,s},\lambda_{N,s}]$; the argument does not depend on the set $\Omega$ or on the weights $w$. Hence $\bar F^{\mathrm{cas}}_{N,\Omega}(t,\cdot)$ extends by continuity to $\cQ_1^{\sS}$, and, for the full cube without external field, the extension is the enriched free energy of \cite{chen2026ising}, which is defined on $\cQ_1^{\sS}$ by the same continuity.

Given $g$ and the cascade, the Gibbs measure of the enriched Hamiltonian on the leaves and the configurations has weights proportional to $\msf v_\al w(\si)e^{\sqrt{2t}\,g\cdot v_\si+x_\al\cdot\si}$; the self-overlap corrections in \eqref{e.cascade.definition} do not affect it. Its marginal on the leaves is
\begin{equation}\label{e.cascade.gibbs}
 \tilde{\msf v}_\al:=\frac{\msf v_\al e^{U(1,x_\al;g)}}{\sum_\be\msf v_\be e^{U(1,x_\be;g)}},
\end{equation}
and, given the leaf $\al$, the configuration has the law
\begin{equation}\label{e.leaf.configuration}
 \P\bigl(\si=\tau\bigm|g,\text{weights},\text{marks},\al\bigr)
 =w(\tau)\exp\bigl(\sqrt{2t}\,g\cdot v_\tau+x_\al\cdot\tau-U(1,x_\al;g)\bigr),
 \qquad\tau\in\Omega,
\end{equation}
which is \eqref{e.sampled.configuration} with $X_1$ replaced by $x_\al$. We say that the leaf $\al$ is selected according to the Gibbs weights when it is sampled from \eqref{e.cascade.gibbs}, and we write $\P^{\mathrm{cas}}$ and $\E^{\mathrm{cas}}$ for the probability and the expectation over the weights, the marks, the selected leaf and the configuration, the realization of $g$ being fixed; this law depends on $g$ and on the step path, and we keep this dependence implicit. The law of the fields along the selected leaf under $\P^{\mathrm{cas}}$ is thus averaged over the environment (by ``environment'' we mean the weights and the marks); for a fixed realization of the environment, the selected leaf has a discrete law.

\begin{lemma}[The selected leaf of a finite-level cascade]\label{l.selected.leaf}
Let $q$ be regular, let $n\ge1$, $r_j:=j/n$, and let $\bar q^{(n)}$ be the step path with the values $q_j:=q(r_j)$ on $[r_{j-1},r_j)$, $1\le j\le n$. Consider the cascade associated with $\bar q^{(n)}$, a leaf $\al$ selected according to the Gibbs weights, and a configuration $\si$ drawn from \eqref{e.leaf.configuration}. For every realization of $g$:
\begin{enumerate}
\item $\E^{\mathrm{cas}}\log\sum_\al\msf v_\al e^{U(1,x_\al;g)}=\E\bigl[U^{(n)}(0,X_0;g)\bigm|g\bigr]$, where $X_0\sim\mcl N(0,2\diag(q(0)))$ is the initial value of the branch, independent of $g$;
\item under $\E^{\mathrm{cas}}$, the fields $(x_{\al_{|0}},\ldots,x_{\al_{|n-1}})$ along the selected leaf and the configuration~$\si$ have the law of $(X^{(n)}_{r_1},\ldots,X^{(n)}_{r_n},\si^{(n)})$, where $X^{(n)}$ is the frozen branch of Lemma~\ref{l.recursion}, and where $\si^{(n)}$ is drawn, given $X^{(n)}$, from the Gibbs measure \eqref{e.sampled.configuration} in the terminal field $X^{(n)}_1$.
\end{enumerate}
\end{lemma}

\begin{proof}
We use the invariance of Poisson--Dirichlet cascades under a tilt of the weights by a function of the marks, \cite[Lemma~3]{pantal07}, through its two consequences \cite[(3.13) and Theorem~6]{pantal07}; the first one is the recursive formula of \cite[Theorem~5.25]{HJbook}, and the formulation with a mark at the root that we use is that of \cite[Appendix~A]{chen2026spherical}. Fix $g$ and the root mark $z_\emptyset$, and regard $\Theta_\al:=U(1,x_\al;g)$ as a function of the marks $(z_{\al_{|1}},\ldots,z_{\al_{|n-1}})$ along the leaf $\al$. This function is Lipschitz, since $|\nabla U(1,\cdot\,;g)|\le\sqrt N$, so that $e^{\Theta_\al}$ is integrable, as required in \cite{pantal07}. With the parameters $m_j=r_j$, the recursion (3.4) of \cite{pantal07} sets $\Theta_{\al_{|n-1}}:=\Theta_\al$ and, for $j$ from $n-1$ down to $1$,
\begin{equation*}
 \Theta_{\al_{|j-1}}:=\frac1{r_j}\log\E_j\exp\bigl(r_j\Theta_{\al_{|j}}\bigr),
\end{equation*}
where $\E_j$ integrates the mark $z_{\al_{|j}}$ only. Since $x_{\al_{|j}}=x_{\al_{|j-1}}+\sqrt{2\diag(q(r_{j+1})-q(r_j))}\,z_{\al_{|j}}$ by \eqref{e.cascade.fields}, a backward induction based on \eqref{e.recursion} gives
\begin{equation}\label{e.recursion.dictionary}
 \Theta_{\al_{|j}}=U^{(n)}(r_{j+1},x_{\al_{|j}};g),\qquad0\le j\le n-1,
\end{equation}
and the weights $\msf v_\al e^{\Theta_\al}/\sum_\be\msf v_\be e^{\Theta_\be}$ of the tilted cascade in \cite{pantal07} are the Gibbs weights \eqref{e.cascade.gibbs}. In particular, $\Theta_{\al_{|0}}=U^{(n)}(r_1,x_\emptyset;g)$ depends on the root mark only, and \cite[(3.13)]{pantal07} gives
\begin{equation*}
 \E^{\mathrm{cas}}\Bigl[\log\sum_\al\msf v_\al e^{U(1,x_\al;g)}\Bigm|z_\emptyset\Bigr]
 =U^{(n)}(r_1,x_\emptyset;g).
\end{equation*}
Since $x_\emptyset=\sqrt{2\diag(q(r_1))}\,z_\emptyset$ has the law of $X_0+\sqrt{2\diag(q(r_1)-q(0))}\,z$, where $z$ is a standard Gaussian vector independent of $X_0$ and of $g$, averaging over $z_\emptyset$ and using the case $j=0$ of \eqref{e.recursion}, in which $T_{r_0,a}=T_{0,a}$ is a plain Gaussian expectation, gives (1). Next, \cite[Theorem~6]{pantal07} states that, conditionally on $z_\emptyset$ and under $\E^{\mathrm{cas}}$, the marks $(z_{\al_{|1}},\ldots,z_{\al_{|n-1}})$ along the selected leaf have the law under which, for $1\le j\le n-1$ and given the marks of smaller depth, $z_{\al_{|j}}$ has density $\exp(r_j[\Theta_{\al_{|j}}-\Theta_{\al_{|j-1}}])$ with respect to the standard Gaussian law. By \eqref{e.recursion.dictionary} and \eqref{e.cascade.fields}, this means that, given $x_{\al_{|j-1}}=x$, the increment $x_{\al_{|j}}-x$ has density \eqref{e.selected.transition} with respect to the centered Gaussian law with covariance $2\diag(q(r_{j+1})-q(r_j))$, while $x_{\al_{|0}}=x_\emptyset$ is a centered Gaussian vector with covariance $2\diag(q(r_1))$. Thus the fields along the selected leaf form the Markov chain of Lemma~\ref{l.recursion}, and $(x_{\al_{|0}},\ldots,x_{\al_{|n-1}})$ has the law of $(X^{(n)}_{r_1},\ldots,X^{(n)}_{r_n})$. Finally, the conditional law \eqref{e.leaf.configuration} of the configuration depends on the environment and on the leaf only through $x_\al=x_{\al_{|n-1}}$, and it is the conditional law of $\si^{(n)}$ given $X^{(n)}_1=x_\al$. This proves (2).
\end{proof}

By Lemmas~\ref{l.recursion} and~\ref{l.selected.leaf}, the frozen branch $X^{(n)}$ is, in law and after averaging over the environment of the cascade, the Gaussian field along a leaf of the cascade with $n-1$ levels selected according to the Gibbs weights, placed at the mesh points and interpolated between them according to the conditional law of the Gaussian path $Z$. Its sampled configuration $\si^{(n)}$ is the configuration drawn from the Gibbs measure of the selected leaf.

We now select two leaves. Given $g$ and the environment, let $(\al^1,\si^1)$ and $(\al^2,\si^2)$ be two independent samples from the Gibbs measure on the leaves and the configurations, so that each $\al^a$ has the law \eqref{e.cascade.gibbs}, and each $\si^a$ is drawn from \eqref{e.leaf.configuration} with $\al=\al^a$; we let $\E^{\mathrm{cas}}$ and $\P^{\mathrm{cas}}$ also integrate the second replica. Let
\begin{equation}\label{e.def.common.depth}
 J:=\max\bigl\{j\in\{0,\ldots,n-1\}:\al^1_{|j}=\al^2_{|j}\bigr\}
\end{equation}
be the depth of the most recent common ancestor of the two selected leaves, so that $J=n-1$ exactly when $\al^1=\al^2$, and the fields $x_{\al^1},x_{\al^2}$ have covariance $2\diag(\bar q^{(n)}(r_J))$. In the cascade, the two replicas thus separate at the level $r_J$. On the side of the branch, we define the Gibbs averages of the frozen problem exactly as in Section~\ref{s.gibbs}: $\si^{(n)}$ is drawn, given $X^{(n)}$, from the Gibbs measure \eqref{e.sampled.configuration} in the terminal field $X^{(n)}_1$, as in Lemma~\ref{l.selected.leaf}, and for $u\in[0,1]$ and a function $b$ of two configurations, $\langle b\rangle^{(n)}_u$ is the integral of $b$ against the twofold product of the conditional law of $\si^{(n)}$ given $\cG_u$, where $\cG_u$ is generated by $g$, $X_0$ and $(B_v)_{v\le u}$. Since the frozen branch is a Markov process, this conditional law is a function of $X^{(n)}_u$ and $g$, and $\E[\langle b\rangle^{(n)}_u\mid g]=\E[b(\si^{(n),1},\si^{(n),2})\mid g]$, where $\si^{(n),1},\si^{(n),2}$ are drawn independently from the Gibbs measures in the terminal fields of two frozen branches that coincide up to the level $u$ and are driven by independent Brownian motions afterwards. Precisely, $X^{(n),1}:=X^{(n)}$, and $X^{(n),2}$ solves the frozen equation with the same initial value $X_0$, driven by the Brownian motion that coincides with $B$ on $[0,u]$ and whose increments after $u$ are those of an independent Brownian motion; the two solutions coincide on $[0,u]$ by pathwise uniqueness. This is the description of two replicas at level~$u$ given in Section~\ref{s.gibbs.terminal}.

\begin{lemma}[Two leaves selected from a finite-level cascade]\label{l.selected.pair}
Consider the setting of Lemma~\ref{l.selected.leaf}, and let $0\le j\le n-1$. For every realization of $g$:
\begin{enumerate}
\item $\P^{\mathrm{cas}}(J=j)=1/n$;
\item conditionally on $J=j$ and under $\E^{\mathrm{cas}}$, the fields along the two selected leaves and the two configurations, that is, the two families $(x_{\al^a_{|0}},\ldots,x_{\al^a_{|n-1}},\si^a)$, $a=1,2$, have the law of the two families $(X^{(n),a}_{r_1},\ldots,X^{(n),a}_{r_n},\si^{(n),a})$, $a=1,2$, where the two frozen branches coincide up to the level $r_{j+1}$ and are driven by independent Brownian motions afterwards, and where the configurations are drawn independently in their terminal fields, as above;
\item consequently, for every $b:\Omega^2\to\R$,
\begin{equation}\label{e.selected.pair}
 \E^{\mathrm{cas}}\bigl[b(\si^1,\si^2)\1_{\{J=j\}}\bigr]
 =\frac1n\,\E\bigl[\langle b\rangle^{(n)}_{r_{j+1}}\bigm|g\bigr].
\end{equation}
\end{enumerate}
\end{lemma}


\begin{proof}
We keep the notation of the proof of Lemma~\ref{l.selected.leaf}: $g$ and the root mark $z_\emptyset$ are fixed, the recursion (3.4) of \cite{pantal07} with the parameters $m_\ell=r_\ell$ produces $\Theta_{\al_{|\ell}}=U^{(n)}(r_{\ell+1},x_{\al_{|\ell}};g)$, the weights of the tilted cascade are the Gibbs weights \eqref{e.cascade.gibbs}, and the factor $W_{\al_{|\ell}}:=\exp(r_\ell[\Theta_{\al_{|\ell}}-\Theta_{\al_{|\ell-1}}])$ of \cite[(3.5)]{pantal07} is, by \eqref{e.recursion.dictionary} and \eqref{e.cascade.fields}, the density \eqref{e.selected.transition} of the increment $x_{\al_{|\ell}}-x_{\al_{|\ell-1}}$ given $x_{\al_{|\ell-1}}$, with respect to its Gaussian law. In \cite{pantal07}, the overlap of two distinct leaves is the index of the first coordinate at which they differ, which is $J+1$ in our notation. Let $Y$ be a bounded function of the marks along two leaves. For $0\le j\le n-2$, \cite[Theorem~7]{pantal07} gives
\begin{equation}\label{e.pantal.pair}
 \E^{\mathrm{cas}}\bigl[Y\1_{\{J=j\}}\bigr]
 =(r_{j+1}-r_j)\,
 \E\Bigl[\prod_{\ell=1}^{j}W_{\al_{|\ell}}
         \prod_{\ell=j+1}^{n-1}W_{\al^1_{|\ell}}W_{\al^2_{|\ell}}\;Y\Bigr],
\end{equation}
where, on the right, $\al^1$ and $\al^2$ are two fixed leaves whose most recent common ancestor $\al_{|j}:=\al^1_{|j}=\al^2_{|j}$ has depth $j$, and the marks are independent standard Gaussian vectors, the two leaves sharing the marks of depth at most $j$. For $j=n-1$, the invariance \cite[(3.14)]{pantal07} of the pair formed by the weights and the marks, which is the content of \cite[Lemma~3]{pantal07}, gives
\begin{equation}\label{e.pantal.same.leaf}
 \E^{\mathrm{cas}}\bigl[Y\1_{\{\al^1=\al^2\}}\bigr]
 =\E\Bigl[\sum_\al\msf v_\al^2\Bigr]\,
  \E\Bigl[\prod_{\ell=1}^{n-1}W_{\al_{|\ell}}\;Y\Bigr],
\end{equation}
with a single fixed leaf $\al$ on the right, since the weights and the marks of the second process in \cite[(3.14)]{pantal07} are independent; and $\E\sum_\al\msf v_\al^2=1-r_{n-1}=1/n$, by \cite[(3.17)]{pantal07} summed over the overlap. Since $r_{j+1}-r_j=1/n$, taking $Y\equiv1$ in \eqref{e.pantal.pair} and \eqref{e.pantal.same.leaf} proves (1).

The density on the right of \eqref{e.pantal.pair}, multiplied by the Gaussian law of the marks and read through \eqref{e.cascade.fields}, is the following law of the fields along the two leaves: the common fields $(x_{\al_{|0}},\ldots,x_{\al_{|j}})$ have the law of $(X^{(n)}_{r_1},\ldots,X^{(n)}_{r_{j+1}})$, since $x_{\al_{|0}}=x_\emptyset$ has the law of $X^{(n)}_{r_1}$ once $z_\emptyset$ is averaged, and since the factors $W_{\al_{|\ell}}$, $\ell\le j$, are the transition densities \eqref{e.selected.transition} of the chain of Lemma~\ref{l.recursion}; and, given the common fields, the two families $(x_{\al^a_{|j+1}},\ldots,x_{\al^a_{|n-1}})$, $a=1,2$, are independent, each with the conditional law of $(X^{(n)}_{r_{j+2}},\ldots,X^{(n)}_{r_n})$ given $X^{(n)}_{r_{j+1}}=x_{\al_{|j}}$, by the same transition densities, each factor integrating to one in the next mark by \eqref{e.recursion}. On the side of the branch, $X^{(n),1}$ and $X^{(n),2}$ coincide up to the level $r_{j+1}$, and, given $\cG_{r_{j+1}}$, their values at the later mesh points are two independent copies of the chain started from $X^{(n)}_{r_{j+1}}$, by the Markov property of the frozen equation and the independence of the two driving noises after $r_{j+1}$. The laws of the fields in (2) thus agree. The conditional law \eqref{e.leaf.configuration} of each configuration depends on the environment and on its leaf only through the terminal field of that leaf, and the two configurations are drawn independently; the same holds for $\si^{(n),1},\si^{(n),2}$. This proves (2), and (3) follows from (1), (2) and the identity $\E[\langle b\rangle^{(n)}_{u}\mid g]=\E[b(\si^{(n),1},\si^{(n),2})\mid g]$ recalled before the statement.
\end{proof}

\subsection{Convergence to the branch}
\label{s.cascade.convergence}

We now let the number of levels tend to infinity.

\begin{proposition}[The branch as the limit of the selected leaves]\label{p.selected.limit}
Let $q$ be regular, and for $n\ge1$ let $\al$ be a leaf selected according to the Gibbs weights and $\si$ the configuration drawn at this leaf, as in Lemma~\ref{l.selected.leaf}. For every realization of $g$, the law under $\P^{\mathrm{cas}}$ of the fields along the selected leaf and of its configuration converges in total variation, as $n\to\infty$, to the conditional law given $g$ of the values of the branch $X$ at the mesh points $r_1,\ldots,r_n$ and of its sampled configuration \eqref{e.sampled.configuration}, that is,
\begin{equation}\label{e.selected.limit}
 \lim_{n\to\infty}\ \sup_{f}\,
 \Bigl|\E^{\mathrm{cas}}\bigl[f(x_{\al_{|0}},\ldots,x_{\al_{|n-1}},\si)\bigr]
       -\E\bigl[f(X_{r_1},\ldots,X_{r_n},\si)\bigm|g\bigr]\Bigr|=0,
\end{equation}
where the supremum is over all measurable $f:(\R^N)^n\times\Omega\to[-1,1]$. The statement remains true after averaging over $g$, with $f$ then also allowed to depend on $g$.
\end{proposition}

\begin{proof}
Let $X^{(n)}$ and $\si^{(n)}$ be the frozen branch of Lemma~\ref{l.recursion} and its sampled configuration, as in Lemma~\ref{l.selected.leaf}(2), so that the law of $(x_{\al_{|0}},\ldots,x_{\al_{|n-1}},\si)$ under $\P^{\mathrm{cas}}$ is the conditional law of $(X^{(n)}_{r_1},\ldots,X^{(n)}_{r_n},\si^{(n)})$ given $g$. The two laws compared in \eqref{e.selected.limit} are thus the images of the conditional laws of $(X^{(n)},\si^{(n)})$ and of $(X,\si)$ given $g$ under the map $(x,\si)\mapsto(x_{r_1},\ldots,x_{r_n},\si)$, and it suffices to show that the conditional law of $(X^{(n)},\si^{(n)})$ given $g$ converges in total variation to the conditional law of $(X,\si)$ given $g$.

The proof of Lemma~\ref{l.joint.density} uses neither the slice nor the value of the weights, and the Girsanov density there is computed at fixed $g$: for every $(\Omega,w)$ and every realization of $g$, the conditional law of $(X,\si)$ given $g$ has density
\begin{equation}\label{e.general.joint.density}
 p(g,x,\si):=w(\si)\exp\Ll(\sqrt{2t}\,g\cdot v_\si+x_1\cdot\si
                         -\int_0^1U(r,x_r;g)\,\d r\Rr)
\end{equation}
at $(x,\si)$ with respect to the product of $\ga_q$ and the counting measure on $\Omega$. Likewise, multiplying the density \eqref{e.frozen.density} of the frozen branch by the Gibbs probability \eqref{e.sampled.configuration} of $\si^{(n)}$ given $X^{(n)}_1$, which contains the factor $e^{-U(1,x_1;g)}$, we see that the conditional law of $(X^{(n)},\si^{(n)})$ given $g$ has density
\begin{equation}\label{e.frozen.joint.density}
 p_n(g,x,\si):=w(\si)\exp\Ll(\sqrt{2t}\,g\cdot v_\si+x_1\cdot\si
                         -\frac1n\sum_{j=1}^nU^{(n)}(r_j,x_{r_j};g)\Rr)
\end{equation}
with respect to the same reference measure, the term $j=n$ of the sum being $U(1,x_1;g)$. For every $g$, every continuous path $x$ and every $\si$, we have $p_n(g,x,\si)\to p(g,x,\si)$ as $n\to\infty$: by Lemma~\ref{l.frozen.convergence}, the sum in \eqref{e.frozen.joint.density} differs by at most $N\|\dot q\|_\infty/n$ from $\frac1n\sum_{j=1}^nU(r_j,x_{r_j};g)$, which is a Riemann sum of the function $r\mapsto U(r,x_r;g)$, continuous on $[0,1]$ by Lemma~\ref{l.wellposed}. For every $g$, both $p_n(g,\cdot,\cdot)$ and $p(g,\cdot,\cdot)$ are probability densities with respect to the same measure, so Scheff\'e's lemma gives
\begin{equation}\label{e.tv.fixed.g}
 \ep_n(g):=\int\sum_{\si\in\Omega}\bigl|p_n(g,x,\si)-p(g,x,\si)\bigr|\,\ga_q(\d x)
 \longrightarrow0,
\end{equation}
and $\ep_n(g)$ is twice the total variation distance between the conditional laws of $(X^{(n)},\si^{(n)})$ and $(X,\si)$ given $g$. The supremum in \eqref{e.selected.limit} is therefore at most $\ep_n(g)$, and tends to zero. Since $\ga$ is the product of the law of $g$ and of $\ga_q$, the total variation distance between the laws of $(g,X^{(n)},\si^{(n)})$ and $(g,X,\si)$ is $\frac12\E\,\ep_n(g)$, which tends to $0$ by dominated convergence, since $\ep_n\le2$; this gives the statement averaged over $g$.
\end{proof}


We now turn to the two replicas of Lemma~\ref{l.selected.pair}, and let the number of levels tend to infinity. The conditional law of the two configurations given the depth of their common ancestor converges, uniformly in the depth, to the law of two replicas at the corresponding level of the branch.

\begin{proposition}[Two replicas at level $r$ as conditioned cascade replicas]\label{p.cascade.pair}
Let $q$ be regular, let $b:\Omega^2\to\R$, and for $n\ge1$ let $(\al^1,\si^1)$, $(\al^2,\si^2)$ and $J$ be as in Lemma~\ref{l.selected.pair}. For every realization of $g$,
\begin{equation}\label{e.cascade.pair.uniform}
 \lim_{n\to\infty}\ \max_{0\le j\le n-1}\ \sup_{u\in[r_j,r_{j+1}]}\,
 \Bigl|\E^{\mathrm{cas}}\bigl[b(\si^1,\si^2)\bigm|J=j\bigr]
       -\E\bigl[\langle b\rangle_u\bigm|g\bigr]\Bigr|=0,
\end{equation}
and consequently, for every continuous $\varphi:[0,1]\to\R$,
\begin{equation}\label{e.cascade.pair.integrated}
 \lim_{n\to\infty}\E^{\mathrm{cas}}\bigl[\varphi(r_J)\,b(\si^1,\si^2)\bigr]
 =\int_0^1\varphi(r)\,\E\bigl[\langle b\rangle_r\bigm|g\bigr]\,\d r .
\end{equation}
Both statements remain true after averaging over $g$; since $\P^{\mathrm{cas}}(J=j)=1/n$ does not depend on $g$, the averaged conditional expectation in \eqref{e.cascade.pair.uniform} is the conditional expectation given $J=j$ under the law that also averages over $g$.
\end{proposition}

\begin{proof}
Fix a realization of $g$, and let $\de_n(g)$ denote the supremum in \eqref{e.selected.limit}, so that $\de_n(g)\to0$ by Proposition~\ref{p.selected.limit}. We proceed in three steps.

\emph{Step 1: stability of two replicas.} For a signed measure $\rho$, let $\|\rho\|$ denote the total mass of $|\rho|$, so that $|\int b\,\d\rho|\le\|b\|_\infty\|\rho\|$. Let $\mu,\mu'$ be probability measures on $\R^N$, and let $\msf k,\msf k'$ be probability kernels from $\R^N$ to $\Omega$. We claim that
\begin{equation}\label{e.kernel.stability}
 \bigl\|\mu(\d x)\,\msf k(x,\cdot)^{\otimes2}-\mu'(\d x)\,\msf k'(x,\cdot)^{\otimes2}\bigr\|
 \le5\,\bigl\|\mu(\d x)\,\msf k(x,\cdot)-\mu'(\d x)\,\msf k'(x,\cdot)\bigr\| .
\end{equation}
Indeed, for probability measures $\rho,\rho'$ on $\Omega$, we have $\|\rho^{\otimes2}-\rho'^{\otimes2}\|\le\|(\rho-\rho')\otimes\rho\|+\|\rho'\otimes(\rho-\rho')\|=2\|\rho-\rho'\|$. Hence the left side of \eqref{e.kernel.stability} is at most
\begin{equation*}
 \|\mu-\mu'\|+\int\mu'(\d x)\,\bigl\|\msf k(x,\cdot)^{\otimes2}-\msf k'(x,\cdot)^{\otimes2}\bigr\|
 \le\|\mu-\mu'\|+2\,\bigl\|\mu'(\d x)\,\msf k(x,\cdot)-\mu'(\d x)\,\msf k'(x,\cdot)\bigr\|,
\end{equation*}
and the last norm is at most $\|\mu'-\mu\|$ plus the right side of \eqref{e.kernel.stability} without the factor $5$. Since $\|\mu-\mu'\|$ is the norm of the difference of the first marginals, it is also at most this right side, and \eqref{e.kernel.stability} follows.

\emph{Step 2: uniform convergence of the two-replica averages.} For $u\in[0,1]$, let $\mu_u$ be the conditional law of $X_u$ given $g$. By the Markov property of the branch, the conditional law of $\si$ given $\cG_u$ is a function of $X_u$: there is a probability kernel $\msf k_u$ from $\R^N$ to $\Omega$, namely $\msf k_u(x,\tau)$ is the probability that the configuration sampled in the terminal field of the branch started from $(u,x)$ equals $\tau$, such that $\P(\si=\tau\mid\cG_u)=\msf k_u(X_u,\tau)$ for every $\tau\in\Omega$. By \eqref{e.response.first} applied to $b=\1_{\{\si=\tau\}}$ and the proof of Lemma~\ref{l.response}, $\msf k_u(x,\tau)$ is the value at $(u,x)$ of the solution of a linear parabolic equation, and $(u,x)\mapsto\msf k_u(x,\tau)$ is continuous and bounded on $[0,1]\times\R^N$. Hence the conditional law of $(X_u,\si)$ given $g$ is $\mu_u(\d x)\,\msf k_u(x,\cdot)$, and, by the definition of $\langle\cdot\rangle_u$,
\begin{equation*}
 \E\bigl[\langle b\rangle_u\bigm|g\bigr]
 =\int\mu_u(\d x)\sum_{\tau,\tau'\in\Omega}b(\tau,\tau')\,\msf k_u(x,\tau)\,\msf k_u(x,\tau').
\end{equation*}
The same holds for the frozen branch, with the conditional law $\mu^{(n)}_u$ of $X^{(n)}_u$ given $g$ and the kernel $\msf k^{(n)}_u$ defined with the frozen equation. At a mesh point $u=r_{j+1}$, $0\le j\le n-1$, the conditional law of $(X^{(n)}_u,\si^{(n)})$ given $g$ is the law of $(x_{\al_{|j}},\si)$ under $\P^{\mathrm{cas}}$, by Lemma~\ref{l.selected.leaf}(2), so the conditional laws of $(X^{(n)}_u,\si^{(n)})$ and $(X_u,\si)$ given $g$ are the images, under the projection onto the $(j+1)$-th field and the configuration, of the two laws compared in \eqref{e.selected.limit}. The norm of their difference is thus at most $\de_n(g)$, and \eqref{e.kernel.stability} gives
\begin{equation}\label{e.pair.uniform.n}
 \max_{1\le j\le n}\Bigl|\E\bigl[\langle b\rangle^{(n)}_{r_j}\bigm|g\bigr]-\E\bigl[\langle b\rangle_{r_j}\bigm|g\bigr]\Bigr|
 \le5\,\|b\|_\infty\,\de_n(g).
\end{equation}

\emph{Step 3: conclusion.} Since $(u,x)\mapsto\msf k_u(x,\tau)$ is continuous and bounded for every $\tau$ and $u\mapsto X_u$ is continuous, the function $u\mapsto\E[\langle b\rangle_u\mid g]=\E[\sum_{\tau,\tau'}b(\tau,\tau')\msf k_u(X_u,\tau)\msf k_u(X_u,\tau')\mid g]$ is continuous on $[0,1]$, by dominated convergence. Let $\omega_g$ denote its modulus of continuity. By Lemma~\ref{l.selected.pair}, $\E^{\mathrm{cas}}[b(\si^1,\si^2)\mid J=j]=\E[\langle b\rangle^{(n)}_{r_{j+1}}\mid g]$, so that, for $u\in[r_j,r_{j+1}]$, \eqref{e.pair.uniform.n} at the mesh point $r_{j+1}$ gives
\begin{equation*}
 \Bigl|\E^{\mathrm{cas}}\bigl[b(\si^1,\si^2)\bigm|J=j\bigr]-\E\bigl[\langle b\rangle_u\bigm|g\bigr]\Bigr|
 \le5\,\|b\|_\infty\,\de_n(g)+\omega_g(1/n),
\end{equation*}
which proves \eqref{e.cascade.pair.uniform}. Next, since $\P^{\mathrm{cas}}(J=j)=1/n$, we have $\E^{\mathrm{cas}}[\varphi(r_J)\,b(\si^1,\si^2)]=\frac1n\sum_{j=0}^{n-1}\varphi(r_j)\,\E^{\mathrm{cas}}[b(\si^1,\si^2)\mid J=j]$, so the previous display gives
\begin{equation*}
 \Bigl|\E^{\mathrm{cas}}\bigl[\varphi(r_J)\,b(\si^1,\si^2)\bigr]
 -\frac1n\sum_{j=0}^{n-1}\varphi(r_j)\,\E\bigl[\langle b\rangle_{r_j}\bigm|g\bigr]\Bigr|
 \le\|\varphi\|_\infty\bigl(5\|b\|_\infty\de_n(g)+\omega_g(1/n)\bigr),
\end{equation*}
and the Riemann sum converges to the right side of \eqref{e.cascade.pair.integrated}, the integrand being continuous. Finally, all the quantities in \eqref{e.cascade.pair.uniform} and \eqref{e.cascade.pair.integrated} are bounded by $\|\varphi\|_\infty\|b\|_\infty$ (with $\varphi\equiv1$ in the first case), so the statements averaged over $g$ follow by dominated convergence.
\end{proof}


\begin{remark}[Disintegration of the derivative formulas]\label{r.pair.derivatives}
For the cascade with $n-1$ levels at the step path $\bar q^{(n)}$, a Gaussian integration by parts gives, as for \eqref{e.cascade.lipschitz}, $\dr_t\bar F^{\mathrm{cas}}_{N,\Omega}(t,\bar q^{(n)})=\E\,\E^{\mathrm{cas}}[\xi(R_N(\si^1,\si^2))]$: the derivatives of the cascade free energy are expectations of two-replica observables under the Gibbs measure, with no conditioning on the separation level. Taking $b=\xi(R_N)$ and $\varphi\equiv1$ in \eqref{e.cascade.pair.integrated}, we see that the limit of this expression is $\int_0^1\E\langle\xi(R_N)\rangle_r\,\d r=\dr_t\bar F_{N,\Omega}(t,q)$. Similarly, taking $b=R_{N,s}$ in \eqref{e.cascade.pair.uniform}, we obtain, for every $r\in[0,1]$ and every species $s$,
\begin{equation}\label{e.p.conditioned}
 p_{N,s}(r)=\E\langle R_{N,s}\rangle_r
 =\lim_{n\to\infty}\E\,\E^{\mathrm{cas}}\bigl[R_{N,s}(\si^1,\si^2)\bigm|J=\min(\lfloor nr\rfloor,n-1)\bigr].
\end{equation}
In words, the quantity $p_N(r)$ is the limit of the mean overlap of two cascade replicas conditioned on separating at level $r$. Theorem~\ref{t.derivatives} thus disintegrates the derivative formulas of the cascade, in which the separation level is uniformly distributed, according to the value of this level.
\end{remark}

\subsection{Agreement of the free energies}
\label{s.cascade.agreement}

\begin{proposition}[Agreement with the cascade]\label{p.cascade}
For every $t\ge0$ and $q\in\cQ_1^{\sS}$, we have
\begin{equation*}
 \bar F_{N,\Omega}(t,q)=\bar F^{\mathrm{cas}}_{N,\Omega}(t,q).
\end{equation*}
\end{proposition}

\begin{proof}
Let $q$ be regular, and let $\bar q^{(n)}$ be as in Lemma~\ref{l.selected.leaf}. Averaging Lemma~\ref{l.selected.leaf}(1) over~$g$, and since $q_n=q(1)$, we obtain
\begin{equation*}
 \bar F^{\mathrm{cas}}_{N,\Omega}(t,\bar q^{(n)})
 =t\xi(\lambda_N)+\lambda_N\cdot q(1)-\frac1N\E\,U^{(n)}(0,X_0;g).
\end{equation*}
As $n\to\infty$, the right side converges to $\bar F_{N,\Omega}(t,q)$ by Lemma~\ref{l.frozen.convergence} and \eqref{e.free.energy.V}. On the left side, $0\le\bar q^{(n)}-q\le q(r_{j+1})-q(r_j)$ on $[r_j,r_{j+1})$, so that $\|\bar q^{(n)}_s-q_s\|_{L^1}\le n^{-1}(q_s(1)-q_s(0))$ for every species $s$, and the left side converges to $\bar F^{\mathrm{cas}}_{N,\Omega}(t,q)$ by \eqref{e.cascade.lipschitz}. This proves the proposition for regular paths. The general case follows from the density of regular paths in $\cQ_1^{\sS}$, by Proposition~\ref{p.extension} and \eqref{e.cascade.lipschitz}. The identity concerns the corrected free energies as a whole: the terminal value $q(1)$ alone is not continuous in $L^1$, and it may be infinite for a general path in $\cQ_1^{\sS}$.
\end{proof}

\begin{remark}[Cost and entropy]\label{r.entropy}
For the finite-level cascade of Lemma~\ref{l.selected.leaf}, the relative entropy of the Gibbs weights \eqref{e.cascade.gibbs} with respect to the cascade weights is
\begin{equation}\label{e.gibbs.variational}
 \Ent(\tilde{\msf v}\mid\msf v):=\sum_\al\tilde{\msf v}_\al\log\frac{\tilde{\msf v}_\al}{\msf v_\al}
 =\sum_\al\tilde{\msf v}_\al U(1,x_\al;g)-\log\sum_\al\msf v_\al e^{U(1,x_\al;g)}.
\end{equation}
Under $\E^{\mathrm{cas}}$, the first term on the right of \eqref{e.gibbs.variational} has mean $\E[U(1,X^{(n)}_1;g)\mid g]$ by Lemma~\ref{l.selected.leaf}(2), and the second has mean $\E[U^{(n)}(0,X_0;g)\mid g]$ by Lemma~\ref{l.selected.leaf}(1). Applying It\^o's formula to the frozen problem as in the proof of Lemma~\ref{l.ito}, we obtain, for every realization of $g$,
\begin{equation}\label{e.leaf.entropy}
\begin{split}
 \E^{\mathrm{cas}}\Ent(\tilde{\msf v}\mid\msf v)
 &=\E\bigl[U(1,X^{(n)}_1;g)-U^{(n)}(0,X_0;g)\bigm|g\bigr]\\
 &=\E\Ll[\int_0^1\rho_n(r)\,\nabla U^{(n)}\cdot D(r)\nabla U^{(n)}(r,X^{(n)}_r;g)\,\d r\Bigm|g\Rr],
\end{split}
\end{equation}
all terms being integrable, since $U(1,\cdot\,;g)$ has linear growth and $X^{(n)}_1$ has Gaussian tails. The right side of \eqref{e.leaf.entropy} is the cost $A_q$ of Definition~\ref{d.free.energy} for the frozen branch. As $n\to\infty$, it converges to $\E[A_q\mid g]=\E[U(1,X_1;g)-U(0,X_0;g)\mid g]$, by \eqref{e.ito.identity}: the second term in \eqref{e.leaf.entropy} converges by Lemma~\ref{l.frozen.convergence}, and the first by Lemma~\ref{l.selected.leaf}(2) and Proposition~\ref{p.selected.limit}, since $|U(1,x;g)|\le|U(1,0;g)|+\sqrt N|x|$ and the second moments of $X^{(n)}_1$ are bounded uniformly in $n$, the drift of the frozen branch being bounded by $2\sqrt N\|\dot q\|_\infty$. Thus the cost is, in expectation, the limit of the relative entropy of the Gibbs weights of the leaves with respect to the cascade weights.
\end{remark}
%
%
%
%
\bibliographystyle{plain}
\bibliography{biased}
\end{document}